\documentclass[final,numrefs,twoadvisors,noinfo]{nddiss2e}

\usepackage{enumitem,linegoal}
\usepackage{lipsum}
\usepackage{geometry}
\usepackage{caption}
\usepackage{proof}
\usepackage{tikz}
\usetikzlibrary{arrows.meta}
\usepackage{caption}
\usepackage{amsmath, amsthm, amssymb}
\usepackage{mathrsfs}
\usepackage{graphicx}
\usepackage{amsbsy}
\usepackage{thmtools}
\usepackage{footmisc}
\usepackage{mathptmx}
\usepackage{mathtools}
\usepackage{natbib}
\usepackage{graphicx}
\usepackage{stmaryrd}
\usepackage{calrsfs}
\usepackage{marginnote}
\usepackage{hyperref}

 \allowdisplaybreaks 
 \usepackage[verbose]{placeins}

\DeclareMathAlphabet{\pazocal}{OMS}{zplm}{m}{n}
\SetMathAlphabet\pazocal{bold}{OMS}{zplm}{bx}{n}

\theoremstyle{definition}

\newtheorem{theorem}{Theorem}
\newtheorem{definition}[theorem]{Definition}
\newtheorem{corollary}{Corollary}[theorem]
\newtheorem{prop}[theorem]{Proposition}
\newtheorem{lemma}[theorem]{Lemma}

\newtheorem{question}[theorem]{Open Question}
\newtheorem{remark}[theorem]{Remark}
\newtheorem{claim}{Claim}[theorem]
\newtheorem{example}[theorem]{Example}

\newenvironment{claimproof}[1]{\par\noindent\textit{Proof of the Claim.}\space#1}{\hfill $\blacksquare$}

\newcommand\A{\mathcal{A}}
\newcommand\G{\mathcal{G}}
\newcommand\F{\mathcal{F}}
\newcommand\I{\mathcal{I}}
\newcommand\J{\mathcal{J}}
\newcommand\U{\mathcal{U}}
\newcommand\T{\mathcal{T}}
\newcommand\W{\mathcal{W}}
\newcommand\M{\mathcal{M}}
\newcommand\ZFUR{\textup{ZFU}_\text{R}}
\newcommand\ZFCUR{\textup{ZFCU}_\text{R}}
\newcommand\GBUR{\textup{GBcU}_\text{R}}
\newcommand\KMUR{\textup{KMcU}_\text{R}}

\newcommand\ACA{\textup{AC}^\mathcal{A}}
\renewcommand{\P}{{\mathbb P}}
\newcommand\Q{\mathbb{Q}}
\newcommand\B{\mathbb{B}}

\newcommand{\barx}{\bar{x}}
\newcommand{\bary}{\bar{y}}
\newcommand{\barz}{\bar{z}}
\newcommand{\barw}{\bar{w}}

\newcommand{\dotx}{\dot{x}}
\newcommand{\doty}{\dot{y}}
\newcommand{\dotz}{\dot{z}}
\newcommand{\dotu}{\dot{u}}

\newcommand{\dotw}{\dot{w}}
\newcommand{\xcheck}{\check{x}}

\newcommand{\forces}{\Vdash}
\newcommand{\nforces}{\nVdash}
\newcommand{\DCK}{\textup{DC}_\kappa\textup{-scheme}}

\newcommand\UB{U^{\mathbb{B}}}
\newcommand\VB{V^{\mathbb{B}}}
\newcommand\UBB{U^{\mathbb{B}}}
\newcommand\BB{\mathbb{B}}

\newcommand\Ucheck{\check{U}}

\def\<#1>{\left\langle#1\right\rangle}
\def\[#1]{\left\llbracket#1\right\rrbracket}
\renewcommand{\restriction}{\mathord{\upharpoonright}}
\newcommand{\Aeq}{\overset{\mathscr{A}}{=}}

\newgeometry{left=1.5in, top=1in, right=1in, bottom=1in}

\begin{document}

\frontmatter             
\title{SET THEORY WITH URELEMENTS}  
\author{Bokai Yao}      
\work{Dissertation}    

\degaward{Doctor of Philosophy}
\advisor{Joel David Hamkins}
\vspace{0pt} 
\secondadvisor{Daniel Nolan} 
\department{Philosophy}        

\maketitle               

\makecopyright


\begin{abstract}
This dissertation aims to provide a comprehensive account of set theory with urelements. In Chapter 1, I present mathematical and philosophical motivations for studying urelement set theory and lay out the necessary technical preliminaries. Chapter 2 is devoted to the axiomatization of urelement set theory, where I introduce a hierarchy of axioms and discuss how ZFC with urelements should be axiomatized. The breakdown of this hierarchy of axioms in the absence of the Axiom of Choice is also explored. In Chapter 3, I investigate forcing with urelements and develop a new approach that addresses a drawback of the existing machinery. I demonstrate that forcing can preserve, destroy, and recover the axioms isolated in Chapter 2 and discuss how Boolean ultrapowers can be applied in urelement set theory. Chapter 4 delves into class theory with urelements. I first discuss the issue of axiomatizing urelement class theory and then explore the second-order reflection principle with urelements. In particular, assuming large cardinals, I construct a model of second-order reflection where the principle of limitation of size fails. \enlargethispage{\baselineskip}
\end{abstract}

\tableofcontents
\listoffigures
\begin{acknowledge}
I thank Joel, my academic idol, for his exceptional supervision of this dissertation. I feel incredibly fortunate to have had the opportunity to learn from him. Our discussions have been a source of pure intellectual enjoyment.

I thank Daniel for his invaluable co-supervision and for providing me with extensive help during my studies at Notre Dame; Tim for teaching me the first course on set theory, which sparked my passion; Paddy for her insightful feedback and encouragement throughout the writing process; Asaf Karagila for his valuable discussions via email.

I would like to acknowledge my mother for her unwavering support since the day I changed my major in college. I am grateful to Xinhe for being an amazing partner in both academia and life, and for the love and joy she brought into my life. I thank my cat Jianjian, the cutest urelement, for allowing me to use his chair.

 \end{acknowledge}




\mainmatter

\chapter{Introduction}
Section \ref{section:UrelementsinSetTheory} presents mathematical and philosophical motivations for studying set theory with urelements. Section \ref{section:basicaxioms} introduces the basic axioms of set theory with urelements. Section \ref{section:interpretingUinV} reviews a well-known method of interpreting urelements in pure set theory, through which various versions of urelement set theory can be interpreted in pure set theory.

\section{Urelements in set theory}\label{section:UrelementsinSetTheory}
 Urelements are members of sets that are not themselves sets. While urelements, as non-reducible mathematical objects, were included in the earlier development of set theory (e.g., in Zermelo's original presentation \cite{zermelo1930grenzzahlen}), most contemporary set theorists decided that their role is superfluous. This is because within a reasonable pure set theory such as ZFC, all mathematical objects can be recovered \textit{up to isomorphism}. However, set theory with urelements are still of mathematical and philosophical interests.

\subsection{The mathematics}
In the pre-forcing era, permutation models were used to establish the independence of the Axiom of Choice in the Zermelo-Fraenkel set theory (ZF) with urelements. This technique was developed by Fraenkel, Mostowski, and Specker \cite{fraenkel1922begriff,Mostowski1939-MOSBDB,Specker1957-SPEZAD}. With the invention of forcing, permutation models became a flexible method of obtaining independence results concerning choice principles, thanks to the Jech-Sochor Embedding Theorem. Urelements play an essential role in alternative set theories such as Kripke–Platek set theory \cite{barwise2017admissible}, Quine's New Foundations \cite{jensen1968consistency}, versions of constructive set theory \cite{Cantini2008-CROCST-4}, and non-well founded set theory \cite{barwise1996vicious}.

However, many questions regarding urelement set theory remain unexplored. Most existing studies of ZF with urelements, such as \cite{jech2008axiom} and \cite{potter2004set}, assume as an axiom that the urelements form a set.\footnote{In \cite{barwise2017admissible}, Barwise studies Kripke–Platek set theory that allows a proper class of urelements.} This assumption is highly unnatural and raises several issues. Firstly, whether the urelements form a set should not be settled by an axiom of set theory. According to many philosophical arguments, certain abstract entities, such as propositions and possible worlds, cannot form a set by their nature.\footnote{For arguments for no set of all propositions, see \cite{Grim1984-GRITIN-5} and \cite{Menzel2012-MENSAW}. The principle that for every cardinal $\kappa$, there is a set of urelements of size $\kappa$ also appears in the discussion of \textit{recombination principles} in modal metaphysics. See \cite{Forrest1984-FORAAA} and \cite{nolan1996recombination}.} Secondly, permitting only a set of urelements conceals a great deal of subtleties of urelement set theory, while allowing a proper class of urelements is needed to understand many set-theoretic axioms and constructions fully. Surprisingly, there does not seem to be a systematic study of ZF(C) with a class of urelements, and many fundamental questions are thus left unanswered. This dissertation aims to address some of these questions by focusing on the following three:
\begin{enumerate}
\item What is the most general formalization of set theory (and class theory) with urelements?
\item How do standard set-theoretic constructions, such as forcing, behave in the context of urelements?
\item Can the existence of urelements affect the strength of strong axioms of infinity?
\end{enumerate}

\subsection{The philosophy}\label{subsection:philosophy}
One philosophical motivation for studying urelement set theory comes from the potential applications of set theory to domains containing non-sets. For example, since the mass-function maps the physical objects to the real numbers, they can only be reduced to sets if set theory allows urelements. Similarly, philosophers often give arguments involving \textit{sets} of propositions, mereological fusions, and possible worlds. Set theory with urelements thus provides a foundation for studying these debates.

Urelements can also have implications on the philosophical conceptions of set. The Naive Comprehension Principle, which says that for every condition $\varphi$ there is a set of things satisfying $\varphi$, seems to be the most natural conception of set. By Russell's paradox, we know that this naive conception is simply contradictory because the Russell set $\{x: x \notin x\}$ cannot exist. If the commonly accepted axioms of ZF can be seen as a response to Russell's paradox, they must be justified on the basis of different conceptions. Three conceptions of set are often discussed in the literature (\cite{maddy1988believing}): \textit{the iterative conception}, \textit{limitation of size}, and \textit{reflection}. 

According to the iterative conception of set, sets are formed in \textit{stages}: at each stage sets of things available on the earlier stages are formed. The Russell set simply cannot exist in this picture because it cannot be formed at any stage. Boolos \cite{boolos1971iterative} shows that a faithful formulation of the iterative conception recovers almost all the axioms of Zermelo set theory. A more general version of the iterative conception also asserts that there are as many stages as possible (see \cite{Martin2001-MARMUO}, which can be adopted to justify the Axiom of Replacement. 

According to the limitation-of-size conception (\cite{hallett1986cantorian}), a collection of objects form a set as long as they are not ``too many''. This is often seen as ``one step back from disaster'' since one might think that the Russell set is paradoxical percisely because it contains too many objects. When formulated in class theory, limitation of size provides justifications for the Axiom of Replacement and the Axiom of Choice. 

The reflection conception (\cite{tait2005constructing}) aims to articulate the idea that the universe of sets is so enormous that it is \textit{indescribable}, i.e., any true statement about the universe of sets is already true in an initial segment of the universe. Reflection principles are often taken as a form of \textit{maximality principle}, and when formulated as an axiom scheme, they provide justifications for many axioms of ZF including the Axiom of Inifinity and Replacement. In the second-order context, reflection principles also produce large cardinals (\cite{bernays1976problem}).

The three conceptions are by no means in conflict, and people often appeal to them simultaneously when justifying the axioms of ZFC. But how do they relate to each other, and which one, if any, is a more robust conception? In the context of pure set theory without urelements, these questions seem to have easy answers. For example, limitation of size is also viewed as a robust maximality principle (e.g., G\"odel \cite{godel1986collected} holds this view) because it allows sets to form as much as possible. Furthermore, since first-order reflection follows as a theorem scheme in ZF, the reflection conception appears to be a consequence of the other two. However, when urelements are included, the situation becomes completely different. In Chapter 4, we will see that limitation of size, when formulated as an axiom in class theory with urelements, becomes highly unnatural. Instead, it is its negation that may count as a maximality principle. Additionally, with urelements, the first-order reflection principle no longer follows from the basic axioms of ZF. In fact, in Chapter 2, I will present evidence suggesting that reflection is a more robust conception of set than the other two.

\section{Basic axioms}\label{section:basicaxioms}
The language of urelement set theory, in addition to $\in$, contains a unary predicate $\A$ for urelements. $Set(x)$ abbreviates $\neg\A(x)$. The standard axioms (and axiom schemes) of ZFC, modified to allow urelements, are as follows.

\begin{itemize}
\item[] (Axiom $\A$) $\forall x (\A(x) \rightarrow \neg \exists y (y \in x))$.
\item[] (Extensionality) $\forall x, y (Set(x) \land Set(y) \land \forall z (z \in y \leftrightarrow z \in x) \rightarrow x = y)$
\item[] (Foundation) $\forall x (\exists y (y \in x) \rightarrow \exists z\in x \ (z \cap x = \emptyset))$
\item[] (Pairing) $\forall x, y \exists z \forall v (v \in z \leftrightarrow v = x \lor v = y )$
\item[] (Union) $\forall x \exists y \forall z (z \in y \leftrightarrow \exists w \in x \ ( z \in w))$.
\item[] (Powerset) $\forall x \exists y \forall z (z \in y \leftrightarrow Set(z) \land z \subseteq x)$
\item [] (Separation) $\forall x, u \exists y \forall z (z \in y \leftrightarrow z \in x \land \varphi(z, u))$
\item[] (Infinity) $\exists s (\exists y \in s \ (Set(y) \land \forall z (z \notin y)) \land \forall x \in s \ (x \cup \{x\} \in s))$
\item[] (Replacement) $\forall w, u (\forall x \in w \ \exists ! y \varphi(x, y, u)   \rightarrow \exists v \forall x \in w \ \exists y \in v \ \varphi(x, y, u))$
\item[] (Collection) $\forall w, u (\forall x \in w \ \exists y \varphi(x, y, u)   \rightarrow \exists v \forall x \in w\  \exists y \in v\  \varphi(x, y, u))$.
\item[] (AC) Every set is well-orderable.\footnote{Over $\ZFUR$, this form of AC is equivalent to other variants such as the principle that every family of non-empty sets has a choice function and Zorn's Lemma. The proofs are all standard and hence omitted.}
\end{itemize}
\begin{definition}
\ \newline ZU = Axiom $\A$ + Extensionality + Foundation + Pairing + Union + Powerset + Infinity + Separation.\\
$\ZFUR = $ ZU + Replacement. \\
$\ZFCUR = $ $\ZFUR$ + AC.\\
ZF = $\ZFUR$ + $\forall x \neg \A (x)$.\\
ZFC =  ZF + AC.
\end{definition}
\noindent Note that there is no axiom asserting that the urelements form a set. The subscript R indicates that the correponding theories are only formulated with Replacement rather than Collection. In fact, it is folklore that $\ZFCUR$ cannot prove Collection when a proper class of urelements is allowed. In the next chapter, the issue of axiomatizing ZFC with urelements will be discussed in depth, and we shall see that it is $\ZFCUR$ + Collection that is a more robust urelement set theory.

\begin{definition}\label{def:kernel,cardinal}
    A set is \textit{transitive} if every member of it is a subset of it. The transitive closure of a set $x$, $trc(x)$, is the smallest transitive set that contains $x$. The \textit{kernel} of an object $x$, $ker(x)$, is the set of the urelements in $trc(\{x\})$. A set is \textit{pure} if its kernel is empty. A set $\alpha$ is an ordinal if it is a transitive pure set well-ordered by $\in$. The powerset of a set $x$, $P(x)$, is the set of all \textit{subsets} of $x$. $x \sim y$ abbreviates ``$x$ is equinumerous with $y$''. A set $\kappa$ is a cardinal if it is an ordinal that is not equinumerous with any ordinal below $\kappa$.
\end{definition}
\noindent Under this definition, the kernel of a urelement is its singleton, which is somewhat nonstandard but will be useful for our purpose. In $\ZFUR$, we shall also (informally) talk about first-order parametrically definable classes.
\begin{definition}
  $U$ is the class of all objects; $V$ is the class of all pure sets; $Ord$ is the class of all ordinals; and $\A$ also denotes the class of all urelements. $A \subseteq \A$ thus means ``$A$ is a set of urelements''.  
\end{definition}
 It is routine to check that $\ZFUR$ proves the transfinite recursion theorem (it is understood that all classes are parametrically definable classes).
\begin{theorem}[$\ZFUR$]
Let $R$ be a well-founded and set-like class relation on a class $X$. Given a class function $F: X \times U \rightarrow U$, there is a unique class function $G: X \rightarrow U$ such that for every $ x \in X $, $G(x) = F(x, G\restriction \{y \in X : \<y, x> \in R\})$.  \qed
\end{theorem}
In $\ZFUR$, every set $x$ has a transitive closure and hence a kernel. For any set of urelements $A$, by transfinite recursion on $Ord$ we define the $V_{\alpha}(A)$-hierarchy as follows.
\begin{itemize}
    \item [] $V_0(A) = A$;
    \item [] $V_{\alpha+1}(A) = P(V_{\alpha}(A)) \cup V_{\alpha}(A)$;
    \item [] $V_{\gamma}(A) = \bigcup_{\alpha < \gamma} V_\alpha(A)$, where $\gamma$ is a limit;
    \item [] $V(A) = \bigcup_{\alpha \in Ord} V_\alpha(A)$.
\end{itemize}
\noindent For every $x$ and set $A$ of urelements, $x \in V(A)$ if and only if $ker(x) \subseteq A$. Every permutation $\pi$ of a set of urelements $A$ can be extended to a definable permutation of the universe $U$ in a canonical way: we let $\pi a = a$ for every urelement $a \notin A$ and $\pi x = \{\pi y : y \in x \}$ for every set $x$ by transfinite recursion. Such $\pi$ preserves $\in$ and is thus an automorphism of the universe $U$. For every set $x$, whenever $\pi$ \textit{point-wise fixes} $ker(x)$, i.e., $\pi a = a$ for every $a \in ker(x)$, $\pi$ also point-wise fixes $x$. Finally, it is a useful fact that $\ZFUR$ proves the following restricted version of Collection.
\begin{itemize}
   \item[] (Collection$^-$) $\forall w, u (\exists A \subseteq \A \ \forall x \in w \ \exists y \in V(A)\ \varphi(x, y, u)   \rightarrow \exists v \forall x \in w\  \exists y \in v\  \varphi(x, y, u))$.
\end{itemize}
\begin{prop}\label{weakcollection}
$\ZFUR \vdash$ Collection$^-$.
\end{prop}
\begin{proof}
For every $x \in w$, let $\alpha_x$ be the least $\alpha$ such that there is some $y \in V_\alpha (A)$ with $\varphi(x, y, u)$ and let $\alpha = \bigcup_{x \in w} \alpha_x$. $V_\alpha (A)$ is then the desired collection set $v$.
\end{proof}

\section{Interpreting $U$ in $V$}\label{section:interpretingUinV}
There is a canonical way of interpreting urelement set theory in pure set theory, which seems to appear first in the appendix of \cite{barwise2017admissible}.

\begin{definition}\label{barwiseinterpretation1}
Let $V$ be a model of ZF and $X$ be a class of $V$. In $V$, we define $V \llbracket X \rrbracket$ by recursion as follows. 
\begin{itemize}
    \item [] $ V \llbracket X \rrbracket = (\{0\} \times X) \cup \{\barx \in V : \exists x (\barx = \<1 , x> \land x \subseteq V \llbracket X \rrbracket)\}.$
\end{itemize}
For every $\barx, \bary \in V \llbracket X \rrbracket$,
\begin{itemize}
    \item [] $\barx \ \bar{\in} \ \bary \text{ if and only if } \exists y (\bary = \<1, y> \land  \barx  \in y);$
    \item [] $\bar{\A}(\barx) \text{ if and only if } \barx \in \{0\} \times X.$
\end{itemize}
$V \llbracket X \rrbracket$ will also denote the model $\<V \llbracket X \rrbracket,\ \bar{\A},\ \bar{\in}>$ for the language of urelement set theory. For every $x \in V$, define $\hat{x} = \<1, \{\hat{y} : y \in x\}>$ and let $\hat{V} = \{\hat{x} : x \in V\}$. For any $\barx. \bary \in V \llbracket X \rrbracket$, let $\overline{\{\barx, \bary\}} = \<1, \{\barx, \bary \}>$, which codes the pair of $\barx$ and $ \bary$ in $V\llbracket X \rrbracket$.
\end{definition}
\noindent That is, we treat $\{0\} \times X$ as the class of urelements and then generate sets by closing under the operation: if $x \subseteq V \llbracket X \rrbracket$, then $\<1, x> \in V \llbracket X \rrbracket$. By an easy induction, one can show that the map $x \mapsto \hat{x}$ is an isomorphism from $\<V, \in>$ to $\<\hat{V}, \bar{\in}>$.

\begin{theorem}\label{thm:V[X]modelsZFU}
Let $V$ be a model of ZF and $X$ be a class of $V$. $V \llbracket X \rrbracket \models$ $\ZFUR$ + Collection.
\end{theorem}
\begin{proof}
For Extensionality, if $\bar{x}, \bar{y}$ are sets in $V \llbracket X \rrbracket$ with the same $\bar{\in}$-members, then $x = \<1, \bar{x}> = \<1, \bar{y}> = y$. Foundation holds because $\bar{\in}$ is a well-founded relation in $V$. 

\textit{Pairing.} For any $\barx, \bary \in V \llbracket X \rrbracket$, $\overline{\{\barx, \bary\}}$ will be the pair of $\barx, \bary$ in $V \llbracket X \rrbracket$.

\textit{Union.} Given any set $\barx$ in $V \llbracket X \rrbracket$, let $\bary = \<1, \{\barz : \exists \bar{w} \ \bar{\in} \ \barx \  (\barz \ \bar{\in} \ \bar{w}) \}>$. $V \llbracket X \rrbracket$ then thinks that $\bary$ is the  union of $\barx$.

\textit{Powerset.}
Given a $\barx = \<1, x>$ with $x \subseteq V \llbracket X \rrbracket$. Define $\bary = \langle 1, \{\langle 1, v\rangle : v \subseteq x\}\rangle$, which is the powerset of $\barx$ in $V \llbracket X \rrbracket$.

\textit{Infinity.} It is easy to check that $\hat{\omega}$ is an inductive set in $V \llbracket X \rrbracket$.

\textit{Separation.} Given a $\barx = \<1, x>$ with $x \subseteq V \llbracket X \rrbracket$ and a parameter $\bar{u} \in V \llbracket X \rrbracket$. Let $y =\{ \barz \in x : V \llbracket X \rrbracket \models \varphi(\barx, \barz, \bar{u} ) \}$. Then $\<1, \bary>$ is the desired subset of $\barx$ in $V \llbracket X \rrbracket$.

\textit{Collection.} Suppose that for every $\barx \ \bar{\in} \ \barw$, there is some $\bary$ such that $V \llbracket X \rrbracket \models \varphi (\barx, \bary, \bar{u})$, where $\bar{u} \in V \llbracket X \rrbracket$. By Collection and Separation in $V$, there is some $v \subseteq V \llbracket X \rrbracket$ such that for every $\barx \ \bar{\in} \ \bar{w}$, there is some $\bary \in v$ with $V \llbracket X \rrbracket \models \varphi (\barx, \bary, \bar{u})$. Then $\<1, v>$ is a desired collection set in $V \llbracket X \rrbracket$.
\end{proof}

\begin{lemma}\label{vhatisov}
Let $V$ be a model of ZF and $X$ be a class of $V$. Then $\hat{V}$ is the class of \textit{all} pure sets in $V \llbracket X \rrbracket$. Therefore, $V$ is isomorphic to $V^{V\llbracket X \rrbracket}$, i.e., the class of all pure sets in $V\llbracket X \rrbracket$.
\end{lemma}
\begin{proof}
I first show that $\hat{V} \subseteq V^{V\llbracket X \rrbracket}$ by an $\in$-induction. Assume that for all $y \in x$, $\hat{y}$ is a pure set in $V \llbracket X \rrbracket$. Since if $V \llbracket X \rrbracket \models \bar{z} \in trc(\hat{x})$ then there is some $y \in x$ such that $V \llbracket X \rrbracket \models \bar{z} \in  trc(\{\hat{y}\})$, it follows that $\hat{x}$ is a pure set in $V \llbracket X \rrbracket$. To show $V^{V\llbracket X \rrbracket} \subseteq \hat{V}$, we used an $\bar{\in}$-induction. Suppose that $\barx$ is a pure set in $V \llbracket X \rrbracket$ and for all $\bar{y} \ \bar{\in} \ \bar{x}$, $\bar{y} \in \hat{V}$. Let $v = \{z \in V : \exists \bar{y} \ \bar{\in} \ \bar{x}(\bar{y} = \hat{z})\}$, which is a set in $V$ because the map $z \mapsto \hat{z}$ is 1-1. Then $\barx = \hat{v}$.
\end{proof} 
\begin{theorem}\label{thm:ACholdsViffACholdsinV[X]}
Let $V$ be a model of ZF and $X$ be a class of $V$. Then
\begin{enumerate}
\item If a set $x \subseteq V \llbracket X \rrbracket$ is well-orderable in $V$, then $\barx = \<1, x>$ is well-oderable in $V \llbracket X \rrbracket$;
\item $V \models $ AC if and only if $V \llbracket X \rrbracket \models $ AC.
\end{enumerate}
\end{theorem}
\begin{proof}
(1) First note that if $x \sim v$ for some set $v \in V$, then $V \llbracket X \rrbracket \models \barx \sim \hat{v}$, where  $\barx = \<1, x>$. This is because any bijection $f$ between $x$ and $v$ in $V$ can be coded by $\bar{f} = \<1, f'>$, where $f' = \{\overline{\<\bary, \hat{w}>} : \bary \in x \land w \in v \land f(\bary) = w\}$ and $\overline{\<\bary, \hat{w}>} = \<1, \{\overline{\{\bary\}}, \overline{\{\bary, \hat{w}\}}\}>$. Then $\bar{f}$ will be a bijection between $\barx$ and $\hat{v}$ in $V\llbracket X \rrbracket$. Thus, if $x \sim \alpha$ for some ordinal $\alpha$, $V \llbracket X \rrbracket \models \bar{x} \sim \hat{\alpha}$ and $\hat{\alpha}$ is an ordinal in $V \llbracket X \rrbracket$ by Lemma \ref{vhatisov}.

(2) The left-to-right direction follows from (1). Suppose that in $V$, some set $x$ cannot be well-ordered. Then $\hat{V}$ thinks that $\hat{x}$ cannot be well-ordered by Lemma \ref{vhatisov}, so AC fails in $\hat{V}$ and hence in $V \llbracket X \rrbracket$.
\end{proof}

\begin{theorem}\label{con(zfc)->con(zfcu)}
The following theories are mutually interpretable.
\begin{enumerate}
    \item ZF.
    \item $\ZFCUR$ + Collection +  $\A \sim \omega$.
    \item $\ZFCUR$ + Collection +  $\A \sim \omega_1$.
    \item $\ZFCUR$ + Collection +  ``for every cardinal $\kappa$, there is a set of $\kappa$-many urelements''.
\end{enumerate}
\end{theorem}
\begin{proof}
It is clear that in any model of $\ZFCUR$, the class of pure sets is a model of ZFC. And if $V$ is a model of ZF, we can first go to its $L$ to have a model of ZFC. To get exactly $\omega$-many urelements, consider the model $L \llbracket  \omega \rrbracket$ in which the set of all urelements $\bar{\A}$ is $\<1, \{0\}\times \omega>$ and $L\llbracket  \omega \rrbracket \models \bar{\A} \sim \hat{\omega}$. Similarly, $L \llbracket  \omega_1 \rrbracket$ will be a model of $\ZFCUR$  with exactly $\omega_1$-many urelements. To get unboundedly many urelements, consider $L\llbracket  Ord \rrbracket$. For every $\omega_\alpha$, $\bar{A} = \<1, \{0\} \times \omega_\alpha>$ will then be a set of urelements of size $\omega_\alpha$ in $L \llbracket Ord \rrbracket$, because $L \llbracket Ord \rrbracket \models \hat{\omega_\alpha} = \omega_\alpha$ by Lemma \ref{vhatisov}.\end{proof}

In fact, we have shown that ZF is semi-bi-interpretable (see \cite{FreireForthcoming-FREBIW-2} for a definition of bi-interpretation) with the other urelement theories listed above since the map $x \mapsto \hat{x}$ is a definable isomorphism from $V$ to the class of pure sets in $V \llbracket  X \rrbracket$ for any class $X$. Also, if $U$ is a model of $\ZFCUR$  where the class of urelements $\A$ is a set, then in $U$ the urelements can be enumerated by a pure set $x$ so there is a definable isomorphism between $U$ and $V\llbracket x \rrbracket$. It is shown in \cite{HamkinsForthcoming-HAMRIS} that no model of $\ZFCUR$ + Collection with a proper class of urelements can be bi-interpretable with a model of ZFC. 
\chapter{Axioms in Set Theory with Urelements}
In this chapter, I investigate the axiomatization of set theory with urelements. Section \ref{section:additionalaxioms} introduces a group of additional axioms in urelement set theory together with the notion of homogeneity. In Section \ref{section:aHierarchyofAxioms}, I show that the group of axioms form a hierarchy over $\ZFCUR$. This gives rise to a natural question: what is ZFC set theory with urelements? In Section \ref{section:WhatisZFCU}, I present some evidence suggesting that a robust version of ZFC with urelements should include Collection as an axiom. I then consider various philosophical positions one might take regarding this issue. In Section \ref{section:ChoicelesZFU}, I explore urelement set theory without the Axiom of Choice. The hierarchy of axioms in $\ZFCUR$, as I shall prove, largely breaks down when sets of urelements are not necessarily well-orderable.

\section{Additional axioms and homogeneity}\label{section:additionalaxioms}
\subsection{Reflection and dependent choice schemes}
Reflection principles in set theory assert that the set-theoretic universe, $V$, is to some extent indescribable: any statement $\varphi$ will become absolute between some $V_\alpha$ and $V$. In particular, ZF proves the following L\'evy-Montague reflection principle.
\begin{itemize}
    \item [] $\forall \alpha \exists \beta > \alpha\forall x_1, ..., x_n \in V_{\beta} (\varphi(x_1, ..., x_n) \leftrightarrow \varphi^{V_\beta}(x_1, ..., x_n))$.
\end{itemize}
\noindent In urelement set theory, one cannot expect the L\'evy-Montague reflection principle to hold, e.g., if there is a proper class of urelements, then no $V_\alpha(A)$ can reflect such statement for any set of urelements $A$. Instead, in the presence of urelements it should be transitive sets that reflect. Namely,
\begin{itemize}
    \item [] (RP) For every set $x$ there is a transitive set $t$ with $x \subseteq t$ such that for every $x_1, ..., x_n \in t$, $\varphi(x_1, ..., x_n) \leftrightarrow \varphi^t(x_1, ..., x_n)$. \end{itemize}

\begin{prop}\label{RP->Collection}
 ZU + RP $\vdash$ Collection.
\end{prop}
\begin{proof}
First note that RP implies that for any finite collection of formulas $\varphi_1, ..., \varphi_n$, there is a transitive set reflecting each $\varphi_i$. For we can let $\psi(v)$ be the formula $(v=1 \land \varphi_1) \lor ... \lor (v=n \land \varphi_n)$, where $v$ is a new free variable. It follows from RP that there is transitive set $t$ set extending $\omega$ that reflects $\psi(v)$. Consequently, $\varphi_i^t \leftrightarrow \varphi_i$ for each $\varphi_i$. Now suppose that $\forall x \in w \exists y \varphi (x, y, u)$. Let $t$ be a transitive set extending $\{w, u\}$ which simultaneously reflects $\forall x \in w \exists y \varphi (x, y, u)$ and $\varphi(x, y, u)$. It follows that $t$ is a desired collection set. 
\end{proof}
There is another seemingly weaker version of RP, which asserts that any true statement is true in some transitive set containing the parameters.
\begin{itemize}
    \item [] (RP$^-$) If $\varphi (x_1, ..., x_n)$, then there is a transitive set $t$ containing $x_1, ..., x_n$ such that $\varphi^t(x_1, ..., x_n)$.
\end{itemize}
\noindent This form of reflection was first introduced by L\'evy \cite{Levy1966principles}. And in \cite{Levy1961principles} L\'evy and Vaught showed that over Zermelo set theory, RP$^-$ does not imply RP.

The Axiom of Dependent Choice (DC) asserts that for any set $x$, if $r \subseteq x \times x$ is a \textit{set relation} without terminal nodes, then there is an infinite sequence $s \in x^\omega$ such that $\<s(n), s (n+1)> \in r$ for every $n$. DC can be generalized to DC$_\kappa$ for any infinite cardinal $\kappa$ as follows (first introduced by L\'evy in \cite{bwmeta1.element.bwnjournal-article-fmv54i1p13bwm}). 
\begin{itemize}
    \item [] (DC$_\kappa$) For every $x$ and $r\subseteq x \times x$, if for every $s \in x^{<\kappa}$, there is some $w \in x$ such that $\<s, w> \in r$, then there is an $f: \kappa \rightarrow x$ such that $\<f\restriction\alpha, f(\alpha)> \in r$ for all $\alpha < \kappa$.
\end{itemize}
$\ZFCUR$ proves $\forall \kappa$DC$_\kappa$ because we can well-order $x$ and construct the desired function $f$ on $\kappa$ by recursion. What will interest us is the following class version of dependent choice.
\begin{itemize}
    \item [] (DC-scheme) If for every $x$ there is some $y$ such that $\varphi(x, y, u)$, then for every $p$ there is an infinite sequence $s$ such that $s(0) = p$ and $\varphi(s(n), s(n+1), u)$ for every $n<\omega$.
\end{itemize} 

\noindent The DC-scheme says that if $\varphi$ defines a class relation without terminal nodes, then there is an infinite sequence threading this relation. Similarly, we can formulate a class version of DC$_\kappa$ for any infinite cardinal $\kappa$.
\begin{itemize}
        \item[]  (DC$_\kappa$-scheme) If for every $x$ there is some $y$ such that $\varphi(x, y, u)$, then there is some function $f : \kappa\rightarrow U$ such that $\varphi(f\restriction \alpha, f(\alpha), u)$ for every $\alpha <\kappa$.
\end{itemize}
DC$_{<Ord}$ holds just in case the DC$_\kappa$-scheme holds for every $\kappa$. As we shall see, over $\ZFCUR$ the DC$_\kappa$-scheme is strictly stronger than DC$_\kappa$, and in fact, $\ZFCUR$ cannot even prove the DC$_\omega$-scheme. Let us verify that the DC-scheme is indeed a reformulation of the DC$_\omega$-scheme.
\begin{prop}
Over $\ZFUR$, the DC$_\omega$-scheme is equivalent to the DC-scheme.
\end{prop}
\begin{proof}
Assume the DC$_\omega$-scheme and suppose that the relation $\varphi(x, y,u)$ has no terminal nodes. Fix any $p$ and define $\psi(x, y, p, u)$ as follows.
\begin{enumerate}
    \item [] $\psi(x, y, p, u)$ $=_\textup{df}$ $(x = \emptyset \land y = p) \lor \exists n  (x = \langle x_0, ..., x_{n} \rangle \land \varphi (x_n, y, u)) \lor (x \neq \emptyset \land x$ is not a finite sequence$)$
\end{enumerate}
Since $\psi(x, y, p, u)$ has no terminal nodes, there is a function $f$ on $\omega$ such that $\psi(f\restriction n, f(n), p, u)$ for all $n$. It then follows that $f(0) = p$ and $\varphi(f(n), f(n+1), u)$ for all $n$.

Now assume the DC-scheme and suppose that $\varphi(x, y,u)$ has no terminal nodes. Fix some $p$ be such that $\varphi(\emptyset, p, u)$. Let $\psi(x, y, p, u)$ be the formula asserting that either $x$ is not a finite sequence, or $x = \emptyset$ and $y = \<p>$, or for some $n > 0$, $x$ is an $n$-sequence and $y$ is an $n+1$-sequence extending $x$ such that $\varphi(x, y(n), u)$. As $\psi(x, y, p, u)$ has no terminal nodes, it follows that there is an $\omega$-sequence $s$ such that $s(0)= \emptyset$ and $\psi(s(n), s(n+1), p, u)$ for every $n$. Let $f(n) = s(n+1)(n)$. Then $\varphi(f\restriction n, f(n), u)$ for every $n$.
\end{proof}

\begin{prop}\label{prop:DCkappaVariants}
Over $\ZFUR$, for every cardinal $\kappa$, the following are equivalent.
\begin{enumerate}
    \item DC$_\kappa$-scheme. 
    \item For every definable class $X$ and $u$, if every $s \in X^{<\kappa}$ has some $y \in X$ such that $\varphi(x, y, u)$, then there is a function $f \in X^{\kappa}$ such that $\varphi(f\restriction \alpha, f(\alpha), u)$ for every $\alpha < \kappa$. 
\end{enumerate}
\end{prop}
\begin{proof}
(2) $\rightarrow$ (1) is immediate. To show (1) $\rightarrow$ (2), assume (1) and the antecedent of (2). Define
\begin{itemize}
    \item []  $\psi(x, y, u) =_\textup{df} x \in X^{<\kappa} \rightarrow \varphi(x, y, u) \land y \in X.$
\end{itemize}
Then there is an $f : \kappa \rightarrow U$ such that $\psi(f\restriction \alpha, f(\alpha), u)$ for every $\alpha < \kappa$. Consider any $\alpha < \kappa$ and suppose that for every $\beta < \alpha$, $f(\beta) \in X$. Then $f(\alpha) \in X$. So $f\restriction \alpha \in X$ for every $\alpha <\kappa$. Thus, $f \in X^\kappa$  and $\varphi(f\restriction\alpha, f(\alpha), u)$ for every $\alpha < \kappa$.\end{proof}
\begin{theorem}[\cite{Gitman2016-GITWIT}]\label{collection+dc->rp}
$\ZFUR +$ DC$_\omega$-scheme + Collection $\vdash$ RP.
\end{theorem}
\begin{proof}
This is first proved in \cite{Gitman2016-GITWIT} in context of ZFC without Powerset, and the point here is that their argument works in $\ZFUR$. Fix a set $w$ and a formula $\varphi(v_1, ..., v_m)$ and let $\varphi_1, ..., \varphi_n$ be all the subformulas of $\varphi$.  For any set $x$ and $y$, we say that $y$ is a $\varphi$-cover of $x$ if (i) $y$ is transitive and $x \subseteq y$, and (ii) for each $1 \leq i \leq n$,  $\forall x_1, ..., x_m \in x \ [\exists v \varphi_i (v, x_1, ..., x_m) \rightarrow \exists v \in y\  \varphi_i (v, x_1, ..., x_m)]$.

Let us construct a $\varphi$-cover for an arbitrary set $x$. For each $1 \leq i \leq n$, define $z_i = \{s \in x^{<\omega}: \exists v \exists x_1, ..., x_m \in x (\varphi_i (v, x_1, ..., x_m) \land s = \langle x_1, ..., x_{m} \rangle)\}$. For each $ i \leq n$, by Collection there is a $y_i$ such that for every $s \in z_i$, there is $v \in y_i$ such that  $\exists x_1, ..., x_m \in x [\varphi_i (v, x_1, ..., x_m) \land s = \langle x, ..., x_m \rangle]$. Set $y = trc(\bigcup_{i\leq n} y_i \cup x)$. $y$ is then a $\varphi$-cover of $x$. Thus, by the DC$_\omega$-scheme, there is an $\omega$-sequence $\langle z_n : n< \omega \rangle$ such that $z_0 = w$ and $z_{n+1}$ is a $\varphi$-cover of $z_n$ for every $n$. Set $t = \bigcup_{n < \omega} z_n$. By a routine induction on the subformulas of $\varphi$, it follows that for every $x_0, ..., x_m \in t$, $\varphi^t (x_1, ..., x_m) \leftrightarrow \varphi (x_1, ..., x_m)$.
\end{proof}

\subsection{Urelement axioms and homogeneity}
\begin{definition}
A set $x$ is \textit{realized} if there is a set of urelements equinumerous with $x$. Let $A$ be a set of urelements.
\begin{enumerate}
    \item A set of urelements $B$ \textit{duplicates} $A$, if $B$ and $A$ are disjoint and equinumerous.
    \item A set of urelements $B$ is a \textit{tail} of $A$, if $B$ is disjoint from $A$ and for every $C \subseteq \A$ disjoint from $A$ there is an injection from $C$ to $B$.
    \item \textit{Duplication holds over} $A$, if every $B \subseteq \A$ disjoint from $A$ has a duplicate $C$ that is also disjoint from $A$;
    \item \textit{Homogeneity holds over} $A$, if whenever $B$ and $C$ are equinumerous and $B \cup C \subseteq \A$ is disjoint from $A$, there is an automorphism $\pi$ such that $\pi B = C$ and $\pi$ point-wise fixes $A$.
\end{enumerate}
\end{definition}
\noindent Intuitively, when homogeneity holds over $A$, the urelements outside $A$ are all indistinguishable from the perspective of $A$. We shall consider the following axioms.
\begin{enumerate}
    \item [] (Plenitude) Every cardinal is realized.
    \item [] (Closure) The supremum of a set of realized cardinals is realized.
    \item [] (Duplication) Every set of urelements has a duplicate.
    \item [] (Tail) Every set of urelements has a tail.
    \item [] ($\ACA$) Every set of urelements is well-orderable.
\end{enumerate}
$\ACA$ is strictly weaker than AC over $\ZFUR$ even if we assume there are infinitely many urelements: let $V$ be a model of ZF $+$ $\neg$AC and consider $V\llbracket Ord \rrbracket$ (Definition \ref{barwiseinterpretation1}); by Theorem \ref{thm:ACholdsViffACholdsinV[X]}, $V\llbracket Ord \rrbracket \models \ACA \land \neg$AC. With $\ACA$, Tail is equivalent to the following.
\begin{itemize}
    \item [] (Tail$^*$) For every $A \subseteq \A$, there is a greatest cardinal $\kappa$ such that $\exists B \subseteq \A \ (B \sim \kappa \land B \cap A = \emptyset)$,
\end{itemize}
where $\kappa$ is called the \textit{tail cardinal} of $A$.
\begin{lemma}[$\ZFUR$]\label{homogeneitylemma}
\ 
\begin{enumerate}
    \item If $A \subseteq A' \subseteq \A$ and duplication holds over $A$, then duplication holds over $A'$.
    \item If duplication holds over $A \subseteq \A$, then homogeneity holds over $A$.
    \item Assume Tail and $\ACA$. Duplication holds over some set of urelements.
    \item Assume $\ACA$. For every $A \subseteq \A$, there is an $A' \subseteq \A$ such that $A \subseteq A'$ and duplication (hence homogeneity) holds over $A'$.
\end{enumerate}
\end{lemma}
\begin{proof}
(1) If $B$ is disjoint from $A'$, then there is another $C$ disjoint from $A$ that is equinumerous with $(A'\setminus A) \cup B$. So $C$ contains a duplicate of $B$ that is disjoint from $A'$.

(2) This is first proved in \cite{HamkinsForthcoming-HAMRIS}. Let $B$ and $C$ be such that $B \sim C$ and $B \cup C$ is disjoint from $A$. If $B$ and $C$ are duplicates, then by swapping them we can get a permutation $\pi$ with $\pi B = C$ that point-wise fixes $A$. If not, then by duplication over $A$, there is a duplicate of $B \cup C$ disjoint from $A \cup B \cup C$ and hence a duplicate $D$ of both $B$ and $C$ that is disjoint from $A$. Thus, there are automorphisms $\pi_1$ and $\pi_2$ such that $\pi_1 B = D$ and $\pi_2 D = C$, and both of them point-wise fix $A$. The composition of $\pi_1$ and $\pi_2$ is the desired automorphism.

(3) Assume Tail and $\ACA$. Let $\kappa$ be the least cardinal that is a tail cardinal of some $D \subseteq \A$. If $B$ is disjoint from $D$, since the tail cardinal of $D \cup B$ is at least $\kappa$, there is another $C$ disjoint from $D \cup B$ that has size at least $\kappa$. Since $B$ is well-orderable and hence has size at most $\kappa$, $C$ contains a duplicate of $B$. Therefore, duplication holds over $D$.

(4) Assume $\ACA$. By (1) and (2), it suffices to show that duplication holds over some set of urelements. Suppose otherwise. Then $\A$ is a proper class, and by (3), some $A \subseteq \A$ has no tail cardinal. Given any infinite $B$ disjoint from $A$, since $B$ is equinumerous with some cardinal $\kappa$, there must be some $C$ of size $\kappa^+$ that is disjoint from $A$. But then C contains a duplicate of $B$ . This shows that duplication holds over $A$ after all, which is a contradiction.
\end{proof}
\noindent In general, homogeneity does not imply duplication: if $\A$ is a set and $\A \setminus A$ has only one urelement, then duplication does not hold over $A$ while homogeneity holds over $A$ trivially. And assuming $\ACA$, duplication holds over $A$ if homogeneity holds over $A$ and $\A \setminus A$ is not a finite non-empty set. To see this, let $B$ be a set of urelements disjoint from $A$, and we may assume $B$ is infinite. By $\ACA$, $B$ can be partitioned into a pair of duplicates $B_1$ and $B_2$. By homogeneity over $A$, there is an automorphism $\pi$ such that $\pi B_1 = B$ and $\pi A =A$. Since $B_1$ has a duplicate disjoint from $A$, it follows that $B$ has one too and hence duplication holds over $A$. In Section \ref{section:ChoicelesZFU}, I will show that $\ACA$ is necessary for homogeneity to hold over some set of urelements. In particular, there are models of $\ZFUR$ + RP + DC$_\omega$-scheme where homogeneity holds over no set of urelements.

\section{A hierarchy of axioms in $\ZFCUR$}\label{section:aHierarchyofAxioms}

\subsection{Implication diagram in $\ZFCUR$}\label{subsection:implicationdiagram}
The main theorem of this section is the following.
\begin{theorem}\label{maintheorem1}
Over $\ZFCUR$, the following implication diagram holds. Moreover, the diagram is complete: if the diagram does not indicate $\varphi$ implies $\psi$, then $\textup{ZFCU}_\text{R}$ $+$ $\varphi \nvdash \psi$ assuming the consistency of ZF.
\end{theorem}
\begin{figure}[hbt!]
\begin{center}
\begin{tikzpicture}
\begin{scope}[every node/.style={}]
    \node (A) at (10, -2.5) {DC$_\omega$-scheme};
    \node (B) at (10, 1.5){Tail};
    \node (C) at (3,3) {Plenitude};
     \node (D) at (6,1.5) { DC$_{<Ord}$ };
    \node (E) at (6, -0.5) {DC$_\kappa$-scheme};
       \node (F) at (12, -2.5) {Closure};
          \node (G) at (11.75, -0.5) {RP};

    \node (H) at (3, 1.5) {\textup{Closure$\land$Duplication}};
    \node (I) at (10, -0.5) {Collection};
        \node (J) at (13, -0.5) {RP$^-$};

    \node (L) at (6, 0.7) {.};
    \node (M) at (6, 0.5) {.};
    \node (N) at (6, 0.3) {.};
    \node (O) at (3, -0.5) {Duplication};
    \node (P) at (6, -1.7) {.};
    \node (Q) at (6, -1.5) {.};
    \node (R) at (6, -1.3) {.};
    \node (S) at (6, -2.5) {DC$_{\omega_1}$-scheme};
    \node (T) at (10, 3) {$\A$ is a set};
    
\end{scope}

\begin{scope}[>={stealth},
              every node/.style={fill=white,circle},
              every edge/.style={draw=black}]
 
    \path [->] (T) edge (B);
    \path [->] (T) edge (D);
    \path [->] (C) edge (H);
    \path [->] (B) edge (I);
    \path [->] (I) edge (A);

     \path [->] (S) edge (A);

    \path [->] (H) edge (I);
    \path [->] (H) edge (O);
    
     \path [->] (C) edge (D);

    \path [->] (D) edge (I);
    \path [->] (D) edge (L);
    \path [->] (N) edge (E);
    \path [->] (E) edge (R);
     \path [->] (P) edge (S);
    
    \path [->] (I) edge (F);

    \path [->] (I) edge (G);
   \path [->] (G) edge (J);
   
       \path [->] (J) edge[bend right=30] (I); 
\end{scope}
\end{tikzpicture}
 \caption{Implication diagram in $\ZFCUR$}
 \label{ZFCUdiagram}
\end{center}
\end{figure}
\FloatBarrier
\noindent The direction from Collection to the DC$_\omega$-scheme was first proved by Schlutzenberg in an answer to a question on Mathoverflow \cite{387471} and the notion of tail cardinal was also implicit in his proof (\cite{HamkinsForthcoming-HAMRIS} also contains a different proof of Collection $\rightarrow$ DC$_\omega$-scheme). My proof of Collection $\rightarrow$ DC$_\omega$-scheme takes a different route and appeals to a key observation that Tail implies Collection, which will also be crucial for later discussions.

The rest of this subsection establishes the implication diagram, while the next subsection proves its completeness. Let us first show that Plenitude implies DC$_{<Ord}$. Given a formula $\varphi(x, y, u)$ with a parameter $u$, for any ordinals $\alpha, \alpha', \kappa, \kappa'$ and a set of urelements $E$, we say that $\<\kappa' \alpha'>$ is a $(\varphi, E)$\textit{-extension} of $\<\kappa, \alpha>$ if (i) $\alpha \leq \alpha'$, and (ii) whenever $A \subseteq \A$ extends $E$ by $\kappa$-many urelements, there is some $B \subseteq \A$ disjoint from $A$ with $B \sim \kappa'$ such that for every $x \in V_{\alpha}(A)$, there is some $y \in V_{\alpha'} (A \cup B)$ such that $\varphi (x, y, u)$. 
\begin{lemma}($\ZFCUR$)\label{phiextension}
Suppose that Plenitude holds and $\varphi(x, y, u)$ defines a relation without terminal nodes. Then every $\langle \kappa, \alpha \rangle$ has a $(\varphi, ker(u))$-extension.
\end{lemma}
\begin{proof}
First note that under $\ZFCUR$ + Plenitude, homogeneity holds over every set of urelements. Fix$\langle \kappa, \alpha \rangle$ and some $A \subseteq \A$ extending $ker(u)$ with $\kappa$-many urelements. For each $x\in V_\alpha(A)$, define $\theta_x$ to be the least cardinal such that there is some $y$ with $\varphi(x, y, u)$ and $ker(y) \sim \theta_x$, and let $\kappa' = Sup\{ \theta_x : x \in  V_\alpha(A)\}$. Fix some infinite $B$ of size $\kappa'$ that is disjoint from A, which exists by Plenitude. Then for every $ x \in V_\alpha(A)$, fix some $y'$ such that $\varphi(x, y', u)$ and $ker(y') \sim \theta_x$. $ker(y') \setminus A$ is equinumerous to a subset of $B$, so by homogeneity over $A$, there is an automorphism $\pi$ that moves $ker(y')$ into $B$ and point-wise fixes $A$. It follows that $\varphi(x, \pi y', u)$ and $\pi y' \in V(A \cup B)$. Thus, each $x \in V_\alpha(A)$ has some $y \in V(A\cup B)$ with $\varphi(x, y, u)$, so there is some large enough $\alpha'$ such that every $x \in V_\alpha(A)$ has some $y \in V_{\alpha'}(A\cup B)$ with $\varphi(x, y, u)$. Furthermore, for every $A'$ extending $ker(u)$ by $\kappa$-many urelements, by homogeneity over $ker(u)$, there is an automorphism $\pi$ with $\pi A = A'$ that point-wise fixes $\ker(u)$; so $\pi B$ will be such that every $x \in V_\alpha(A')$ has some $y \in V_{\alpha'}(A' \cup \pi B)$ with $\varphi(x, y, u)$. Therefore, $\<\kappa', \alpha'>$ is indeed a $(\varphi, ker(u))$-extension of $\<\kappa, \alpha>$.
\end{proof}

\begin{theorem}\label{Plenitude->DCS}
 $\ZFCUR$ $\vdash$Plenitude $\rightarrow$ DC$_{<Ord}$.
\end{theorem}
\begin{proof}
Suppose that Plenitude  holds and  $\varphi(x, y, u)$ defines a relation without terminal nodes with some parameter $u$. Consider any infinite cardinal $\kappa$. To prove the DC$_\kappa$-scheme, we first find a set $\barx$ that is closed under $<\kappa$-sequences and the relation $\varphi$; we can then apply DC$_\kappa$ to get a desired function on $\kappa$. Let $\delta$ be a cardinal with $\textup{cf}(\delta) = \kappa$. We first define a $\delta$-sequence of pairs of ordinals $\langle \langle \lambda_\alpha, \gamma_\alpha \rangle : \alpha < \delta \rangle$ by recursion as follows. Let $A_0$ be a set of urelements that extends $\ker(u)$ by $\lambda_0$-many urelements and $\gamma_0$ be an ordinal with cf$(\gamma_0)  \geq \kappa$. For each ordinal $\alpha < \delta$, we let $\langle \lambda_{\alpha+1} , \gamma_{\alpha+1} \rangle$ be the lexicographical-least $(\varphi, ker(u))$-extension of $\<\lambda_\alpha ,\gamma_\alpha>$ with cf$(\gamma_\alpha) \geq \kappa$, which exists by the previous lemma. And we take the union at the limit stage.

By Plenitude, we can fix a $\delta$-sequence of sets of urelements $\langle A_\alpha : \alpha < \delta \rangle$, where $A_\alpha$ extends $\bigcup_{\beta < \alpha} A_\beta \cup ker(u)$ by $\lambda_\alpha$-many urelements. Let $\barx = \bigcup_{\alpha < \delta} V_{\gamma_\alpha} (A_\alpha)$. For any $x \in V_{\gamma_\alpha} (A_\alpha)$, There is some $B$ disjoint from $A_\alpha$ witnessing the fact that $\<\lambda_{\alpha +1}, \gamma_{\alpha+1}>$ is a $(\varphi, ker(u))$-extension of $\<\lambda_\alpha ,\gamma_\alpha>$. And by homogeneity over $A$, it follows that $A_{\alpha +1} \setminus A_\alpha$ works as such witness as well; so there is some $y \in V_{\gamma_{\alpha +1}}(A_{\alpha +1})$ with $\varphi(x, y, u)$, and such $y$ lives in $\barx$. $\barx$ is also closed under $<\kappa$-sequences since $\textup{cf}(\delta) = \kappa$ and each $V_{\gamma_\alpha} (A_\alpha)$ is closed under $<\kappa$-sequences. Thus, if $s \in \barx^{<\kappa}$, there is some $y \in \barx$ such that $\varphi (s, y, u)$. By DC$_\kappa$, there exists a function $f$ on $\kappa$ such that $\varphi (f\restriction \alpha, f(\alpha), u)$ for all $\alpha < \kappa$. Hence, the DC$_\kappa$-scheme holds.
\end{proof}

\begin{lemma} \label{Plenitude->Collection}
$\ZFUR + \ACA$ $\vdash$ Closure $\land$ Duplication $\rightarrow$ Collection
\end{lemma}
\begin{proof}
Fix some set $w$ such that $\forall x \in w \exists y \varphi (x, y, u)$. For every $x \in w$, let $\theta_x$ be the least $\theta$ realized by the kernel of some $y$ such that $\varphi(x, y, u)$, and define $\theta$ as the supremum of all such $\theta_x$. Let $A \subseteq \A$ be such that $ker(w) \cup \ker(u) \subseteq A$ and duplication holds over $A$,  which exists by Lemma \ref{homogeneitylemma} (4). By Closure and Duplication, there is a $B \subseteq \A$ of size $\theta$ that is disjoint from $A$. Then for every $x \in w$, fix a $y'$ such that $\varphi(x, y', u)$ with the smallest kernel. By homogeneity over $A$, there is an autormophism that moves $ker(y')$ into $A \cup B$ without moving any urelements in $A$. Therefore, every $x \in w$ has a $y \in V(A \cup B)$ such that $\varphi(x, y, u)$. Then Collection holds by applying Proposition \ref{weakcollection}.
\end{proof}

\begin{lemma}\label{tail->collection}
$\ZFUR + \ACA$ $\vdash$ Tail $\rightarrow$ Collection
\end{lemma}
\begin{proof}
Assume that every set of urelements has a tail. Suppose that every $x \in w $ has some $y$ with $\varphi (x, y, u)$. Let $A \subseteq \A$ be such that $ker(w) \cup ker(u) \subseteq A$ and duplication holds over $A$ and $B$ be a tail of $A$. For every $x \in w$ and $y$ such that $\varphi(x, y, u)$, $B$ must contain a subset that is equinumerous with $ker(y) \setminus A$. By homogeneity over $A$, there is an automorphism that moves $ker(y)$ into $A \cup B$ while point-wise fixing $A$. Therefore, every $x \in w$ has some $y \in V(A \cup B)$ such that $\varphi (x, y, u)$ and hence Collection holds by Proposition \ref{weakcollection}.
\end{proof}

\begin{lemma}[$\ZFCUR$]\label{Tailkappa->DCkappa}
Let $\kappa$ be an infinite cardinal and suppose that every set of urelements has a tail of size at least $\kappa$. Then the DC$_\kappa$-scheme holds.
\end{lemma}
\begin{proof}
First assume that $\kappa$ is regular. Suppose that $\varphi(x, y, u)$ defines a relation without terminal nodes with a parameter $u$. Let $A$ be a set of urelements extending $ker(u)$ over which duplication holds and $B$ be a tail of $A$. Since $B$ has size at least $\kappa$, $B$ can be partitioned into $\kappa$-many pieces $\{B_\eta : \alpha < \kappa\}$, where each $B_\eta$ is equinumerous with $B$. Let $\beta$ be a ordinal such that cf$(\beta) = \kappa$. We define a $\kappa$-sequence of ordinals $\<\gamma_\alpha : \gamma < \kappa>$ above $\beta$ by recursion, where $\gamma_\alpha$ is the least ordinal such that 
\begin{itemize}
    \item [] (i) $\gamma_\alpha > \bigcup_{\eta < \alpha} \gamma_\eta$ and cf$(\gamma_\alpha) = \kappa$;
    \item [] (ii) for every $x$ in $ \bigcup_{\eta < \alpha} V_{\gamma_\eta}(\bigcup_{\eta < \alpha}B_\eta \cup A)$, there is a $y \in V_{\gamma_\alpha}(\bigcup_{\eta \leq \alpha} B_\eta \cup A)$ with $\varphi(x, y, u)$.
\end{itemize}
Such $\gamma_\alpha$ exists because homogeneity holds over $\bigcup_{\eta < \alpha}B_\eta \cup A$ and each $B_\alpha$ is a tail of of $A$. Let $x = \bigcup_{\alpha < \kappa} V_{\gamma_\alpha}(\bigcup_{\eta \leq \alpha}B_\eta \cup A)$. $x$ is then closed under $\varphi(x, y, u)$. And since $x$ is the union of an increasing $\kappa$-sequence of sets and each $\gamma_\alpha$ has cofinality $\kappa$, it follows that $x^{<\kappa} \subseteq x$. We can then apply DC$_\kappa$ to $x$ to get a desired $\kappa$ sequence, so the DC$_\kappa$-shceme holds.

Suppose $\kappa$ is singular. Then for every regular $\lambda < \kappa$, the argument in the previous paragraph shows that the DC$_\lambda$-scheme holds. But this implies the DC$_\kappa$-scheme by a standard argument as in \cite[Theorem 8.1]{jech2008axiom}.
\end{proof}

\noindent To show that Diagram \ref{ZFCUdiagram} holds, it remains to prove the following.

\begin{lemma}\label{easyimplication}Over $\ZFCUR$, the following implications hold.
\begin{enumerate}
    \item $\mathcal{A}$ is a set $\rightarrow$ DC$_{<Ord}$. 
    \item DC$_{<Ord}$ $\rightarrow$ Collection
    \item RP$^-$ $\rightarrow$ Collection.
    \item Collection $\rightarrow$ Closure
    \item Collection $\rightarrow$ DC$_\omega$-scheme.
    \item Collection $\rightarrow$ RP.
\end{enumerate}
\end{lemma}
\begin{proof}

(1) This is proved by a standard argument, which I include for completeness. Assume $\mathcal{A}$ is a set and $\forall x \exists y \varphi (x, y, u)$. Fix any $\kappa$ and let $\delta$ be such that cf$(\delta)= \kappa$. We define a $\delta$-sequence of ordinals $\langle \gamma_\alpha : \alpha < \delta \rangle$, where $\gamma_\alpha$ is the least ordinal of cofinality $\kappa$ such that $\forall x \in \bigcup_{\eta <  \alpha} V_{\gamma_\eta} (\mathcal{A}) \exists y \in V_{\gamma_\alpha}(\mathcal{A}) \varphi (x, y, u)$. Then the DC$_\kappa$-scheme holds by applying DC$_\kappa$ to $\bigcup_{\alpha < \delta} V_{\gamma_\alpha} (\mathcal{A})$.

(2) This is because under DC$_{<Ord}$, either $\A$ is a set or Plenitude holds, but Collection holds either way by Proposition \ref{weakcollection} and Lemma \ref{Plenitude->Collection}.

(3) Suppose that RP$^-$ holds. It suffices to show that Tail holds by Lemma \ref{tail->collection}. We may assume that  $\mathcal{A}$ is not a set and Plenitude fails by Lemma \ref{Plenitude->Collection}.  Fix some $A \subseteq \A$ and let $\kappa$ be the least cardinal not realized by some $B \subseteq \A$ that is disjoint from $A$. Then there is a transitive set $t$ reflecting the statement that $\forall \lambda < \kappa \exists B (B \sim \lambda \land B \cap A = \emptyset )$. We may assume that $t$ extends $\{\kappa, A\}$ and is closed under pairs. $C = \bigcup \{B \in t : B \subseteq \A \land B \cap A = \emptyset\}$ is then a tail of $A$.

Now assume Collection. 

(4) Let $x$ be a set of realized cardinals. Then there is a set $y$ such that for every $\kappa \in x$, there is some $A \in y$ such that $A \sim \kappa$. Let $B = \bigcup \{A: A \in y\}$. Then the cardinality of $B$ is at least the supremum of $x$ and hence Closure holds. 

(5) Observe that Collection + $\neg$Plenitude implies Tail. Given a set $A$ of urelements, let $w$ be the set of cardinals realized by some $B \subseteq \A$ disjoint from $A$. Then there is some $v$ such that for every $\lambda \in w$, there is some $B \in v$ such that $B \sim \lambda$ and $B \cap A = \emptyset$. $C = \bigcup \{B \in v: B \cap A = \emptyset\}$ is then a tail of $A$. Now to show the DC$_{\omega}$-scheme holds, we may assume that Plenitude fails by Theorem \ref{Plenitude->DCS}. Then every set of urelements must have an infinite tail, so the DC$_{\omega}$-scheme follows from Lemma \ref{Tailkappa->DCkappa}.

(6) RP holds by (5) and Theorem \ref{collection+dc->rp}.
\end{proof}

\subsection{Independence results}
I now proceed to show that Diagram \ref{ZFCUdiagram} is complete by an easy method of building inner models of $\ZFCUR$, which was implicitly used in \cite{levy1969definability} and \cite{felgner1976choice}.
 \begin{definition}[$\ZFUR$]\label{normalideal}
A (definable) class $\I$ of sets of urelements is an $\A$-\textit{ideal} if 
\begin{enumerate}
    \item $\A \notin \I$ (if $\A$ is a set);
    \item if $A, B \in \I$, then $A \cup B \in \I$; 
    \item if $A \in \I$ and $B \subseteq A$, then $B \in \I$;
    \item for every $a \in \A$, $\{a\} \in \I$.
\end{enumerate}
Given an $\A$-ideal $\I$, $U^\I = \{x \in U : ker(x) \in \I\}$, i.e., the class of objects whose kernel is in $\I$.
 \end{definition}

\begin{lemma}[$\ZFUR$]\label{idealpermutation}
Let $\I$ be an $\A$-ideal. Then for every $a, A$ such that $a \in A \in \I$, there is a permutation $\pi$ of $\mathcal{A}$ such that (i) $\pi \I = \I$, (ii) $\pi a \neq a$ and (iii) $\pi$ point-wise fixes $A \setminus \{a\}$.
\end{lemma}
\begin{proof}
Fix some $a^* \in \mathcal{A} \setminus A$. Let $\pi$ be a permutation that only swaps $a$ and $a^*$. To see that $\pi \I = \I$, let $B \in \I$. Without lost of generality, we may assume $a \in B$ and $a^* \notin B$. Then $\pi B = (B \setminus \{a\}) \cup \{a^*\}$, which is in $\I$. Also, $B = \pi ((B \setminus \{a\}) \cup \{a^*\})$. Therefore, $\pi \I = \I$.
\end{proof}

\begin{theorem}[$\ZFUR$]\label{smallkernelmodel}
Let $\I$ be an $\A$-ideal.
\begin{enumerate}
\item $U^\I \models \ZFUR$ + ``$\A$ is a proper class'';
\item $U^\I \models $ AC if $U \models $ AC.
\end{enumerate}
\end{theorem}
\begin{proof}
It is clear that $U^\I$ is transitive and contains all the urelements and pure sets. Thus, $U^\I$ satisfies Foundation, Extensionality, Infinity, and $\A$ is a proper class in  $U^\I$. It is also immediate that $U^\I$ satisfies Pairing, Union, Powerset and Separation. When AC holds in $U$, it holds in $U^\I$ because for a given set $x$ in $U^\I$, any bijection in $U$ between $x$ and an ordinal has the same kernel as $x$ and hence also lives in $U^\I$. It remains to show that Replacement holds in $U^\I$.

Suppose that $U^\I \models \forall x \in w \exists ! y \varphi (x, y, u)$ for some $w, u \in U^\I$. Let $v = \{ y \in U^\I : \exists x \in w\ \varphi^{U^\I}(x, y, u) \}$, which is a set in $U$. It suffices to show that $ker(v) \subseteq  ker(w) \cup ker(u)$. Suppose not. Then there are some $y$ and $a$ such that $y \in v$, $a \in ker(y)$\footnote{Note that by our convention (Definition \ref{def:kernel,cardinal}) if $y$ is an urelement, $ker(y) = \{y\}$.} and $a \notin ker(w) \cup ker(u)$. Let $A = ker(w) \cup ker(u) \cup ker(y)$, which is in $\I$. By Lemma \ref{idealpermutation}, there is an automorphism $\pi$ such that (i) $\pi \I = \I$, (ii) $\pi a \neq a$ and (iii) $\pi$ point-wise fixes $A \setminus \{a\}$. So $\pi$ point-wise fixes $w$ and $u$. Since $y \in v$, there is some $x \in w$ with $\varphi^{U^\I}(x, y, u)$. It follows that $\varphi^{U^\I}(x, \pi y, u)$, but $\pi y \neq y$ because $\pi a$ is in $ker(\pi y)$ but not in $ker(y)$, which contradicts the uniqueness of $y$.
\end{proof}

\begin{theorem}\label{zfcurindependece}
Assume the consistency of ZF. Over $\ZFCUR$,
\begin{enumerate}
     \item (Closure $\land$ Duplication) $\nrightarrow$ (Plenitude $\lor$ DC$_{\omega_1}$-scheme);
    \item  Collection $\nrightarrow$ Duplication;
    \item  Duplication $\nrightarrow$ (Closure $\lor$ DC$_\omega$-scheme);
    \item  Closure $\nrightarrow$ DC$_\omega$-scheme;
    \item  DC$_\kappa$-scheme $\nrightarrow$ Closure, where $\kappa$ is any infinite cardinal;
    \item (Collection $\land$ DC$_\kappa$-scheme) $\nrightarrow$ DC$_\lambda$-scheme, where $\kappa < \lambda$ are infinite cardinals.
\end{enumerate}
Hence, Diagram \ref{ZFCUdiagram} is complete.
\end{theorem}
\begin{proof}
In each case, $U$ is a model of $\ZFCUR$. These models exist if ZF is consistent by Theorem \ref{con(zfc)->con(zfcu)}.

(1) Assume that in $U$, $\mathcal{A} \sim \omega_1$. Let $\I_1$ be the ideal of all countable subsets of $\mathcal{A}$. In $U^{\I_1}$, $\omega$ is the greatest realized cardinal. It is clear that Closure Duplication hold while Plenitude fails. The DC$_{\omega_1}$-scheme fails in $U^{\I_1}$ because every kernel can be properly extended but there cannot be a function $f$ on $\omega_1$ such that $ker(f\restriction \alpha) \subsetneq ker(f(\alpha))$ for all $\alpha < \omega_1$, as the kernel of such $f$ would be uncountable.  

(2) Assume that in $U$, $\mathcal{A} \sim \omega_1$. Fix an $A \subseteq \A$ such that $A \sim \omega_1$ and $\A \setminus A \sim \omega_1$. Let $\I_2 = \{ B \subseteq \mathcal{A}:  B  \setminus A \text{ is countable} \}$. For every $B \in U^{\I_2}$, let $\lambda =$Max$\{|A \setminus B|, \omega\}$, where $|A\setminus B|$ is the cardinality of $A\setminus B$. $\lambda$ is then the tail cardinal of $B$. So Collection holds in $U^{\I_2}$ by Lemma \ref{tail->collection}. Duplication fails because $A$ has no duplicates in $U^{\I_2}$.

(3) Assume that in $U$, $\mathcal{A} \sim \omega$. Let $\I_3$ be the ideal of finite subsets on $\mathcal{A}$. It is cleat that in $U^{\I_3}$ Duplication holds and Closure fails. The DC$_\omega$-scheme also fails in $U^{\I_3}$ because set of urelements can be properly extended but there is no infinite increasing sequence of sets of urelements.

(4) Assume that in $U$, $\mathcal{A} \sim \omega$ and fix an infinite and co-infinite $A \subseteq \A$. Let $\I_4 = \{ B \subseteq \mathcal{A}: B \setminus A \text{ is finite}\}$. Closure holds in $U^{\I_4}$ because $\omega$ is the greatest realized cardinal. The DC$_\omega$-scheme fails in $U^{\I_4}$ since every set of urelements can be properly extended by another set of urelements disjoint from $A$. but there cannot be a corresponding infinite sequence.

(5) Let $\kappa$ be an infinite cardinal. Assume that in $U$, $\mathcal{A} \sim \omega_{\kappa^+}$. Let $\I_5$  be the set of sets of urelements of size less than $\omega_{\kappa^+}$. Closure fails in $U^{\I_5}$ because $ \omega_{\kappa^+}$ is not realized while every cardinal below it is realized. To show that the DC$_\kappa$-scheme holds, suppose that for every $x \in U^{\I_5}$, there is some $y \in U^{\I_5}$ such that $\varphi^{U^{\I_5}}(x, y, u)$. $U^{\I_5}$ is closed under $\kappa$-sequences. Since DC$_{<Ord}$ holds in $U$ by Lemma \ref{easyimplication}, in $U$ there is a function $f: \kappa \rightarrow U^{\I_5}$ such that $\varphi^{U^{\I_5}}(f\restriction\alpha, f(\alpha), u)$ for every $\alpha < \kappa$, and $f$ lives in $U^{\I_5}$.

(6) It suffices to show that for any $\kappa$, $\ZFCUR$ + Collection + the DC$_\kappa$-scheme does not prove the DC$_{\kappa^+}$-scheme. Assume that in $U$, $\A \sim \kappa^+$ and let $\I_6$ be the ideal of all sets of urelements of size less than $\kappa^+$. By an argument as before, the DC$_{\kappa^+}$-scheme fails in $U^{\I_6}$. Every set of urelements in $U^{\I_6}$ has tail cardinal $\kappa$, so Collection holds by Lemma \ref{tail->collection} and the DC$_\kappa$-scheme holds by Lemma \ref{Tailkappa->DCkappa}.
\end{proof}
Combining the $U^\I$-construction with the $V\llbracket X \rrbracket$-construction in Definition \ref{barwiseinterpretation1}, we can further establish the mutual interpretability between various extensions of $\ZFCUR$.
\begin{corollary}\label{zfcumutualinter}
 The following theories are pairwise mutually interpretable.  
 \begin{enumerate}
     \item $\ZFCUR$ + Plenitude + $\neg$CH.
     \item $\ZFCUR$ + Collection + ``every set of urelements is countable'' + ``$\A$ is not a set'' + CH
     \item $\ZFCUR$ + Collection +  DC$_{\omega_1}$-scheme + ``$\A$ is not a set'' + $\neg$Plenitude.
 \end{enumerate}
 \end{corollary}
\begin{proof}
\

\textit{(1) interprets (2)}. In (1), we can go to its constructible universe $L$ and consider $L\llbracket \omega_1 \rrbracket$, which will be a model of $\ZFCUR$ + ``$\A \sim \omega_1$''. $L\llbracket \omega_1 \rrbracket \models$ CH because $L$ is isomorphic to the pure sets of $L\llbracket \omega_1 \rrbracket$ by Lemma \ref{vhatisov}. In $L\llbracket \omega_1 \rrbracket$, let $\I$ be the ideal of countable subsets of $\A$ (i.e., $\{0\} \times \omega_1$). By Theorem \ref{smallkernelmodel} and Lemma \ref{tail->collection}, $L\llbracket \omega_1 \rrbracket^\I$ is a model of theory (2). In particular, $L\llbracket \omega_1 \rrbracket^\I \models $ CH because it has the same pure sets as $L\llbracket \omega_1 \rrbracket$.

\textit{(2) interprets (1)}. It is known that given any model $V$ of ZF, we can construct a definable interpreted model $W$ of ZFC + $\neg$CH by the Boolean ultrapower construction (see \cite[Theorem 7]{FreireForthcoming-FREBIW-2}). So we can simply consider $W\llbracket Ord \rrbracket$ for such $W$. $W\llbracket Ord \rrbracket \models \neg$CH by Lemma \ref{vhatisov}. The rest of theorem can be proved by using the same method.
\end{proof}

\section{What is ZFC with urelements?}\label{section:WhatisZFCU}
$\ZFCUR$ thus proves none of the axioms in Diagram \ref{ZFCUdiagram}.\footnote{The situation here is very similar to the axiomatizations of certain fragments of ZFC. For example, in both ZFC without Powerset and intuitionistic ZF, Replacement does not imply Collection over the remaining axioms (see \cite{Zarach1996:ReplacmentDoesNotImplyCollection} and \cite{FRIEDMAN19851} respectively). And when ZFC without Powerset is formulated with only Replacement, as shown in \cite{Gitman2016-GITWIT}, it turns out to have various pathological models, all of which can be excluded by Collection. For this reason, it is argued in \cite{Gitman2016-GITWIT} that ZFC without Powerset \textit{should} be axiomatized with Collection.} A natural response at this point is to view $\ZFCUR$ as an inadequate way of formalizing ZFC with urelements. In particular, Replacement seems to be too weak in the context of urelements. Then, what is ZFC with urelements?

Three results suggest that ZFC with urelements should be formulated with \textit{Collection} instead. The first piece of evidence can be found in Theorem \ref{maintheorem1}: ZCU + Collection (since Collection trivially implies Replacement over ZU) yields desirable consequences such as the DC$_\omega$-scheme and the Reflection Principle.

Second, Collection is also essential for applying standard constructions to models of ZFC with urelements. Let $U$ be a model of $\ZFCUR$ and $F, x \in U$ be such that $U \models (F$ is an ultrafilter on $x)$. One can form an internal ultrapower of $U$ as usual. Namely, for every $f, g \in U$ such that $U \models $ ($f, g$ are functions on $x$), define
\begin{itemize}
    \item [] $f =_F g \text{ if and only if } U \models (\{y \in x: f(y) = g(y) \} \in F);$
    \item [] $[f] = \{h \in U : (h \text{ is a function on } x)^U \land h =_F f\};$
    \item [] $U/F = \{ [h] : h \in U \land (h \text{ is a function on } x)^U\}$.
\end{itemize}
For every $[f], [g] \in U/F$, define
\begin{itemize}
    \item [] $[g] \hat{\in} [f] \text{ if and only if } U \models (\{y \in x : g(y) \in f(y)\} \in F);$
    \item [] $\hat{\A}([f]) \text{ if and only if } U \models (\{y \in x : \mathcal{A}( f(y))\} \in F).$
\end{itemize}
Then the internal ultrapower is the model $\<U/F, \ \hat{\in}, \ \hat{\A}>$ (denoted by $U/F$). The \L o\'s theorem holds for $U/F$ if for every $\varphi$ and $[f_1], ..., [f_n] \in  U/F$, $U/F \models \varphi ([f_1], ..., [f_n])$ if and only if $U \models (\{y \in x : \varphi(f_1(y), ..., f_n(y))\} \in F).$
When $V \models $ ZFC, the \L o\'s theorem holds for all internal ultrapowers of $V$, which is commonly used in the study of large cardinals.

\begin{theorem}\label{thm:collection<->losthm}
Let $U$ be a model of $\ZFCUR$. The following are equivalent.
\begin{enumerate}
    \item The \L o\'s theorem holds for all internal ultrapowers of $U$.
    \item $U \models$ Collection.
\end{enumerate}
\end{theorem}
\begin{proof}
The proof of (2) $\rightarrow$ (1) is standard, and the point here is that the use of Collection is essential. 

For (1)$\rightarrow$(2), suppose that Collection fails in $U$. Then by Theorem \ref{maintheorem1}, it follows that both Plenitude and Tail fail in $U$. In $U$, fix some $A \subseteq \A$ without a tail cardinal and let $\kappa$ be the least cardinal not realized by any $B\subseteq \A$ that is disjoint from $A$, which is an infinite limit cardinal $U$. Let $F \in U$ be an ultrafilter on $\kappa$ containing all the unbounded subsets of $\kappa$. Suppose \textit{for reductio} that the \L o\'s theorem  holds for $U/F$. Let $id$ be the identity function on $\kappa$ and $c_A$ be the constant function sending every $\alpha < \kappa$ to $A$. Since $U \models (\{ \alpha < \kappa : \exists B \subseteq \A \ (B \sim \alpha \land B \cap A =\emptyset)\} \in F)$, by  the \L o\'s theorem, $U/F \models \exists B \subseteq \A (B \sim [id] \land B \cap [C_A] =\emptyset)$. Thus, there is some $g \in U$ such that 
$$U/F \models [g] \subseteq \A \land [g] \sim [id] \land ([g] \cap [C_A] =\emptyset).$$
Let $x= \{\alpha < \kappa : g(\alpha) \subseteq \A \land g(\alpha) \sim \alpha \land (g(\alpha) \cap A = \emptyset) \}$, which is in $F$ by the \L o\'s theorem again. Then $D = \bigcup_{\alpha\in x} g(\alpha)$ has size $\kappa$ and is disjoint from $A$---contradiction.
\end{proof}

Third, as we shall see in Chapter 3 (Theorem \ref{collection<->fullness}), over $\ZFCUR$ Collection is also equivalent to the principle that every (properly defined) forcing relation has the property of fullness, which is a property one would expect every forcing relation to have when AC is assumed. Hence, it is safe to say ZU + Collection + AC is a more robust theory than $\ZFCUR$. The following notation is thus justified, which has been adopted in \cite{HamkinsForthcoming-HAMRIS}. 
\begin{definition}
ZFCU = ZU + Collection + AC.
\end{definition}
\noindent However, $\ZFCUR$ (or $\ZFUR$) should not be discarded for two reasons that will be made clear. For one thing, $\ZFUR$ suffices for the basic forcing machinery and hence serves as a natural theory where one can study forcing with urelements. For another, since models of $\ZFCUR$ are easier to obtain, sometimes it is more convenient to start with a model of $\ZFCUR$ (or $\ZFUR$) and then establish Collection in the model (e.g., see the proof of Theorem \ref{thm:RPdoesnotproveHomogeneity}).

So far I have only offered \textit{extrinsic justifications} for Collection as an axiom, i.e., justifications based on its consequences. Can Collection be justified \textit{intrinsically} on the basis of a certain conception of set? Let us consider the three conceptions mentioned in Section \ref{section:UrelementsinSetTheory}. To begin with, it is unclear if the iterative conception of set is able to provide such justification: after all, it is a theorem of $\ZFUR$ that every set is in some $V_\alpha(A)$, in which case sets are indeed formed stage by stage. Regarding limitation of size, if we formulate it as a second-order axiom, Collection indeed follows (Proposition \ref{prop:GBU->RPandDCord}) because limitation of size implies that there is a global well-ordering. But it seems that a natural justification for Collection should not commit to any form of second-order choice principle. The reflection conception, however, provides a straightforward justification for Collection (Proposition \ref{RP->Collection}). Given that Collection is an attractive axiom, this, in turn, suggests that the reflection conception of set is more robust than the other two.

There is an alternative view regarding the question of what is ZFC with urelements. That is, in urelement set theory we turn out to have more ``axiomatic freedom'' in the sense that there are equally reasonable ways to axiomatize ZFC set theory with urelements even though they differ in strength; and it is this axiomatic freedom that prompts us to have a deeper understanding of the subject matter (see \cite{sep-set-theory-constructive} for a discussion on a similar view regarding intuitionistic set theory). A fact supporting this view is that even ZFCU has models that are somehow ``unnatural'': there can be models of ZFCU with a proper class of urelements where every set of urelements is only countable. This situation might conflict a standard conception of proper class, i.e., proper classes are \textit{big} in the sense that their sets are unbounded. Can there be a natural axiom securing this conception of proper class in urelement set theory? One can indeed formulate this conception as a second-order assertion called \textit{the Injection Principle} (see \cite[pp. 138-140]{jech2008axiom}), which says that every set can be injectively mapped into every proper class. $\ZFUR$ + Injection Principle proves AC, and under $\ZFUR$ + Injection Principle, either $\A$ is a set or Plenitude holds. Thus, $\ZFUR$ + Injection Principle has many desirable consequences by Theorem \ref{maintheorem1}. Yet the problem with Injection Principle is precisely that it is not neutral to AC. As a result, we cannot appeal to principles of this sort in a choiceless urelement set theory, to which I now turn.

\section{Urelement set theory without choice}\label{section:ChoicelesZFU}
The mutual interpretability between ZFC and ZFCU shown in Theorem \ref{con(zfc)->con(zfcu)} indicates certain redundancy of urelements when set theory is treated as a foundation: if every set of urelements is equinumerous with a pure set, then we may simply identify these urelements with objects in $V$. In other words, urelement set theory would become a more interesting foundational theory if sets of urelement are not necessarily well-orderable. Moreover, the assumption that every set of urelements, regardless what they are, is well-orderable seems to be rather restrictive as it excludes the existence of certain objects (mathematical or otherwise) \textit{a priori}. This provides motivations for studying urelement set theory in the absence of AC.

How should ZF with urelement be axiomatized? Regarding Diagram \ref{ZFCUdiagram}, it is natural to consider which implications still hold when AC is dropped. What further complicates this issue is the fact that different formulations of Plenitude and Tail come apart without choice.
\begin{itemize}
    \item [] (Plenitude) Every cardinal\footnote{In the choiceless context, by ``cardinals'' I always mean the \textit{well-ordered} cardinals---ordinals that are not equinumerous with any ordinal below themselves. Note that the general notion of cardinality, unlike in ZF, is not definable in $\ZFUR$, as shown in \cite{levy1969definability}.} is realized.
    \item [] (Plenitude$^+$) Every set $x$ is realized.
    \item [] (Tail) Every $A \subseteq \A$ has a tail.
    \item [] (Tail$^*$) For every $A \subseteq \A$, there is a greatest cardinal $\kappa$ such that $\exists B \subseteq \A \ (B \sim \kappa \land B \cap A = \emptyset)$.
     \item [] (Tail$^+$) Every $A \subseteq \A$ has a well-ordered tail.
\end{itemize}
I shall first show that the following implication diagram holds in $\ZFUR$.
\begin{figure}[hbt!]
\centering
\begin{tikzpicture}
\begin{scope}[every node/.style={}]
    \node (B) at (5, 1.5){RP};
    \node (C) at (-2,2.5) {Plenitude};
        \node (D) at (1,3) {Plenitude$^+$};

          \node (G) at (3, -0.5) {RP$^-$};

    \node (H) at (1, 1.5) {\textup{Closure$\land$Duplication}};
    \node (I) at (5, -0.5) {Collection};

    \node (O) at (1, -0.5) {Duplication};
      \node (F) at (-2, -0.5) {Closure};

    \node (T) at (5, 3) {$\A$ is a set};
    \node (R) at (7, 3) {Tail$^+$};
    \node (S) at (7, 1.5) {Tail};
      \node (U) at (9, 1.5) {Tail$^*$};
    
\end{scope}

\begin{scope}[>={stealth},
              every node/.style={fill=white,circle},
              every edge/.style={draw=black}]
 
    \path [->] (T) edge (B);
    \path [->] (T) edge (S);
    \path [->] (D) edge (C);
    \path [->] (D) edge (H);
     \path [->] (C) edge (F);
    \path [->] (B) edge (I);
    \path [->] (B) edge (G);

    \path [->] (R) edge (S);
\path [->] (R) edge (U);
    \path [->] (R) edge (B);
    \path [->] (H) edge (O);
    \path [->] (H) edge (F);
\end{scope}
\end{tikzpicture}
\caption{Implication diagram in $\ZFUR$}
\label{ZFUdiagram}
\end{figure}
\FloatBarrier 
\noindent It is still unknown if the diagram is complete in $\ZFUR$. I shall prove several independence results in the next subsection and summarize some key open questions at the end of this chapter. 

To show Diagram \ref{ZFUdiagram} holds in $\ZFUR$, we utilize the following theorem proved in \cite{HamkinsForthcoming-HAMRIS}.

\begin{theorem}[\cite{HamkinsForthcoming-HAMRIS}]\label{homo->RP}
$\ZFUR$ + Collection + $\ACA$ $\vdash$ RP. \qed
\end{theorem}
\noindent Note that if $\A$ is a set, then the usual proof of RP in ZF works; and RP implies Collection by Proposition \ref{RP->Collection}. So it remains to show the following.
\begin{theorem}\label{thm:Plenitude+->duplication}
Over $\ZFUR$,
\begin{enumerate}
    \item Tail$^+$ $\rightarrow$ RP;
    \item Plenitude$^+$ $\rightarrow$ Duplication.
\end{enumerate}
\end{theorem}
\begin{proof}
(1) Assume Tail$^+$. Consider a well-ordered tail for $\emptyset$. Then every set of urelements can be injectively mapped to this set of urelements, so $\ACA$ holds. Then it follows from Lemma \ref{tail->collection} that Collection holds, so we can apply Theorem \ref{homo->RP}.

(2) Assume Plenitude$^+$. Suppose \textit{for reductio} that some $A \subseteq \A$ cannot be duplicated. Consider any ordinal $\alpha$. Then there is a bijection $f$ from $A \times \alpha$ to some set of urelements $B$. It follows that for every $\beta < \alpha$, $A \cap f[A_\beta]$ is non-empty, where $A_\beta = A \times \{\beta\}$, which produces an injection from $\alpha$ to $P(A)$. This shows that every ordinal can be mapped injectively into $P(A)$, contradicting Hartog's Theorem.
\end{proof}

\subsection{Permutation models}
 Now I proceed to prove some independence results concerning Diagram \ref{ZFUdiagram} by constructing suitable models. Note that the $V\llbracket X \rrbracket$-construction in Definition \ref{barwiseinterpretation1} is not flexible enough for such task. Firstly, $V\llbracket X \rrbracket$ always satisfies Collection (Theorem \ref{thm:V[X]modelsZFU}), while we wish to show that, for instance, Plenitude does not imply Collection over $\ZFUR$. Secondly, it folllows from Theorem \ref{thm:ACholdsViffACholdsinV[X]} that if $\ACA$ fails in $V\llbracket X \rrbracket$, then AC fails for the pure sets of $V\llbracket X \rrbracket$. But we may wish to construct models where AC only fails \textit{outside} $V$. The method I shall utilize is a combination of the $U^\I$-construction in Definition \ref{smallkernelmodel} and the technique of \textit{permutation models} due to Fraenkel \cite{fraenkel1922begriff}, Mostowski \cite{mostowski1939begriff}, and Specker \cite{Specker1957-SPEZAD}. Since the standard text book on permutation models \cite{jech2008axiom} only considered permutation models with a \textit{set} of urelements, here I shall consider a more general construction allowing a proper class of urelements.
\begin{definition}[$\ZFCUR$]\label{permutationmodeldef}
Let $A$ be a set of urelements and $\G_A$ be a \textit{group} of permutations of $A$. For every $x$, define $sym(x) = \{\pi \in \G_A : \pi x = x\}$; if $x$ is a set, define $fix(x) = \{\pi \in \G_A : \pi y = y \text{ for all } y\in x \}$. A \textit{normal filter} $\F$ on $\G_A$ is a non-empty set of subgroups of $\G_A$ which contains $sym(a)$ for every urelement $a \in \A$ and is closed under supergroup, finite intersection, and conjugation (i.e., for all $\pi \in \G_A$ and $H \in \F$, $\pi H \pi^{-1} \in \F$). An object $x$ is \textit{symmetric} (with respect to $\F$) if $sym(x) \in \F$. The \textit{permutation model} $W$, determined by $A$, $\G_A$, and $\F$, is the class of all \textit{hereditarily symmetric} objects, i.e., $W = \{x \in U : x \text{ is symmetric} \land x \subseteq W\}$.
\end{definition}
\noindent If $x$ is symmetric, then so is $\pi x$ for every $\pi \in \G_A$ because $sym(\pi x) = \pi \circ sym(x) \circ \pi^{-1}$ and $\F$ is closed under conjugation. By an $\in$-induction, it follows that if $x \in W$, then $\pi x \in W$ for every $\pi \in \G_A$. Therefore, every $\pi \in \G_A$ is an automorphism of $W$.
\begin{theorem}[$\ZFCUR$]\label{fundamentalthmPM}
Let $A$, $\G_A$ and $\F$ model be as in Definition \ref{permutationmodeldef}. And let $W$ be the resultant permutation model. Then 
\begin{enumerate}
    \item $W \models \ZFUR$;
    \item $W \models$ Collection if $U \models$ Collection.
\end{enumerate}
\end{theorem}
\begin{proof}
(1) Since $W$ is transitive and contains all the pure sets, Extensionality, Foundation, and Infinity all hold in $W$, and AC holds for the pure sets of $W$. Union holds in $W$ because for any set $x \in W$, $sym(x) \subseteq sym(\bigcup x)$. If $x, y \in W$, $sym(x) \cap sym(y) \subseteq sym(\{x, y\})$, so $W$ satisfies Pairing.

\textit{Powerset.} Let $x \in W$ be a set. It suffices to show that $P^W(x) = \{y \in W : y \subseteq x\}$ is symmetric. If $\pi \in sym(x)$ and $y \in P^W(x)$, then $\pi y \subseteq x$ and $\pi y \in W$, so $\pi y \in P^W(x)$. This shows that $sym(x) \subseteq sym(P^W(x))$, and hence $P^W(x)$ is symmetric.

\textit{Separation.} Let $x \in W$ be a set. It suffices to show that the set $v = \{y \in W : y \in x \land \varphi^W(y, u)\}$ is symmetric, where $u$ is a parameter in $W$. If $\pi \in sym(x) \cap sym(v)$ and $y \in v$, it follows that $\pi y$ is in $W \cap x$ and $\varphi^W(\pi y,x,u)$ since $\pi$ is an automorphism of $W$. So $sym(x) \cap sym(v) \subseteq sym(v)$ and hence $v$ is symmetric.

\textit{Replacement.} Suppose that $W \models \forall x \in w \exists ! y \varphi(x, y, u)$, where $w, u \in W$. Let $v = \{y \in W : \exists x \in w \ \varphi^W(x, y, u)\}$, which is a set by Replacement in $U$. It suffices to show that $sym(w)\cap sym(u) \subseteq sym(v)$. If $\pi \in sym(w)\cap sym(u)$ and $y \in v$, then $\varphi^W(x, y, u)$ for some $x \in w$ and so $\varphi^W(\pi x, \pi y, u)$ for some $\pi x \in w$; thus, $\pi y \in v$ and hence $sym(w)\cap sym(u) \subseteq sym(v)$.

(2) Suppose that $U \models$ Collection and that $W \models \forall x \in w \exists y \varphi(x, y, u)$ for some $x, u \in W$. So in $U$ there is a $v$ such that $\forall x \in w \exists y \in v (y \in W \land \varphi^W(x, y, u))$. Let $B = A \cup ker(v)$. $B \in W$ because every $\pi \in \G_\A$ point-wise fixes all urelements outside $A$. Thus, $W \models \forall x \in w \exists y \in V(B) \ \varphi(x, y, u)$, and this suffices for Collection to hold in $W$ by Proposition \ref{weakcollection}. 
\end{proof}
\begin{definition}\label{Gnormalideal}
Given an $A\subseteq\A$ and a group $\G_A$ of permutations of $A$, $I \subseteq P(A)$ is a $\G_A$-\textit{normal ideal} on $A$ if and only if
\begin{enumerate}
    \item $A \notin I$;
    \item if $E_1, E_2 \in I$, then $E_1 \cup E_2 \in I$; 
    \item if $E_1 \in I$ and $E_2 \subseteq E_1$, then $E_2 \in I$;
    \item for every $a \in A$, $\{a\} \in I$;
    \item for every $E \in I$ and $\pi \in \G_A$, $\pi E \in I$.
\end{enumerate}
\end{definition}
\noindent If $I$ is a $\G_A$-normal ideal on $A$, it is not hard to verify that $$\F = \{ H \subseteq \G_A: H \text{ is a subgroup of }\G_A \text{ and } fix(E) \subseteq H \text{ for some } E \in I \}$$ is a normal filter on $\G_A$. Thus, in this case $A$, $\G_A$, and $I$ will generate a permutation model $W$. For every $x$ in such $ W$, there will be some $E \in I$, called a \textit{support} of $x$, such that $fix(E) \subseteq sym(x)$. This concludes the basic setup of permutation models. 

\begin{example}[The Basic Fraenkel Model]\label{exp:BasicFModel}
Let $U$ be a model of ZFCU in which $\A$ is a countably infinite set, $\G_\A \in U$ be the group of all permutations of $\A$, and $I$ be the ideal of finite subsets of $\A$. In the resulting permutation model $W$, although $\A$ is still the set of all urelements, no set of urelements is equinumerous with $\omega$ because any injection from $\omega$ to some $A\subseteq\A$ would have a finite support, which is impossible.
\end{example}
\begin{corollary}
Assume the consistency of ZF. There is a model of $\ZFUR$ in which 
\begin{enumerate}
    \item $\A$ is a set;
    \item Closure fails;
    \item Tail$^*$ fails.
\end{enumerate}
\end{corollary}
\begin{proof}
Consider the Basic Fraenkel Model.
\end{proof}

\subsection{Independence results}
Many implications in Diagram \ref{ZFCUdiagram} fail in $\ZFUR$. In particular, I shall prove that over $\ZFUR$,
\begin{enumerate}
    \item Plenitude $\nrightarrow$ (Duplication $\lor$ Collection);
    \item Tail$^*$ $\nrightarrow$ (Collection $\lor$ Tail);
    \item (Plenitude $\land$ Duplication) $\nrightarrow$ Collection;
    \item (RP $\land$ DC$_\omega$-scheme) $\nrightarrow$ Homogeneity holds over some $A\subseteq\A$.
\end{enumerate}
\begin{theorem}\label{plenitudenvdashcollection}
Assume the consistency of ZF. There is a model of $\ZFUR$ in which 
\begin{enumerate}
    \item Plenitude holds;
    \item Collection fails;
    \item Duplication fails.
\end{enumerate}
\end{theorem}
\begin{proof}
Let $U$ be a model of ZFCU + Plenitude. In $U$, fix a countable set of urelements $A \subseteq \A$ and enumerate it with $\omega \times \omega$. So $A = \bigcup_{n<\omega} A_n$, where each row $A_n$ is an infinitely countable sequence of urelements. Let $\G_A$ be the group of permutations of $A$ that preserve each $A_n$, i.e, a permutation $\pi$ of $A$ is in $\G_A$ just in case $\pi A_n = A_n$ for every $n<\omega$. Let $I = \{E \subseteq A : E \text{ is finite} \}$, which is a $\G_A$-normal ideal on $A$. Now let $W$ be the permutation model generated by $A$, $\G_A$ and $I$. By Theorem \ref{fundamentalthmPM}, $W \models \ZFUR$ + Collection.

To get the failure of Collection, we go to an inner model of $W$ by using the construction in Definition \ref{normalideal}. Since the sequence $\<A_n : n <\omega>$ is in $W$, say a $B\subseteq \A$ \textit{is finitely contained in} $A$ if $A \cap B$ is a subset of the union of finitely many $A_n$. Define $\I = \{B \subseteq \A : B  \text{ is finitely contained in } A \}$, which is an $\A$-ideal. This produces an inner model $W^\I = \{ x \in W : ker(x) \in \I\}$, and by Theorem \ref{smallkernelmodel}, $W^\I \models \ZFUR$. 

Note that if $B$ is a set of urelements in $W^\I$ that is disjoint from $A$, then $B$ is well-orderable in $W^\I$ because every $\pi \in G_A$ point-wise fixes $B$, so its well-ordering in $U$ is preserved through the constructions. It follows that  $W^\I \models$ Plenitude since in $U$ we can find arbitrarily large sets of urelements disjoint from $A$.
\begin{lemma}\label{amorphous}
In $W^I$, each $A_n$ is \textit{amorphous}, i.e., it is infinite but is not a union of two disjoint infinite sets.
\end{lemma}
\begin{proof}
Suppose for \textit{reductio} that $A_n = B_1 \cup B_2$ for two infinite disjoint sets $B_1$, $B_2$ in $W^\I$. Let $E_1 \in I$ be a support of $B_1$. Since both $B_1 \setminus E_1$ and $B_2 \setminus E_1$ are non-empty, we can pick an urelement from each of them and let $\pi$ be a permutation in $\G_A$ that only swaps these two urelements. It follows that $\pi B_1 \neq B_1$ and $\pi \in fix(E_1)$, contradicting the fact that $E_1$ supports $B_1$.\end{proof}

\begin{lemma}
$W^\I \models \neg$Collection.
\end{lemma}
\begin{proof}
The following holds in $W^\I$ since $W^\I$ contains $A_0 \cup ... \cup A_{n-1}$ for each $n$.
\begin{align}
    \forall n < \omega \ \exists D \subseteq\A \ (D \text{ is a union of } n \text{ disjoint amorphous sets}).
\end{align}
Suppose \textit{for reductio} that Collection holds in $W^\I$. Then there is a set $v \in W^\I$ such that
\begin{align}
    \forall n < \omega \ \exists D \in v \ (D \subseteq \A \land D \text{ is a union of } n \text{ disjoint amorphous sets}).
\end{align}
And $ker(v) \cap A$ is contained in an $m$-block of $A_n$, $A_{n_1} \cup ... \cup A_{n_m}$, for some finite number $m$. By (2.2), there is a set of urelements $D \in v$ such that $D = D_1 \cup ... \cup D_{m+1}$, where $D_1, ..., D_{m+1}$ are disjoint amorphous sets. For each $k \leq m+1$, $D_k \cap A$ must be infinite because any set of urelements disjoint from $A$ is well-orderable in $W^\I$. Since $D_k \cap A = D_k \cap (A_{n_1} \cup ... \cup A_{n_m})$, it follows that for each $k \leq m+1$, there is an $l \leq m$ such that $D_k \cap A_{n_l}$ is infinite. However, no two $D_k$ and $D_{k'}$ can have an infinite intersection with the same $A_{n_l}$ by Lemma \ref{amorphous}. This is a contradiction because it amounts to having an injection from $m+1$ to $m$.
\end{proof}
It remains to show that $W^\I \models$ $\neg$Duplication. It suffices to show that in $W^\I$, for each $A_n$ and any infinite set of urelements $B$, if $B \cap A_n$ is finite, then there is no injection from such $B$ to $A_n$. Suppose \textit{for reductio} that $f$ is an injection from $B$ to $A_n$ in $W^\I$. Let $E \in I$ be a support of $f$. Then there must be two urelements $a, b \in A_n \setminus (E\cup B)$. Let $\pi \in fix(E)$ swap only $a$ and $b$. It then follows that $f(b) = \pi f(a) = f(a)$, contradicting the injectivity of $f$. 
\end{proof}

\begin{theorem}
Assume the consistency of ZF. There is a model of $\ZFUR$ in which 
\begin{enumerate}
    \item Tail$^*$ holds;
    \item Collection fails;
    \item Tail fails.
\end{enumerate}
\end{theorem}
\begin{proof}
Let $U$ be a model of ZFCU + $\A \sim \omega_1$, where $\A$ is enumerated with $\omega \times \omega_1$. So $\A = \bigcup_{n < \omega} A_n$, where each row $A_n$ is uncountable. Let $G_\A$ be the group of permutations of $\A$ that preserve each row $A_n$ and $I = \{E \subseteq \mathcal{A} : E \text{ is countable} \}$. This generates a permutation model $W$. In $W$, define $\I = \{B \subseteq \mathcal{A} : B \text{ is finitely contained in } \A\}$.
\begin{lemma}\label{noomega1set}
In $W^\I$, no $A_n$ contains two disjoint subsets that are both uncountable. Hence, in $W^I$ there is no set of urelements of size $\omega_1$.
\end{lemma}
\begin{proof}
 Suppose \textit{for reductio} that in $W^\I$, $B_1, B_2 \subseteq A_n$ are disjoint and uncountable. Let $E \in I$ be a support of $B_1$. Then we can pick an urelement from each $B_1 \setminus E$ and $B_2 \setminus E$ respectively. A permutation that swaps only these two urelements will then fix $B_1$, which is a contradiction.
\end{proof}
\noindent Note that for each $n$, every countable subset of $A_n$ remains countable in $W^\I$. So for every set of urelements $B \in W^\I$, we can always find $\omega$-many urelements outside the finite block containing $B$. It follows from Lemma \ref{noomega1set} that Tail$^*$ holds in $W^\I$.

Suppose \textit{for reductio} that Collection holds in $W^\I$. Since for every $n<\omega$, there is a set of urelements that is a disjoint union of $n$ uncountable sets, it follows that there is some $v$ such that for every $n<\omega$, there is a set of urelements $D \in v$ that is a disjoint union of $n$ uncountable sets. Then $ker(v) \subseteq A_{n_1} \cup .... \cup A_{n_m}$ for some $m < \omega$. And there is some $D\subseteq A_{n_1} \cup  ....\cup A_{n_m}$ such that $D = D_1 \cup ... D_{m+1}$, where for each two $k, l \leq m+1$, $D_k$ and $D_l$ are disjoint and uncountable. For each $k \leq m+1$, there is a $l \leq m$ such that $D_k \cap A_{n_l}$ is uncountable; but no two $D_k$ and $D_{k'}$ can have an uncountable intersection with the same $A_{n_l}$ by Lemma \ref{noomega1set}. This is a contradiction, so Collection fails in $W^\I$. 

To see that Tail fails in $W^\I$, fix an $A_n$. For any $B$ that is disjoint from $A_n$, there is a row $A_m$ that is disjoint from $B \cup A_n$. But it is clear that there cannot be an injection in $W^\I$ from $A_m$ to $B$. \end{proof}

\begin{theorem}
Assume the consistency of ZF. There is a model of $\ZFUR$ in which 
\begin{enumerate}
    \item Plenitude holds;
    \item Duplication holds;
    \item Collection fails.
\end{enumerate}
\end{theorem}
\begin{proof}
Let $U$ be a model of ZFCU + Plenitude. In $U$, fix a set of urelements $A$ of size $\omega$ and enumerate it with $\omega \times \omega$, i,e, $A = \bigcup_{n<\omega} A_n$. We then identify each row $A_n$ with the rationals $\<\Q, <_\Q >$ such that $A_n = \{a^n_j : j \in \Q\}$. Define $\G_A$ as the following permutation group such that for every permutation $\pi$ of $A$,
\begin{itemize}
    \item [] $\pi \in \G_A$ if and only if there is an automorphism $\rho$ of $\<\Q, <_\Q>$ such that $\pi(a^n_j) = \pi (a^n_{\rho(j)})$ for every $n < \omega$ and $j \in \Q$.
\end{itemize}

\noindent That is, each $\pi \in \G_A$ follow a same automorphism of $\Q$ at each row $A_n$. Define $I = \{E \subseteq A : E \text{ is finite} \}$, which is a $\G_A$-normal ideal on $A$, and let $W$ be the permutation model generated by $A$, $\G_A$, and $I$. And as before, let $\I = \{B \subseteq \A : B \cap A \text{ is finitely contained in } A\}$. Plenitude holds in $W^\I$ because every set of urelements disjoint from $A$ is well-orderable in $W^\I$.
\begin{lemma}\label{Anareduplicates}
In $W^\I$, for every distinct $n, m <\omega$, $A_n$ and $A_m$ are duplicates.
\end{lemma}
\begin{proof}
Recall that two sets of urelements are said to be \textit{duplicates} if they are disjoint and equinumerous. It suffices to show that the bijection $f : A_n \rightarrow A_m$ that maps $a^n_i$ to $a^m_i$ for every $i \in \Q$ is symmetric. Consider any $\pi \in \G_\A$ and let $\rho$ be the automorphism of $\Q$ such that $\pi a^k_i = \pi a^k_{\rho i}$ for every $k < \omega$ and $i \in \Q$. For every $\<a^n_i, a^m_i> \in f$, $\<\pi a^n_i, \pi a^m_i > = \<a^n_{\rho i}, a^m_{\rho i}> \in f$. Therefore, $f$ is symmetric.
\end{proof}
\noindent For any set $B$ of urelements in $W^\I$, $B \setminus A$ can be easily duplicated outside $A$ because it is well-orderable and Plenitude holds. For $B \cap A$, since it is contained in an $m$-block of $A_n$ for some finite number $m$, we can find another disjoint $m$-block of $A_n$. By Lemma \ref{Anareduplicates}, it follows that these two blocks are duplicates, and in particular, $B \cap A$ will have a duplicate inside the disjoint block. This shows that Duplication holds in $W^\I$.

The failure of Collection in $W^\I$ will be proved through the following four lemmas.

\begin{lemma}\label{Anhasnoduplicates}
In $W^\I$, no $A_n$\footnote{Note that unlike in Theorem \ref{plenitudenvdashcollection}, $A_n$ is no longer amorphous: every $\pi \in \G_A$ that fixes $a^n_i$ will have to fix $\{a^n_j : j \leq i\}$ and $\{a^n_j : j > i\}$, making these two disjoint intervals symmetric.} contains a pair of infinite duplicates for any $n<\omega$. Hence, no infinite subset of $A_n$ is well-orderable. 
\end{lemma}
\begin{proof}
Suppose otherwise. Then in $W^\I$, there is some injection $f$ from $B$ to $B'$, where $B, B'$ are infinite duplicates and $B, B' \subseteq A_n$ for some $n$. Let $E \in I$ be a support of $f$. Since $E$ is finite, it follows that there must be some $a^n_i \in B \setminus E$ and $a^n_j \in B' \setminus E$ such that $f(a^n_i) = a^n_j$. Then we can find an open interval of $\Q$ that contains $j$ but no rational indexes that appeared in $E  \cup \{a^n_i\}$. Any non-trivial automorphism of $\Q$ that only moves points in this interval will generate a $\pi \in G_\A$. Clearly, $\pi \in fix(E \cup \{a^n_i\})$; so $f(a^n_i) = \pi a^n_j \neq a^n_j$, which is a contradiction.
\end{proof}
For any urelement $a^n_i \in A_n$, we say that a function $f$ \textit{vertically fixes} $a^n_i$ if $f (a^n_i) = a^m_i$ for some $m$.
\begin{lemma}\label{alignment}
 In $W^\I$, if $B\subseteq A_n$ and $B' \subseteq A_m$ are infinite, where $n \neq m$, then for every injection $f$ from $B$ to $B'$, $f$ vertically fixes infinitely many $a^n_i \in B$.
\end{lemma}
\begin{proof}
Let $E \in I$ be a support of $f$. We show that for every $a^n_i \in B \setminus E$, $f(a^n_i) = a^m_i$. Suppose \textit{for reductio} that $f(a^n_i) = a^m_j$ and $i \neq j$. Then we can find an interval of $\Q$ that contains $i$ but no rational indexes that have appeared in $E \cup \{a^m_j\}$, which will give us a $\pi \in fix(E \cup \{a^m_j\})$. But then $f(\pi a^n_i) = a^m_j$ and $\pi a^n_i \neq a^n_i$, contradicting the assumption that $f$ is injective.\end{proof}

For a set of urelements $B$, we say that it has a \textit{nice n-partition} if  there are $B_1, ... B_n$ such that 
\begin{itemize}
    \item [] (i) $B = B_1 \cup ... \cup B_n$;
    \item [] (ii) for each two $k, l \leq n$, $B_k, B_l$ are non-well-orderable duplicates.
\end{itemize}
\begin{lemma}\label{partition}
In $W^\I$, there is no set of urelements $B$ such that for every $n < \omega$, $B$ has a subset with a nice $n$-partition.
\end{lemma}
\begin{proof}
Suppose \textit{for reductio} that $B \in W^\I$ is such set. Then we may assume that $B \cap A \subseteq A_1, ..., A_n$ for some finite number $n$. Let $C \subseteq B$ be a set with a nice $n+1$-partition such that $C = C_1 \cup ... \cup C_{n+1}$, where for each two $k, l \leq n$, $C_k, C_l$ are non-well-orderable duplicates. And for each $1 \leq k < n+1$, let $f_k \in W^\I$ be a bijection from $C_k$ to $C_{k+1}$. Define an $(n+1)$-sequence of pairs $\<m_{1}, D_1, >, ..., \<m_{n+1}, D_{n+1}>$ by recursion as follows.
\begin{itemize}
   \item [] $m_1$ is the least number $\leq n$ such that $C_1 \cap A_{m_1}$ is infinite; and
    \item [] $D_1 = \{a \in C_1 \cap A_{m_1} : f_1  \text{vertically fixes } a \}$.
\end{itemize}   
Suppose that $\<m_{k-1}, D_{k-1}>$ has been defined. Then
 \begin{itemize}
    \item [] $m_{k}$ is the least number $\leq n$ such that $f_{k-1}[D_{k-1}] \cap A_{m_k}$ is infinite; and
    \item [] $D_{k} = \{ a \in f_{k-1} [D_{k-1}] \cap A_{m_k}: f_k \text{ vertically fixes } a\}$.
\end{itemize}
\begin{claim}
For each $k \leq n+1$, $m_k$ exists and $D_k \subseteq C_k$ is infinite. Hence, the sequence above is well-defined.
\end{claim}
\begin{claimproof}
First,  $m_1$ exists. For $C_1$ is non-well-orderable so $C_1 \cap A$ must be infinite, and since $C_1 \cap A \subseteq B \cap A \subseteq A_1 \cup ... \cup A_n$, it follows that $C_1 \cap A_m$ is infinite for some $m \leq n$ and hence $m_1$ exists. Second, observe that $D_1$ is infinite. $f_1[C_1 \cap A_{m_1}]$ is non-well-orderable because $C_1 \cap A_{m_1}$ is by Lemma \ref{Anhasnoduplicates}; so infinitely many urelements in $f_1[C_1 \cap A_{m_1}]$, which is a subset of $C_2$, live in $A$ and hence in some $A_l$, where $l \leq n$. Since $C_1$ and $C_2$ are disjoint, $l \neq m_1$ by Lemma \ref{Anhasnoduplicates}. This means that $f_1$ moves infinitely many urelements in $C_1 \cap A_{m_1}$ to another row $A_l$, and by Lemma \ref{alignment}, $f_1$ must vertically fix infinitely many urelements in  $C_1 \cap A_{m_1}$. Therefore, $D_1$ is infinite.

Now suppose that $\<m_{k-1}, D_{k-1}>$ exists and $D_{k-1} \subseteq C_{k-1}$ is infinite. Then $f_{k-1}[D_{k-1}]$, which is a subset of $C_{k} \cap A$, must have an infinite intersection with some $A_l$ for $l \leq n$. Hence $m_k$ exists. By the same reasoning as in the last paragraph, it follows that $f_k$ must vertically fix infinitely many urelements in $f_{k-1}[D_{k-1}] \cap A_{m_k}$, and hence $D_k$ is infinite.
\end{claimproof}

\begin{claim}
For every natural number $l > 0$, whenever $k, k+l \leq n+1$ for some natural number $k$, then for every $a^{m_{k+l}}_i \in D_{k+l}$, there is some $a^{m_k}_i \in D_k$ such that $a^{m_{k+l}}_i = f_{k+l-1} \circ ...\circ f_k(a^{m_k}_i)$.
\end{claim}
\begin{claimproof}
By an easy induction on $l$. When $l = 1$, the claim follows by the definition of $D_k$. Suppose it holds for $l - 1 (> 0)$ and consider $k, k+l \leq n+1$. For every $a^{m_{k+l}}_i \in D_{k+l}$, there is an $a^{m_{k+l-1}}_i \in D_{k+l-1}$ such that $f_{k+l-1}(a^{m_{k+l-1}}_i) = a^{m_{k+l}}_i$. By the induction hypothesis, there is an  $a^{m_k}_i \in D_k$ such that $a^{m_{k+l-1}}_i = f_{k+l-2} \circ ...\circ f_k(a^{m_k}_i)$, and so the claim follows.\end{claimproof}

It then follows from the claim that if $k < k' \leq n+1$, then $m_k \neq m_k'$. Otherwise, by the previous claim we can find a $a^{m_{k}}$ in both $D_k$ and $D_k'$, which are disjoint. However, this is a contradiction, because for each $k \leq n+1$, $m_k \leq n$. The lemma is thus proved. \end{proof}

Now suppose \textit{for reductio} that Collection holds in $W^\I$. In $W^\I$, we have 
\begin{align}
     \forall m < \omega \exists B \subseteq \A (B \text{ has a nice } m \text{-partition}).
\end{align}
This is because the union of any $m$-block of $A_n$ has has a nice $m$-partition, since each pair of $A_n$ are duplicates by Lemma \ref{Anareduplicates}. Then by Collection, 
\begin{align}
   \exists v \forall m < \omega \exists B \in v (B \subseteq \A \land B \text{\ has a nice } m\text{-partition}).
\end{align}
This means that for every $m$, $ker(v)$ has a subset with a nice $m$-partition, which contradicts Lemma \ref{partition}. This completes the proof of the theorem. \end{proof}

The next theorem explains why all the $\ZFCUR$ arguments fail without choice: homogeneity, which is crucial for those arguments to go through, may fail globally when $\ACA$ is not assumed.
\begin{theorem}\label{thm:RPdoesnotproveHomogeneity}
Assume the consistency of ZF. There is a model of $\ZFUR$ in which
\begin{enumerate}
    \item the DC$_\omega$-scheme holds;
    \item RP holds;
    \item homogeneity holds over no set of urelements.
\end{enumerate}
\end{theorem}
\begin{proof}
Let $U$ be a model of ZFCU in which $\A$ is a set of size $\omega_1$  and enumerate $\A$ with $\omega_1 \times \omega_1$, i.e., $\A = \bigcup_{\alpha < \omega_1} A_\alpha$, where each row $A_\alpha$ is uncountable. Let $G_\A$ be the group of permutations of $\A$ that preserve $A_\alpha$ for each $\alpha < \omega_1$ and $I = \{E \subseteq \mathcal{A} : E \cap A_\alpha \text{ is countable for each } \alpha < \omega_1\}$. This generates a permutation model $W$. In $W$, let $\I = \{B \subseteq \mathcal{A} : B \text{ is countably contained in } \A\}$.

As before, it follows that in $W^\I$, for any $\alpha, \beta <\omega_1$ such that $\alpha \neq \beta$, there is no injection from $A_\alpha$ to $A_\beta$. The next lemma says that every permutation of $\A$ that only swaps some $A_\alpha$ fixes $W^\I$.

 \begin{lemma}\label{Wfix}
In $U$, if $\sigma$ is a permutation of $\mathcal{A}$ such that for every $\alpha < \omega_1$, $\sigma A_{\alpha} = A_{\beta}$ for some $\beta$ (consequently,  $A_{\alpha} = \sigma A_{\gamma}$ from some $\gamma$), then $\sigma (W^\I) = W^\I$. 
\end{lemma}
\begin{proof}
Since $W^\I$ is a class of $U$ defined by $\mathcal{G}_\A, \ I$, and $\I$, it suffices to show that they are all fixed by $\sigma$. If $\pi \in \G_\A$, then for every $A_\alpha$, since $A_\alpha = \sigma A_\beta$ for some $\beta$ and $\pi A_\beta = A_\beta$, by automorphism it follows that $(\sigma \pi) (\sigma A_\beta) = \sigma A_\beta$; so $(\sigma \pi )A_\alpha = A_\alpha$ and hence $\sigma \pi \in \mathcal{G}_\A$.\footnote{Note here $\sigma \pi$ is not $\sigma \circ \pi$ but $\{\langle \sigma a, \sigma (\pi a) \rangle : a \in \mathcal{A} \rangle\}$} This shows that $\sigma \mathcal{G}_\A = \mathcal{G}_\A$. If $E \in I$, then for every $A_\alpha$, since $A_\alpha = \sigma A_\beta$ for some $\beta$ and $E \cap A_\beta$ is countable, $\sigma E \cap A_\alpha$ is countable and hence $\sigma E \in I$. Therefore, $\sigma I= I$. Similarly, if $B$ is contained in countably many $A_\alpha$, then so is $\sigma B $. Hence, $\sigma \I= \I$ and the lemma is proved.
\end{proof}
\begin{lemma}
$W^\I \models $ DC$_\omega$-scheme.
\end{lemma}
\begin{proof}
Suppose that in $W^\I$, $\varphi(x, y, u)$ defines a relation without terminal nodes. By the DC$_\omega$-scheme in $U$, there is an infinite sequence $\<x_n : n<\omega>$ in $U$ such that $x_n \in W^\I$ and $W^\I \models \varphi(x_n, x_{n+1}, u)$ for every $n$. By AC in $U$, for each $n$ we can choose an $E_n$ which is a support of $x_n$. Then $\bigcup_{n < \omega} E_n$ is in $I$, and as $\bigcup_{n < \omega} E_n $ supports $\langle x_n : n < \omega \rangle$, the sequence is also in $W$. Furthermore, since each $ker(x_n)$ is in $\I$ and $\I$ is countably closed, it follows that the kernel of this sequence is also in $\I$. Therefore, the  DC$_\omega$-scheme holds in $W^\I$. 
\end{proof}

\begin{lemma}
$W^\I \models $ Collection.
\end{lemma}
\begin{proof}
Suppose that $W^\I\models \forall x \in w \exists y \varphi(x, y, u)$ for some $w, u \in W^\I$. Let $\{ A^1_{n} : n < \omega \}$ be a countable block containing $ker(w) \cup ker(u)$ and $\{ A^2_n : n < \omega \}$ be a countable block that is disjoint from each $A^1_{n}$. Set $A = \bigcup_{n<\omega} (A^1_n \cup A^2_n )$. It suffices to show that $W^\I\models \forall x \in w \exists y \in V(A) \varphi(x, y, u)$. Suppose that $x \in w$ and fix some $y \in W^\I$ such that $W^\I \models \varphi(x, y, u)$. Let $\{ A^3_n : n <\omega\}$ be another disjoint countable block that contains $ker(y) \setminus A$. Back in $U$, find another block $\{A^4_n : n <\omega \}$ disjoint from all of these three blocks. Using AC in $U$ we can define a permutation $\sigma$ of $\mathcal{A}$ as follows. Let $\sigma$ move each $A^2_n$ to $A^2_{2n}$, each $A^3_n$ to $A^2_{2n+1}$, each $A^4_{2n}$ to $A^3_n$ and each $A^4_{2n+1}$ to $A^4_n$. By Lemma \ref{Wfix}, it follows that $\sigma (W^\I) = W^\I$. Since $ W^\I \models \varphi(x, y, u)$ and $\sigma$ fixes $x$ and $u$, it follows that $W^\I \models \varphi(x, \sigma y, u) \land \sigma y \in  V(A)$.\end{proof}
By Theorem \ref{collection+dc->rp}, the last two lemmas jointly imply that $W^\I \models $ RP. Finally, for every $A \subseteq \A$ in $W^\I$, there will be some $A_\alpha$ and $A_\beta$ such that $A$, $A_\alpha,$ and $A_\beta$ are pair-wise disjoint. But there is no injection from $A_\alpha$ to $A_\beta$. Therefore, homogeneity holds over no set of urelements in $W^\I$, which completes the proof.\end{proof}

\subsection{Open questions}
\begin{question}
Is Diagram \ref{ZFUdiagram} complete over $\ZFUR$? In particular, does $\ZFUR$ prove any of the following?
\begin{enumerate}
    \item Plenitude$^+ \rightarrow$ Collection.
    \item Tail $\rightarrow$ Collection.
    \item Collection $\rightarrow$ RP$^-$.
    \item RP$^- \rightarrow$ Collection.
\end{enumerate}
\end{question}
\noindent My conjecture is that $\ZFUR$ + Collection cannot prove RP$^-$ (and therefore RP). This conjecture has several implications from a philosophical perspective. Some may view this result as another compelling evidence in support of the reflection conception of sets, as it highlights the limitations of the iterative conception in the choiceless context. Thus, the iterative conception may be seen as an incomplete picture of sets. Conversely, the independence of RP could lead others to question the reflection conception in the context of urelements. And the fact that ZFCU proves RP suggests that RP is a form of choice principle, which may not be natural in urelement set theory. Furthermore, an alternative perspective is that these independence results demonstrate another instance of axiomatic freedom, emphasizing that there may not necessarily be a ``correct" system. 
\chapter{Forcing with Urelements}
In this chapter, I investigate forcing in the context of urelement set theory. In Section \ref{section:twomethods}, I review two standard methods of forcing: forcing via countable transitive models and forcing via Boolean-valued models. Section \ref{section:forcingoverctm} explores forcing over countable transitive models of $\ZFUR$. To overcome a significant drawback of the existing approach regarding the fullness property, I develop a new forcing machinery. The main results are as follows: (i) Over $\ZFCUR$, Collection is equivalent to the principle that every forcing relation, defined in the new way, is full. (ii) Forcing over $\ZFUR$ preserves $\ZFUR$ together with several axioms introduced in Section \ref{section:additionalaxioms}. (iii) Forcing can also destroy the $\DCK$ and recover Collection. (iv) Ground model definability fails when the ground model contains a proper class of urelements. (v) The new forcing machinery generates the same forcing extensions as the old one. In Section \ref{section:forcingBVM}, I provide a brief overview of some fundamental results about Boolean-valued models with urelements, established in joint work with Wu \cite{wu2022}. Based on these results, I consider how Boolean ultrapowers can be applied to arbitrary models of ZFCU.

\section{Two methods}\label{section:twomethods}
Forcing is a widely used technique in contemporary set theory. It has led to landmark results such as the independence of the Continuum Hypothesis and philosophical analyses such as the multiverse conception of set (\cite{Hamkins2012-HAMTSM} and \cite{Antos2015-ANTMCI}) and set-theoretic potentialism \cite{Hamkins2022-HAMTML}.

In the presence of urelements, it is natural to ask how forcing behaves. For example, will a forcing extension of a model of $\ZFUR$ always be a model of $\ZFUR$? Which of the axioms isolated in Chapter 2 will be preserved by forcing? Will the presence of urelements affect the forcing machinery? These questions are the focus of this chapter.

There are two main approaches to forcing. The first is the countable transitive model approach (CTM), which involves starting with a countable transitive model $M$ of set theory and a forcing poset $\P$ in $M$, and then extending $M$ by an $M$-generic filter $G$ over $\P$. The resulting forcing extension $M[G]$ is a countable transitive model of set theory where various set-theoretic statements, such as the Continuum Hypothesis, may fail or hold depending on the choice of the poset $\P$.

The second approach is based on Boolean-valued models (BVM). Given a complete Boolean algebra $\B$, a Boolean-valued model $M^\B$ for a first-order language $\mathcal{L}$ consists of a domain of $\B$-names together with a $\B$-valued truth assignment $\llbracket \ \rrbracket_\B$, which assigns a $\B$-value to each assertion in $\mathcal{L}$ about the $\B$-names in a way that obeys the axioms of first-order logic. If $V$ is an arbitrary model of set theory, we can form $V^\B$ as a definable class in $V$ by a complete Boolean algebra $\B$ in $V$. By the fundamental theorem of $\VB$, all the axioms of set theory will have value $1$ in $V^\B$ for every $\B$. This allows us to prove the independence of a statement $\varphi$ from set theory by finding some $\B$ such that $\llbracket \varphi \rrbracket_\B \neq 1$.

The CTM approach to forcing assumes the existence of a countable transitive model of ZFC, which is a stronger assumption than the mere consistency of ZFC. This drawback is usually addressed by working with a countable transitive model of a finite fragment of ZFC and appealing to the reflection principle in the meta-theory. On the other hand, the BVM approach to forcing does not require the ground model to be either transitive or countable. From any model of ZFC, one can construct a concrete model of ZFC + $\varphi$ by taking the quotient structure of certain Boolean-valued models. The fullness of $V^\B$ ensures that the Łoś Theorem holds, and the quotient structure can be seen as a definable class in $V$. As a result, the BVM approach is able to establish mutual interpretability between ZFC (see \cite{FreireForthcoming-FREBIW-2} for more on this) and its various extensions, while the CTM approach only establishes equiconsistency. Furthermore, the BVM approach provides a naturalistic account of forcing and allows one to force over any model of ZFC, avoiding the need for countable transitive models (see \cite{hamkins2012well}).

However, in the context of forcing with urelements, the theoretic virtues of the BVM approach are not as clear as in the classical case. While $\ZFUR$ is sufficient for the basic machinery of forcing, if the ground model $U$ is only a model of $\ZFUR$, the quotient structure may fail to be constructed due to the essential role of AC. For instance, if AC fails in $U$, there may not be any non-principal ultrafilters on a complete Boolean algebra $\B \in U$. Moreover, even if such an ultrafilter $F$ exists, the Boolean-valued model $U^\B$ may not be full without AC, and therefore the Łoś theorem may not hold for $U^\B/F$. In fact, even if $U$ is a model of $\ZFCUR$, $U^\B$ may still fail to be full, as the fullness of every properly defined $U^\B$ is equivalent to Collection (see \cite{wu2022}). Thus, a naturalist account of forcing with urelements is only possible when the ground model satisfies ZFCU. Consequently, the CTM approach becomes useful when we are interested in forcing over $\ZFUR$. An investigation of both approaches of forcing in the context of urelements is thus justified.

\section{Forcing over countable transitive models of $\ZFUR$}\label{section:forcingoverctm}
In this section, I will investigate poset forcing over countable transitive models of $\ZFUR$. Basic knowledge of forcing in ZFC covered in \cite[Ch.\ VII]{kunen2014set} will be assumed. Our meta theory, accordingly, will be some suitable urelement set theory such as $\ZFUR$. Notably, forcing with urelements has been studied in several places including \cite{blass1989freyd}, \cite{hall2002characterization}, and \cite{Hall2007-ERIPMA}. However, in all of these studies it is assumed that the urelements form a set,\footnote{For instance, in \cite{hall2002characterization}, Hall shows that if $N \subseteq M$ are countable transitive model of $\ZFUR$ with a set of urelements, then $N$ is a permutation model of $M$ only if $M$ is a certain forcing extension of $N$.} so it only needs some trivial adjustments to show that forcing preserves the axioms. But when a proper class of urelements is allowed, it becomes interesting to see which of the axioms introduced in Section \ref{section:additionalaxioms} are preserved by forcing. And this, as we shall see, will require some new arguments based on the earlier results.

\subsection{The existing approach and its problem}
 In pure set theory, given a forcing poset $\P$ with the maximal element $1_\P$, by transfinite recursion we define: $\dotx $ is a $\P$-name if and only if $\dotx$ is a set of ordered-pairs $\<\doty, p>$, where $\doty$ is a $\P$-name and $p \in \P$. Then every set $x$ in $V$ will have a canonical name  $\check{x} = \{ \<\check{y}, 1_\P> : y \in x \}$. In particular, $\emptyset$ be the its own name. To generalize this definition in urelement set theory, a natural idea, adopted in  \cite{blass1989freyd}, \cite{hall2002characterization} and \cite{Hall2007-ERIPMA}, is to treat each urelement as a different copy of $\emptyset$, which yields the following definition.
\begin{definition}\label{oldpnames}
Let $\P$ be a forcing poset. $\dot{x}$ is a $\P$-name$_\#$ if and only if either $\dot{x}$ is an urelement, or $\dot{x}$ is a set of ordered-pairs $\langle \dot{y}, p \rangle$, where $\dot{y}$ is a $\P$-name$_\#$ and $p \in \P$. $U^{\P}_\# = \{ \dot{x} \in U : \dotx \text{ is a } \P\text{-name}_\#\}$.
\end{definition}
\noindent This definition turns out to have a drawback, and the subscript $\#$ is meant to indicate this fact. To reveal its problem, we need to develop the basics of this approach. Let $\mathcal{L}^\P_\#$ be forcing language which contains $\{=, \in, \A \}$ as the non-logical symbols and each $\P$-name$_\#$ as a constant symbol. For each formula $\varphi (v_1, ..., v_n) \in \mathcal{L}^\P_\#$ and $\dotx_1, ..., \dotx_n \in U^{\P}_\#$, one can define the forcing relation $p \forces^\P_\# \varphi(\dotx_1, ..., \dotx_n)$ by recursion as follows (the superscript $\P$ will be omitted when it is clear from the context).

\begin{definition}\label{oldforcinglrelation}
Let $\P$ be a forcing poset. For every $\dotx, \doty, ... \in U^{\P}_\#$ and $p \in \P$,
\begin{enumerate}
    \item $p \forces_\# \dotx \in \doty$ if and only if $\{ q \in \P : \exists \<\dotz, r> \in \doty (q \leq r \land q \forces_\# \dotz = \dotx)\}$ is dense below $p$.
    \item $p \forces_\# \dotx \subseteq \doty$ if and only if whenever $\<\dotz, r> \in \dotx$ and $q \leq p, r$, then $q \forces_\# \dotz \in \doty$.
    \item $p \forces_\# \dotx = \doty$ if and only if either (i) $\dotx$ and $\doty$ are the same urelement, or (ii) $p \forces_\# \dotx \subseteq \doty$ and $p \forces_\# \doty \subseteq \dotx$.
    \item $p \forces_\# \A (\dotx)$ if and only if $\dotx$ is an urelement.
    \item $p \forces_\# \neg \varphi$ if and only if there is no $q \leq p$ such that $q \forces_\# \varphi$.
    \item $p \forces_\# \varphi \land \psi$ if and only if $p \forces_\# \varphi$ and $p \forces_\# \psi$.
    \item $p \forces_\# \exists x \varphi$ if and only if $\{q \in \P : \text{ there is some } \dotz \in U^{\P}_\# \text{ such that }q \forces_\# \varphi(\dotz)\}$ is dense below $p$.
\end{enumerate}
\end{definition}

Now let $M$ be a countable transitive model of $\ZFUR$ and $\P \in M$ be a forcing poset. Given an $M$-generic filter $G$, for every $\dotx \in M^\P_\# = M \cap U^\P_\#$ , we define 

\begin{equation*}
       \dotx_G =
    \begin{cases*}
       \dotx & if $\A(x)$ \\
\{\dot{y}_G : \exists p \in G \langle \dot{y}, p \rangle \in \dot{x} \}     & otherwise 
    \end{cases*}
\end{equation*}
$M[G]_\# = \{\dotx_G : \dotx \in M^\P_\#  \}$ is then a transitive model that includes $M$ with the same ordinals and urelements of $M$. Moreover, one can easily prove the forcing theorem for $(\forces_\#)^M$.

\begin{theorem}[The Forcing Theorem for $\forces_\# $]
Let $M$ be a countable transitive model of $\ZFUR$, $\P \in M$ be a forcing poset, and $G$  be an $M$-generic filter over $\P$. For every $\dotx_0, ..., \dotx_n \in M^\P_\#$, 
\begin{align*}
    M[G]_\# \models \varphi (\dotx_{0_{G}}, ..., \dotx_{n_{G}})\text{ if and only if } \exists p \in G  (p \forces_\# \varphi (\dot{x}_0, ..., \dot{x}_n))^M.
\end{align*}
\end{theorem}
\begin{proof}
By the definition of $\forces_\#$, $M[G]_\# \models \A (\dotx_{G})$ just in case  $p \forces_\#  \A(\dotx)$ for every $p \in \P$, so the urelement predicate causes no problem. And when one of $\dotx_0$ and $\dotx_1$ is an urelement, their $G$-valuations are identical only if $\dotx_0$ and $\dotx_1$ are the same urelement. The rest of the theorem can then be proved by standard text-book arguments as in Kunen\cite[Ch.VII]{kunen2014set}.
\end{proof}
\noindent In fact, one can proceed to show that  $M[G]_\# \models \ZFUR$. 

However, one important feature of forcing is missing in this approach. The following is a standard theorem of ZFC (i.e., $\ZFCUR$ + ``there is no urelements'').

\begin{theorem}[ZFC]
Let $\P$ be a forcing poset. Its forcing relation is \textit{full}, i.e., whenever $p \forces \exists y \varphi(y,\dotx_0, ..., \dotx_n)$, then $p \forces \varphi (\doty, \dotx_0, ..., \dotx_n)$ for some $\doty \in V^\P$. \qed
\end{theorem}
\begin{remark}
If $U$ contains two urelements and $\P$ contains a maximal antichain with at least two elements, its forcing relation $\forces_\#$ is not full.
\end{remark}
\begin{proof}
Suppose that $\<p_i : i \in I>$, where $I$ has at least two elements, is a maximal antichain, and let $\<a_i : i \in I>$ be some urelements such that at least two of them are distinct. Consider the $\P$-name$_\#$ $\dotx = \{\<a_i, p_i> : i \in I \}$. It is routine to check that $1_\P \forces_\# \exists y (y \in \dotx)$. Suppose \textit{for reductio} that $1_\P \forces_\# \doty \in \dotx$ for some $\doty \in U^\P_\#$. Then there will be two distinct urelements $a_i$ and $a_j$ such that both ``$\doty = a_i$'' and ``$\doty = a_j$'' are forced, which is impossible by the definition of $\forces_\#$.
\end{proof}
\noindent The reason of why this happened is that $U^\P_\#$ contains too few names. Recall the following standard theorem in ZF.
\begin{theorem}[ZF]
Let $\P$ be a forcing poset. Then if $f$ is a function from an antichain of a forcing poset $\P$ to $V^\P$, then there is a $\doty \in V^\P$, called \textit{a} \textit{mixture} of $f$, such that $p \forces f(p) = \doty$ for every $p \in dom(f)$. \qed
\end{theorem}
\noindent As we have seen, this does not hold for $U^\P_\#$ because we cannot even mix two urelements.
\subsection{A new forcing machinery with urelements}
To have mixtures, we would want all names to be sets of ordered pairs. Also, in any forcing extension, if a name is collapsed into an urelement, there must be a unique one; furthermore, no name should be collapsed into a member of this urelement. This motivates the following new definition of names.
\begin{definition}\label{newpnames}
Let $\P$ be a forcing poset. $\dot{x}$ is a $\P$-name if and only if (i) $\dot{x}$ is a set of ordered-pairs $\langle y, p \rangle$ where $p\in \P$ and $y$ is either a $\P$-name or an urelement, and (ii) whenever $\langle a, p \rangle, \langle y, q\rangle \in \dot{x}$, where $a$ is an urelement and $a \neq y$, then $p$ and $q$ are incompatible (written as $p \bot q$). For every urelement $a$, $\check{a} = \{\<a, 1_\P>\}$; for every set $x$, $\check{x} = \{\<\check{y}, 1_\P> : y \in x\}$. $U^\P = \{ \dotx : \dotx \text{ is a } \P \text{-name}\}$.
\end{definition}
\noindent Note that $U^\P$, unlike $U^\P_\#$, contains no urelements. In particular, $\{\<a, 1>\}$ is not the canonical name of $\{a\}$ but the canonical name of $a$ itself. And when $\<a, p> \in \dotx$ for some urelement $a$, this indicates that $a$ will be \textit{identical to}, rather than \textit{a member of}, $\dotx_G$ for any generic filter $G$ containing $p$. Now we proceed to define the new forcing relation.

\begin{definition}\label{def:newforcingrelation}
Let $\P$ be a forcing poset. The forcing language $\mathcal{L}^\P$ contains $\{\subseteq, =, \in, \A, \overset{\mathscr{A}}{=}\}$ as the non-logical symbols and every $\P$-name in $U^\P$ as a constant symbol. For every $\dotx_1, \dotx_2, ... \in U^\P$, $p \in P$ and $\varphi \in \mathcal{L}^\P$,
\begin{enumerate}
    \item $p \forces \A(\dotx_1)$ if and only if $\{q \in \P : \exists \<a , r> \in \dotx_1 \ (\A(a) \land q \leq r)\}$ is dense below $p$.
    \item $p \forces \dotx_1 \Aeq \dotx_2$ if and only if $\{q \in \P : \exists a, r_1, r_2 (\A(a) \land \<a, r_1> \in \dotx_1 \land \<a, r_2> \in \dotx_2 \land q \leq r_1, r_2 )\} \cup \{q \in \P : \forall \<a_1, r_1> \in \dotx_1 \ (\A(a_1) \rightarrow q \bot r_1) \land \forall \<a_2, r_2> \in \dotx_2 \ (\A(a_2) \rightarrow q \bot r_2)\}$ is dense below $p$. 
    \item $p \forces \dotx_1 \in \dotx_2$ if and only if $\{ q \in \P : \exists \<\doty, r> \in \dotx_2 (\doty \in U^\P \land q \leq r \land q \forces \doty = \dotx_1)\}$ is dense below $p$.
    \item $p \forces \dotx_1 \subseteq \dotx_2$ if and only if for every $\doty \in U^\P$ and $r, q \in \P$, if $\<\doty, r> \in \dotx_1$ and $q \leq p, r$, then $q \forces \doty \in \dotx_2$.
    \item $p \forces \dotx_1 = \dotx_2 $ if and only if $p \forces \dotx_1 \subseteq \dotx_2$, $p \forces \dotx_2 \subseteq \dotx_1$ and $p \forces \dotx_1 \Aeq \dotx_2$.
    \item $p \forces \neg \varphi$ if and only if there is no $q \leq p$ such that $q \forces \varphi$.
    \item $p \forces \varphi \land \psi$ if and only if $p \forces \varphi$ and $p \forces \psi$.
    \item $p \forces \exists x \varphi$ if and only if $\{q \in \P : \text{ there is some } \dotz \in U^{\P} \text{ such that }q \forces \varphi(\dotz)\}$ is dense below $p$.
\end{enumerate}
\end{definition}
\begin{lemma}\label{forcingbasic}
Let $\P$ be a forcing notion and $p, q \in \P$.
\begin{enumerate}
    \item If $p \forces \varphi$ and $q \leq p$, then $q \forces \varphi$.
    \item If $\{r \in \P : r \forces \varphi \}$ is dense below $p$, $p \forces \varphi$.
    \item $1_\P \forces \dotx = \dotx$ for every $\dotx \in U^\P$.
    \item If $\dotx \in U^\P$ and $\<\doty, r> \in \dotx$, where $\doty$ is a $\P$-name, then $r \forces \dotx \in \doty$.
    \item If $p \forces \varphi(\dotx, \dotu_1, ..., \dotu_n)$ and $p \forces \dotx = \doty$, then $p \forces \varphi(\doty, \dotu_1, ..., \dotu_n)$.
    \item Let  $\varphi$ be an atomic formula in the language of urelement set theory. For every $x, y$, $\varphi(x, y)$ if and only $p \forces \varphi(\check{x}, \check{y})$. \qed
\end{enumerate}
\end{lemma}
\noindent The following lemma verifies that $U^\P$, unlike $U^P_\#$, is closed under mixtures.
\begin{lemma}\label{mixinglemmaforUP}
Let $\P$ be a forcing poset. Then for every function $f : dom(f) \rightarrow U^\P$, where $dom(f)$ is an antichain in $\P$, there is a $\dotx \in U^\P$ ( \textit{a} \textit{mixture} of $f$) such that $p \forces f(p) = \dotx$ for every $p \in dom(f)$.
\end{lemma}
\begin{proof}
Define $\dotx$ as follows.
\begin{align*}
    \dotx = \bigcup_{p \in dom(f)}\{\<y, r> \in dom(f(p)) \times \P :  \exists q \  ( \<y, q> \in f(p) \land r \leq p, q )\}.
\end{align*}
\noindent Let us first check that $\dotx$ is a $\P$-name satisfying the incompatibility condition. Suppose that $\<a, r_1>, \<y, r_2> \in \dotx$ and $a \neq y$ for some urelement $a$. Then there is some $p \in dom(f)$ and some $q$ such that $\<a, q_1> \in f(p)$ and $r_1 \leq p, q_1$. If $y$ is also in $dom(f(p))$, then there is some $q_2$ such that $\<y, q_2> \in f(p)$ and $r_2 \leq q_2$, which means $r_1$ and $r_2$ are incompatible. If $y \in dom(p')$ for some $p' \in dom(f)$ distinct from $p$, then $\<y, q> \in f(p')$ for some $q'$ and $r_2 \leq p', q'$. Then $r_2$ is incompatible with $r_1$ because $dom(f)$ is an antichain.

Consider any $p \in dom(f)$. It remains to show that $p \forces \dotx = f(p)$. $p \forces \dotx \subseteq f(p)$ because if $\<\doty, r> \in \dotx$ and $q \leq p, r$, then $q \forces \doty \in f(p)$ since $r \forces \doty \in f(p)$. To show $p \forces f(p) \subseteq \dotx$, suppose that $\<\doty, q> \in f(p)$ and $r \leq p, q$. Then $\<\doty, r> \in \dotx$ and hence $r \forces \doty \in \dotx$. 

Finally, we show that $p \forces \dotx \Aeq f(p)$. Fix a condition $s\leq p$. 

\noindent \textit{Case 1.} There is some $\<b, r> \in \dotx$ for some urelement $b$ such that $s$ and $r$ are compatible. Then for some $p' \in dom(f)$ and $\<b, r'> \in f(p')$, $r \leq p', r'$. It follows that $p$ and $p'$ are compatible and hence $p = p'$. Thus, $s$ has an extension in the set $\{q \in \P : \exists a, r_1, r_2 (\A(a) \land \<a, r_1> \in \dotx \land \<a, r_2> \in f(p) \land q \leq r_1, r_2 )\}$. 

\noindent \textit{Case 2.} There is no $\<a, r> \in \dotx$ such that $a$ is an urelement and $s$ is compatible with $r$. Then for every $\<a, r'> \in f(p)$, $s$ cannot be compatible with $r'$ either since otherwise $\<a, s'>$ will be in $\dotx$ for some $s' \leq s, r'$. Thus, $s$ is in the set $\{q \in \P : \forall \<a_1, r_1> \in \dotx \ (\A(a_1) \rightarrow q \bot r_1) \land \forall \<a_2, r_2> \in f(p) \ (\A(a_2) \rightarrow q \bot r_2)\}$. This shows that  $p \forces \dotx \Aeq f(p)$.\end{proof}

\begin{theorem}\label{thm:fulness<->collection}
Over $\ZFCUR$, the following are equivalent.
\begin{enumerate}
    \item Collection.
    \item Fullness Principle: for every forcing poset $\P$, its forcing relation $\forces$ is full.
\end{enumerate}
\end{theorem}
\begin{proof}
The argument for $(1) \rightarrow  (2)$ is standard given that we know $U^\P$ is closed under mixtures for every $\P$. So fix some forcing poset $\P$ and suppose that $p \forces \exists y \varphi$ (with parameters suppressed). By AC, there is a maximal antichain $I$ in the subposet $\{q \in \P : q \leq p \land \exists \dotx \in U^\P q \forces \varphi(\dotx) \}$. By Collection and AC, there is a function $f : I \rightarrow U^\P$ such that for every $q \in I$, $q \forces \varphi(f(q))$. It follows from Lemma \ref{mixinglemmaforUP} that there is some $\doty \in U^\P$ such that $q \forces \doty = f(q)$ and hence $q \forces \varphi(\doty)$ for every $q \in I$. It is then routine to check that $p \forces \varphi(\doty)$.

$(2) \rightarrow (1).$ Assume $(2)$ and suppose that $\forall x \in w \exists y \varphi (x, y, u)$ for some set $w$ and parameter $u$. Define the forcing poset $\mathbb{W}$ to be $w \cup \{w\}$, where for every $p, q \in \mathbb{W}$, $p \leq q$ if and only if $p = q $ or $q = w$.

\begin{claim}
For every $\doty \in U^{\mathbb{W}}$ and $p \in w$, there is some $y$ such that $ker(y) \subseteq ker(\mathbb{W}) \cup ker(\doty)$ and $p \forces \doty = \check{y}$.
\end{claim}
\begin{claimproof}
By induction on the rank of $\doty$. We may assume that there is no urelement $a$ such that $\<a, p> \in \doty$ or $\<a, 1_{\mathbb{W}}> \in \doty$, since otherwise $p \forces \doty = \check{a}$ for such $a$. Define
\begin{align*}
  y = \{z : \exists \dotz \in dom(\doty) \cap U^\mathbb{W} \ (ker(z) \subseteq ker(\dotz) \cup ker(\mathbb{W}) \land p \forces \dotz \in \doty \land p \forces \dotz = \check{z})\}.
\end{align*}
$y$ is a set by Lemma \ref{forcingbasic} (6), and it is clear that $ker(y) \subseteq ker(\mathbb{W}) \cup ker(\doty)$. $p \forces \doty \Aeq \check{y}$ by the assumption. To show that $p \forces \doty \subseteq \check{y}$, observe that if $\<\dotz, r> \in \doty$ and $q \leq p, r$, then $p = q$ and $p \forces \dotz \in \doty$; and since by the induction hypothesis $p \forces \dotz = \check{z}$ for some $z \in y$, it follows that $p \forces \dotz \in \check{y}$. $p \forces \check{y} \subseteq \doty$ because for any $z \in y$,  $p \forces \dotz = \check{z}$ and $p \forces \dotz \in \doty$ for some $\mathbb{W}$-name $\dotz \in dom(\doty)$, and hence $p \forces \check{z} \in \doty$.
\end{claimproof}

\noindent By induction on the complexity of formulas and Lemma \ref{forcingbasic} (6), it follows that for every $p \in \mathbb{W}$ and $x_1, ..., x_n$, $\psi(x_1, ..., x_n)$ if and only if $p \forces \psi(\check{x_1}, ..., \check{x_n})$ for any formula $\psi$ in the language of urelement set theory.

Next we define a $\P$-name $\dotx = \{\<\check{z}, p> : \neg \A (p) \land z \in p \land p \in w\} \cup \{\<a, a> : a \in w \land \A(a) \}$.
\begin{claim}
For every $p \in w$, $p \forces \check{p} = \dotx.$
\end{claim}
\begin{claimproof}
If $p = a $ for some urelement $a$, then it is clear from the definition that $p \forces \check{a} = \dotx$. If $p$ is a set, then $p \forces \check{p} \Aeq \dotx$ because if $\<a, a> \in \dotx$ then $p \bot a$. $p \forces \check{p} \subseteq \dotx$ because whenever $z \in p$, $p \forces \check{z} \in \dotx$. $p \forces \dotx \subseteq \check{p}$ because if $\<\check{z}, r> \in \dotx$ and $q \leq p, r$, then $r= q= p$ and $z \in p$, so $p \forces \check{z} \in \check{p}$.
\end{claimproof}

Now for every $p \in w$, there is some $y$ such that $\varphi(p, y, u)$; so $p \forces \varphi(\check{p}, \check{y}, \check{u})$ and hence $p \forces \varphi(\dotx, \check{y}, \check{u})$. By the definition of $\forces$, this means $1_\mathbb{W} \forces \exists y \varphi(\dotx, y, \check{u})$. Since $\forces$ is full for $U^\mathbb{W}$, it follows that there is a $\doty \in U^\mathbb{W}$ with $1_\mathbb{W} \forces \varphi(\dotx, \doty, \check{u})$; so for every $p \in w$, $p \forces \varphi(\check{p}, \doty, \check{u})$. By the first claim $p \forces \doty = \check{y}$ for some $y$ such that $ker(y) \subseteq ker(\doty) \cup ker(\mathbb{W})$. Let $A = ker(\doty) \cup ker(\mathbb{W})$. For every $p \in w$, there is some $y \in V(A)$ such that $p \forces \varphi( \check{p}, \check{y}, \check{u})$ and hence $\varphi(p, y, u)$. This suffices for Collection by Proposition \ref{weakcollection}.
\end{proof}
\noindent It is folklore that Fullness Principle implies AC. I include a proof of this for completeness.

\begin{theorem}\label{thm:fulness->AC}
$\ZFUR \vdash$ Fullness Principle $\rightarrow$ AC.
\end{theorem}
\begin{proof}
Let $w$ be a set of non-empty sets. Define the forcing poset $\mathbb{W}$ to be $w \cup \{w\}$, where for every $p, q \in \mathbb{W}$, $p \leq q$ if and only if $p = q $ or $q = w$. As before, since $\mathbb{W}$ is trivial, for every $\doty \in U^{\mathbb{W}}$ and $p \in w$, there is some $y$ such that $p \forces \doty = \check{y}$. And consequently, for every $x_1, ... x_n$, $\psi(x_1, ..., x_n)$ if and only if $p \forces \psi(\check{x_1}, ..., \check{x_n})$. Define $\dotx = \{\<\check{z}, p> :  z \in p \land p \in w\}$. The same argument as before shows that $p \forces \check{p} = \dotx $ for every $p \in w$. For every $p \in w$, since it is non-empty, there is some $z \in p$ such that $p \forces \check{z} \in \check{p}$ and hence $p \forces \check{z} \in \dotx$. This shows that $1_\mathbb{W} \forces \exists y (y \in \dotx)$. By Fullness Principle, there is a $\doty \in U^\mathbb{W}$ such that $p \forces \doty \in \dotx$ for every $p \in w$. Thus, $p \forces \doty \in \check{p}$ for every $p \in w$. Now define a function $f$ on $w$ such that $f(p) = z$ if and only if $p \forces \check{z} = \doty$. It follows that $f$ is a choice function on $w$.\end{proof}

\begin{corollary}
Over $\ZFUR$,
\begin{enumerate}
    \item Fullness Principle $\rightarrow$ RP;
    \item RP $\nrightarrow$ Fullness Principle.
\end{enumerate}
\end{corollary}
\begin{proof}
The implication follows from Theorem \ref{thm:fulness->AC}, Theorem \ref{thm:fulness<->collection} and Lemma \ref{easyimplication} (6). The implication cannot be reversed by the Basic Fraenkel Model (Example \ref{exp:BasicFModel}).
\end{proof}

\subsection{Forcing extensions and the forcing theorem}
\begin{definition}\label{m[g]def}
Let $M$ be a countable transitive model of $\ZFUR$, $\P \in M$ be a forcing poset and $G$ be an $M$-generic filter over $\P$.
\begin{enumerate}
    \item $M^\P = U^\P \cap M$
   \item For every $\dot{x} \in M^\P$, 
        \subitem (i) $\dotx_G = a$ if  $\mathcal{A}(a)$ and  $\langle a, p \rangle \in \dot{x}$ for some $p \in G$; 
        \subitem (ii) $\dotx_G = \{ \dot{y}_G: \langle \doty , p \rangle \in \dotx \text{ for some } \doty \in M^\P \text{ and } p \in G \}$ otherwise.
     \item $M[G] = \{\dotx_G : \dot{x} \in M^\P \}$. 
\end{enumerate}
\end{definition}
\noindent Note that $\dotx_G$ is well-defined by clause (ii) in Definition \ref{newpnames}. It is shown in \ref{subsection:m[g]=m[g]sharp} that $M[G]$ is in fact the same as $M[G]_\#$.
\begin{lemma}
Let $M$ be a countable transitive model of $M$, $\P \in M$ be a forcing poset, and $G$ be an $M$-generic filter over $\P$. Then
\begin{enumerate}
    \item  $M \subseteq M[G]$;
    \item  $G \in M[G]$;
    \item  $M[G]$ is transitive;
    \item  $Ord \cap M = Ord \cap M[G]$;
    \item  For every transitive model $N$ of $\ZFUR$ such that $G\in N$ and $M \subseteq N$, $M[G] \subseteq N$;
    \item $\mathcal{A} \cap M = \mathcal{A} \cap M[G]$.
\end{enumerate}
\end{lemma}
\begin{proof}
(1)--(5) are all proved by standard text-book arguments as in \cite[Ch.VII]{kunen2014set}. (6) is clear by the construction of $M[G]$ because every urelement in $M[G]$ must come from $ker(\dotx)$ for some $\dotx \in M^\P$. 
\end{proof}

\begin{lemma}\label{McoversM[G]}
$ker(\dot{x}_G) \subseteq ker(\dot{x})$, for every $\dot{x} \in M^{\P}$. Hence, every set of urelements in $M[G]$ is a subset of some set of urelements in $M$.
\end{lemma}
\begin{proof}
By induction on the rank of $\dot{x}$, and we may assume that $\dotx_G$ is a set. Since $ker(\dot{x}_G) \subseteq \bigcup \{ ker(\dot{y}_G) : \dot{y} \in dom(\dot{x})\}$ and by the induction hypothesis $ker(\dot{y}_G) \subseteq ker(\dot{y}) \subseteq ker(\dot{x})$ for every $\dot{y} \in dom(\dot{x})$, the result follows.
\end{proof}

\begin{theorem}[The Forcing Theorem for $\forces$]
Let $M$ be a countable transitive model of $\ZFUR$, $\P \in M$ be a forcing poset. Then for every $\dotx_1, ..., \dotx_n \in M^\P$,

\begin{enumerate}
    \item For every $M$-generic filter $G$ over $\P$, $M[G] \models \varphi (\dotx_{1_G}, ..., \dotx_{n_G})$ if and only if $\exists p \in G (p \forces \varphi(\dotx_1, ..., \dotx_n))^M$.
    \item For every $p \in \P$, $(p \forces \varphi(\dotx_1, ..., \dotx_n))^M$ if and only if for every $M$-generic filter $G$ over $\P$ such that $p \in G$, $M[G] \models \varphi (\dotx_{1_G}, ..., \dotx_{n_G})$. 
\end{enumerate}
\end{theorem}
\begin{proof}
(2) is an easy consequence of (1) and the proof of the Boolean cases and quantifier case of (1) is the same as in \cite[Chapter VII. Theorem 3.5]{kunen2014set}. So it remains to show that the atomic cases for (1) hold.

\noindent \textit{Case 1}. $\varphi(\dotx_1, \dotx_2)$ is $\dotx_1 \in \dotx_2$. The argument is the same as in \cite[Chapter VII, Theorem 3.5]{kunen2014set}.

\noindent \textit{Case 2}. $\varphi(\dotx)$ is $\A(\dotx)$. Suppose that $\dotx_G$ is some urelement $b$. Then $\<b, p> \in \dotx$ for some $p \in G$, so $\{q \in \P : \exists \<a , r> \in \dotx \ (\A(a) \land q \leq r)\}$ is dense below $p$ and hence $p \forces \A(\dotx)$. Suppose that $p \forces \A(\dotx)$ for some $p \in G$. Then there is some $q \in G$ such that $\<b, r> \in \dotx$ for some $r \geq q$ and urelement $b$. Thus, $\dotx_G = b$.

\noindent \textit{Case 3}. $\varphi(\dotx_1, \dotx_2)$ is $\dotx_1 = \dotx_2$. For the left-to-right direction of (1), suppose that $\dotx_{1_G} = \dotx_{2_G}$.

\textit{Subcase 3.1}. $\dotx_{1_G} = \dotx_{2_G} = b$ for some urelement $b$. Then $\<b, s_1> \in \dotx_1$ and $\<b, s_2> \in \dotx_2$ for some $s_1, s_2 \in G$. Fix some $p \in G$ such that $p \leq s_1, s_2$. Observe first that $p \forces \dotx_1 \subseteq \dotx_2$ and $p \forces \dotx_2 \subseteq \dotx_1$ trivially hold: for any $\P$-name $\doty$ and $r \in \P$ such that $\<\doty, r> \in \dotx_1 (\text{or } \dotx_2)$, $p$ must be incompatible with $r$ because $r$ is incompatible with $s_1(\text{or }s_2)$. Moreover, $p \forces \dotx_1 \Aeq \dotx_2$ because $\{q \in \P : \exists a, r_1, r_2 (\A(a) \land \<a, r_1> \in \dotx_1 \land \<a, r_2> \in \dotx_2 \land  q \leq r_1, r_2 )\}$ is clearly dense below $p$. Hence, $p \forces \dotx_1 = \dotx_2$.

\textit{Subcase 3.2}. $\dotx_{1_G}$ is a set.  We first use a standard text-book argument to show that there is some $p \in G$ such that $p \forces \dotx_1 \subseteq \dotx_2$ and $p \forces \dotx_2 \subseteq \dotx_1$. Define:
\begin{itemize}
    \item [] $D_1= \{p \in \P : p \forces \dotx_1 \subseteq \dotx_2 \land p \forces \dotx_2 \subseteq \dotx_1 \}$
    \item [] $D_2 = \{p \in \P : \exists \<\doty_1, q_1> \in \dotx_1 \ (p \leq q_1 \land \forall \<\doty_2, q_2> \in \dotx_2 \ \forall r \leq q_2 \ (r \forces \doty_1 = \doty_2 \rightarrow p \bot r )) \}$
    \item [] $D_3 = \{p \in \P : \exists \<\doty_2, q_2> \in \dotx_2 \ (p \leq q_2 \land \forall \<\doty_1, q_1> \in \dotx_1 \ \forall r \leq q_1 \  (r \forces \doty_2 = \doty_1 \rightarrow p \bot r )) \}$
\end{itemize}
 If $p \nVdash \dotx_1 \subseteq \dotx_2$, then there are $\<\doty_1, q_1> \in \dotx_1$ and $r \leq p, q_1$ such that $r \nforces \doty_1 \in \dotx_2$; so there is an $s \leq r$ such that for every $\<\doty_2, q_2> \in \dotx_2$ and $s' \leq q_2$, if $s' \forces \doty_1 = \doty_2$, then $s \bot s'$. Hence, $s \leq p$ and $s \in D_2$. Similarly, if $p \nforces \dotx_2 \subseteq \dotx_1$, then $p$ will have an extension in $D_3$. This shows that $D_1 \cup D_2 \cup D_3$ is dense. However, $G \cup (D_2 \cup D_3)$ must be empty. Suppose \textit{for reductio} that $p \in G\cap D_2$. Fix some $\<\doty_1, q_1> \in \dotx_1$ with $p \leq q_1$ that witnesses $p \in D_2$. It follows that $\doty_1{_G} = \doty_2{_G}$ for some $\<\doty_2, q_2> \in \dotx_2$ with $q_2 \in G$. By the induction hypothesis, there is some $r \in G$ such that $r \leq q_2$ and $r \forces \doty_1 = \doty_2$. But $p$ must be incompatible with such $r$, which is a contradiction. The same argument shows that $G \cap D_3$ is empty. Therefore, there is some $p \in G$ such that $p \forces \dotx_1 \subseteq \dotx_2$ and $p \forces \dotx_2 \subseteq \dotx_1$.

Now I wish to find some $q \in G$ such that $q \forces \dotx_1 \Aeq \dotx_2$. Define:
\begin{itemize}
    \item [] $E_1 = \{ q \in \P : \forall r \leq q \ [\forall \<a_1, s_1> \in \dotx_1 \ (\A(a) \rightarrow r \bot s_1) \land \forall \<a_2, s_2> \in \dotx_2 \ (\A (a_2) \rightarrow r \bot s_2)]\}$.
    \item [] $E_2 = \{ q \in \P : \exists \<a, r> \in \dotx_1 \ (\A(a) \land q \leq r) ) \}$.
    \item [] $E_3 = \{ q \in \P : \exists \<a, r> \in \dotx_2 \ (\A(a) \land q \leq r)\}$.
\end{itemize}
$E_1 \cup E_2 \cup E_3$ is dense. But if there is some $q \in G \cap (E_2 \cup E_3)$, either $\dotx_{1_G}$ or $\dotx_{2_G}$ would be an urelement. Thus there is some $q \in G \cap E_1$ such that the set
\begin{itemize}
    \item [] $\{r \in \P : \forall \<a_1, s_1> \in \dotx_1 \ (\A(a_1) \rightarrow r \bot s_1) \land \forall \<a_2, s_2> \in \dotx_2 \ (\A (a_2) \rightarrow r \bot s_2)\}$
\end{itemize}
is dense below $q$. Therefore, $q \forces \dotx_1 \Aeq \dotx_2$. A common extension of $p$ and $q$ in $G$ will then force $\dotx_1 = \dotx_2$.

For the right-to-left direction of (1) in Case 3, suppose that for some $p \in G$, $p \forces \dotx_1 = \dotx_2$.

\textit{Subcase 3.3}. $\dotx_{1_G} = b$ for some urelement $b$. Then $\<b, r> \in \dotx_1$ for some $r \in G$. Define:
\begin{itemize}
    \item [] $F_1 = \{ q \in \P : \exists a, s_1, s_2 (\A (a) \land \<a, s_1> \in \dotx_1 \land \<a, s_2> \in \dotx_2 \land q \leq s_1, s_2) \}$.
    \item [] $F_2 = \{q \in \P : \forall \<a, s_1> \in \dotx_1 \ (\A(a) \rightarrow q \bot s_1) \land \forall \<a, s_2> \in \dotx_2 \ (\A(a) \rightarrow q \bot s_2 )\}$.
\end{itemize}
Since $p \forces \dotx_1 \Aeq \dotx_2$, $F_1 \cup F_2$ is dense below $p$. But clearly $F_2 \cap G$ is empty as $\<b, r> \in \dotx_1$, so there is some $q \in F_1 \cap G$. It follows that $\<b, s_1> \in \dotx_1$ and $\<b, s_2> \in \dotx_2$ for some $s_1, s_2 \in G$. Therefore, $\dotx_{2_G} = b = \dotx_{2_G}$.

\textit{Subcase 3.4}. $\dotx_{1_G}$ is a set. Suppose \textit{for reductio} that $\dotx_{2_G}$ is some urelement $b$ and so $\<b, r> \in \dotx_2$ for some $r \in G$. Since $p \forces \dotx_1 \Aeq \dotx_2$, it follows that there are some urelement $a$ and  $s \in G$ such that $\<a, s> \in \dotx_1$. This implies that $\dotx_{1_G} = a$, which is a contradiction. Hence, $\dotx_{2_G}$ is a set, so it remains to show that $\dotx_{1_G}$ and $\dotx_{2_G}$ have the same members. If $\doty_G \in \dotx_{1_G}$, then $\<\doty, r> \in \dotx_1$ for some $r \in G$. So there is some $q \in G$ with $q \leq p, r$, and since $p \forces \dotx_1 \subseteq \dotx_2$, $q \forces \doty \in \dotx_2$. By the induction hypothesis, $\doty_G \in \dotx_{2_G}$. The same argument will show that $\dotx_{2_G} \subseteq \dotx_{1_G}$. \end{proof}

\subsection{The fundamental theorem of forcing with urelements}
\begin{theorem}\label{forcingpreservesZFCU}
Let $M$ be a countable transitive model of $\ZFUR$, $\P \in M$ be a forcing poset, and $G$ be an $M$-generic filter over $\P$. Then
\begin{enumerate}
    \item $M[G]$ is a countable transitive model of ZU;
    \item $M[G] \models$ AC if $M \models$ AC;
    \item $M[G] \models \ACA$ if  $M \models \ACA$;
    \item $M[G] \models$ Collection if $M \models$ Collection.
\end{enumerate}
\end{theorem}
\begin{proof}
The proof of (1) and (2) are the same as in Kunen \cite[Ch.VII]{kunen2014set} and hence omitted. (3) follows from Lemma \ref{McoversM[G]} that every set of urelements in $M[G]$ is covered by some set of urelements in $M$. 

For (4), suppose that $M[G] \models \forall v \in \dotw_G\ \exists y \varphi(v, y, \dot{u}_G)$ for some $\dotw_G$ and $\dot{u}_G$. In $M$, define 
\begin{align*}
    x = \{\langle \dotx, p \rangle \in (dom(\dotw) \cap M^\P) \times \P :  \exists \dot{y} \in M^\P p \forces \varphi(\dotx, \dot{y}, \dot{u}) \}.
\end{align*}
By Collection in $M$, there is a set of $\P$-names $v$ such that for every $\langle \dotx, p \rangle \in x$, there is a $\dot{y} \in v$ with $p \forces \varphi (\dotx, \dot{y}, \dot{u})$. Define $\dot{v}$ to be $v \times \{1_\P\}$. It is now routine to check that $M[G] \models \forall x \in \dot{w}_G \ \exists y \in \dot{v}_G \ \varphi(x, y, \dot{u}_G)$.
\end{proof}
A more difficult question is whether forcing preserves Replacement when the ground model $M$ does not satisfy Collection. When $M$ is a model of ZF, the standard argument for $M[G] \models $ Replacement appeals to Collection in $M$. But this move is not allowed when $M$ only satisfies $\ZFCUR$. A new argument is thus needed.

\begin{definition}\label{purification}
Let $\P$ be a forcing poset and $A$ be a set of urelements. For every urelement $a$, let $\overset{A}{a} = a$. For every $\dotx \in U^\P$, we define the $A$\textit{-purification of }$\dotx$, $\overset{A}{\dot{x}}$, as follows.
\begin{align*}
   \overset{A}{\dot{x}} = \{\langle \overset{A}{y}, p \rangle : \langle y , p \rangle \in \dot{x} \land ( y \in U^\P \lor y \in A) \}. 
\end{align*}
\end{definition}
\noindent That is, we get $\overset{A}{\dot{x}}$ by hereditarily throwing out the urelements used to build $\dot{x}$ that are not in $A$.
\begin{prop}
Let $\P$ be a forcing poset and $A$ be a set of urelements such that $ker(\P) \subseteq A$. For every $\dotx \in U^\P$, $\overset{A}{\dot{x}} \in U^\P$ and $ker(\overset{A}{\dot{x}}) \subseteq A$.
\end{prop}
\begin{proof}
By induction on the rank of $\dotx$. To show that $\overset{A}{\dot{x}}$ is always a $\P$-name, we only need to check the incompatibility condition in Defnition \ref{newpnames} holds. Suppose that $\<a, p>, \<y, q> \in \overset{A}{\dot{x}}$, where $a$ is an urelement and $y \neq a$. If $y$ is another urelement in $dom(\dotx)$, then $p$ and $q$ are incompatible; otherwise $y$ is some $\overset{A}{\dotz}$, where $\<\dotz, q> \in \dotx$ and $\dotz$ is a $\P$-name, then $p$ and $q$ are incompatible because no urelement is a $\P$-name. $ker(\overset{A}{\dot{x}}) \subseteq A$ because $ker(\overset{A}{\dot{x}})$ is contained in $\bigcup_{y \in dom(\dotx)}ker(\overset{A}{\dot{y}}) \cup ker(\P)$, which is a subset of $A$ by the induction hypothesis.
\end{proof}

\begin{lemma}\label{iso->id}
Let $\bar{x}$ and $\bar{y}$ be transitive sets. Every surjective $\in$-isomorphism from $\barx$ to $\bary$ that fixes $\emptyset$ and every urelement is the identity map.
\end{lemma}
\begin{proof}
Let $f$ be such surjective $\in$-ismorphism. We show that $x = f(x)$ for every set $x \in \barx$ by $\in$-induction. Fix some $x \in \bar{x}$. By the induction hypothesis, $x \subseteq f(x)$. If $y\in f(x)$, then $y = f(w)$ for some $w \in x$ so by the induction hypothesis $y \in x$. Therefore, $x = f(x)$.
\end{proof}

\begin{theorem}\label{forcingpreservesreplacement}
Let $M$ be a countable transitive model of $\ZFUR$, $\P \in M$ be a forcing poset and $G$ be $M$-generic over $\P$. Then $M[G] \models$ Replacement.
\end{theorem}
\begin{proof}
Suppose that for some $\dot{w}_G$ and $\dot{u}_G$ in $M[G]$, $M[G] \models \forall x \in \dot{w}_G\ \exists ! y \varphi(x, y, \dot{u}_G)$. Let $A = ker(\dot{w}) \cup ker(\P) \cup ker(\dot{u})$. By Theorem \ref{forcingpreservesZFCU}, we may assume $M$ does not satisfy Collection and hence has a proper class of urelements. 

\begin{lemma}\label{keylemmarep}
For every $\dot{v}_G \in \dot{w}_G$, there exist $p \in G$ and $\mu' \in M^{\P}$ such that $p \forces \varphi(\dot{v}, \mu', \dot{u})$ and $ker(\mu') \subseteq A$.
\end{lemma}
\begin{proof}
Fix a $\dot{v}_G \in \dot{w}_G$ for some $\dot{v} \in dom (\dot{w}) \cap M^\P$. Since $M[G] \models \exists ! y (\dot{v}_G, y, \dot{u}_G)$, there is a $\P$-name $\mu$ and a $p \in G$ such that $ p \forces \varphi (\dot{v}, \mu, \dot{u}) \land \forall z (\varphi (\dot{v}, z, \dot{u}) \rightarrow \mu = z)$.
\begin{claim}\label{claim1}
 For every $M$-generic filter $H$ over $\P$ such that $p \in H$, $ker(\mu_H) \subseteq A$.
\end{claim}
\begin{claimproof}
Suppose not. Then there is some $b \in ker(\mu_H) \setminus A$. Since $M$ has a proper class of urelements, there is some urelement $c \in M$ such that $c \notin A \cup ker(\mu)$. In $M$, let $\pi$ be the automorphism that only swaps $b$ and $c$. Since $\pi$ point-wise fixes $A$, it follows that 
\begin{align*}
    p \forces \varphi (\dot{v}, \pi \mu, \dot{u}) \land \forall z (\varphi (\dot{v}, z, \dot{u}) \rightarrow \pi \mu = z).
\end{align*}
Thus, $M[H] \models \mu_H = (\pi\mu)_H$. Since $b \in ker(\mu_H)$, $\pi b  \in ker(\pi\mu_H)$($\pi$ is viewed as an automorphism of the background universe); but $\pi b = c \notin ker(\mu)$ and $ker(\mu_H) \subseteq ker(\mu)$, so $\pi b \notin ker(\mu_H)$, which is a contradiction.
\end{claimproof}

\noindent Note that we cannot hope to show that $ker(\mu) \subseteq A$ in general. For if $\mu^*$ is some $\P$-name such that  $\mu^* = \mu \cup \{\langle \{\<b, 1_\P>\}, q \rangle\}$, where $b$ is an urelement not in $A$ and $q$ is not compatible with $p$, we would still have $p \forces \mu = \mu^*$.

\begin{claim}\label{claim2}
Let $H$ be an $M$-generic filter over $\P$ such that $p \in H$. For every $\dot{x}, \dot{y} \in M^{\P}$, if $\dot{x}_H, \dot{y}_H \in trc(\{\mu_H\})$, then $\dot{x}_H = \dot{y}_H$ if and only if $(\overset{A}{\dot{x})}_H = (\overset{A}{\dot{y}})_H$.
\end{claim}
\begin{claimproof}
If $\dotx_H = \doty_H = a$ for some urelement $a$, then by Claim \ref{claim1} $a \in A$. Then it is easy to check that $(\overset{A}{\dot{y})}_H = (\overset{A}{\dot{x}})_H = a$. If $(\overset{A}{\dot{y})}_H = (\overset{A}{\dot{x}})_H = b$ for some urelement $b$, then $b \in A$ and it follows that $\dotx_H = \doty_H = b$.

So suppose $\dotx_H = \doty_H$ are sets in $trc(\{\mu_H\})$ and the claim holds for every $\dot{z} \in dom(\dot{x}) \cup dom(\dot{y})$. Clearly, $(\overset{A}{\dot{x})}_H$ and $(\overset{A}{\dot{y}})_H$ must also be sets. If $\overset{A}{\dot{z}}_H \in \overset{A}{\dot{x}}_H$ for some $\dotz \in M^\P \cap dom(\dotx)$, we have $\dot{z}_H \in \dot{y}_H = \dot{x}_H$. So there is some $\dot{w} \in M^\P \cap dom(\doty)$ such that $\dot{w}_H = \dot{z}_H$. $\dot{z}_H \in trc(\{\mu_H\})$ so by the induction hypothesis $\overset{A}{\dot{z}}_H = \overset{A}{\dot{w}}_H \in (\overset{A}{\dot{y}})_H$. This shows that $\overset{A}{\dot{x}}_H \subseteq \overset{A}{\dot{y}}_H$, and we will have $\overset{A}{\dot{x}}_H = \overset{A}{\dot{y}}_H$ by the same argument.

Now suppose that $\dotx_H, \doty_H \in trc(\{\mu_H\})$ and $\overset{A}{\dot{x}}_H = \overset{A}{\dot{y}}_H$ are sets. Then $\dotx_H$ and $\doty_H$ must be sets. For if, say, $\dotx_H = a$ for some urelement $a$, then $a \in A$ by Claim \ref{claim1}, which implies that $\overset{A}{\dot{x}}_H = a$. Let $\dot{z}_H \in \dot{x}_H$ for some $\dotz \in M^\P \cap dom(\dotx)$. Then $\overset{A}{\dot{z}}_H \in \overset{A}{\dot{y}}_H$ and so $\overset{A}{\dot{z}}_H = \overset{A}{\dot{w}}_H$ for some $\dot{w}_H \in \dot{y}_H$. By the induction hypothesis, it follows that $\dot{z}_H = \dot{w}_H$. This shows that $\dot{x}_H \subseteq \dot{y}_H$ and consequently, $\dot{x}_H = \dot{y}_H$.
\end{claimproof}

\begin{claim}\label{claim3}
$p \forces \overset{A}{\mu} = \mu$.
\end{claim}
\begin{claimproof}
 Let $H$ be an $M$-generic filter on $\P$ that contains $p$. We show that $\overset{A}{\mu}_H = \mu_H$. Let $f$ be the function on $trc(\{\mu_H\})$ that sends every $\dot{y}_H$ to $\overset{A}{\dot{y}}_H$, which is is well-defined by Claim \ref{claim2}. By Lemma \ref{iso->id}, it suffices to show that $f$ maps $trc(\{\mu_H\})$ onto $trc(\{\overset{A}{\mu}_H\})$, preserves $\in$ and fixes all the urelements.

\textit{ $f$ preserves $\in$}. Consider any $\dot{y}{_H}, \dot{x}{_H} \in trc(\{\mu_H\})$. Suppose that $\dot{y}{_H} \in \dot{x}{_H}$. Then $\dot{y}{_H} = \dot{z}_H$ for some $\dot{z} \in M^\P \cap dom(\dot{x})$ so $ \overset{A}{\dot{z}}_H \in \overset{A}{\dot{x}}_H$; by Claim \ref{claim2}, it follows that $\overset{A}{\dot{y}}_H = \overset{A}{\dot{z}}_H \in \overset{A}{\dot{x}}_H$. Suppose that $\overset{A}{\dot{y}}_H \in \overset{A}{\dot{x}}_H$. Then $\overset{A}{\dot{y}}_H = \overset{A}{\dot{z}}_H$ for some $\dot{z}_H \in \dot{x}_{H}$ so $\dot{y}{_H} = \dot{z}_H \in \dot{x}_{H}$ by Claim \ref{claim2} again.

\textit{$f$ maps $trc(\{\mu_H\})$ onto $trc(\{\overset{A}{\mu}_H\})$}. If $\dot{y}_H \in trc(\{\mu_H\})$, then $\dot{y}_H \in \dot{y}_1{_H} \in ... \in \dot{y}_n{_H} \in \mu_H$ for some $n$. Since $f$ is $\in$-preserving, it follows that $\overset{A}{\dot{y}}_H \in \overset{A}{\dot{y}_1}_H \in ... \in \overset{A}{\dot{y}_n}_H \in \overset{A}{\mu}_H$ and hence $\overset{A}{\dot{y}}_H \in trc(\{{\overset{A}{\mu} }_H\})$. To see it is onto, let $x \in x_1 \in ... \in x_n \in \overset{A}{\mu}_H$. Then $x = \overset{A}{\dot{y}}_H \in \overset{A}{\dot{y}_1}_H  \in ... \in \overset{A}{\dot{y}_n}_H \in \overset{A}{\mu}_H$, but then $\dot{y}_H \in \dot{y}_1{_H} \in ... \in \dot{y}_n{_H} \in \mu_H$ and hence  $\dot{y}_H \in trc(\{\mu_H\})$.

\textit{$f$ fixes all the urelements in  $trc(\{\mu_H\})$.} Suppose $\dotx_H = a \in trc(\{\mu_H\})$ for some urelement $a$. Then by Claim \ref{claim1}, $a \in A$ and hence $\overset{A}{\dot{x}}_H = a$. \end{claimproof}

\noindent The lemma is now proved by letting $\mu'$ be $\overset{A}{\mu}$.
\end{proof}
Now in M, we define
\begin{align*}
    \bar{w} = \{\langle \dot{v}, p\rangle \in (dom(\dot{w}) \cap M^\P) \times \P : \exists \mu \in M^{\P} (ker(u) \subseteq A \land p \forces \varphi(\dot{v}, \mu, \dot{u})) \}.
\end{align*}
For every $\langle \dot{v}, p \rangle \in \bar{w}$, let $\alpha_{ \dot{v}, p }$ be the least $\alpha$ such that there is some $\mu \in V_\alpha(A) \cap M^{\P}$ such that $p \forces \varphi(\dot{v}, \mu, \dot{u})$. Let $\beta = Sup_{\langle \dot{v}, p \rangle \in \bar{w}} \alpha_{\dot{v}, p}$ and set $\rho = (V_\beta(A) \cap M^{\P}) \times \{1_\P\}$. It remains to show that $M[G] \models \forall x \in \dot{w}_G\ \exists y \in \rho_G\ \varphi(x, y, \dot{u}_G)$. Let $\dot{v}_G \in \dot{w}_G$. By Lemma \ref{keylemmarep}, there is some $p \in G$ such that $\langle \dot{v}, p \rangle \in \bar{w}$. So there is some  $\P$-name $\mu \in dom(\rho)$ such that $p \forces \varphi(\dot{v}, \mu, \dot{u})$. Thus, $M[G] \models \varphi(\dot{v}_G, \mu_G, \dot{u}_G)$ and $\mu_G \in \rho_G$.\end{proof}

\begin{theorem}[The Fundamental Theorem of Forcing with Urelements]\label{fundamentalthmofforcing}
Let $M$ be a countable transitive model of $\ZFUR$, $\P \in M$ be a forcing poset and $G$ be an $M$-generic fitler over $\P$. Then
\begin{enumerate}
    \item $M[G] \models$ $\ZFUR$.
    \item $M[G] \models \ZFCUR$ if $M\models$ $\ZFCUR$.
    \item $M[G] \models$ ZFCU if  $M \models $ ZFCU. 
    \item $M[G] \models$ Plenitude if $M \models $ Plenitude.
    \item $M[G] \models $ Duplication if $M \models $ Duplication.
    \item $M[G] \models$ Plenitude$^+$ if $M \models $ Plenitude$^+$.
    \item $M[G] \models$ Tail if $M \models$ Tail.
    \item $M[G] \models $ DC$_{<Ord}$ if $M \models $ DC$_{<Ord}$.
    \item $M[G] \models$ RP$^-$ if $M \models$ RP$^-$.
    \item $M[G] \models$ RP if $M \models$ RP.
    \item $M[G] \models $ Closure if $M \models $ Closure + $\ACA$.
\end{enumerate}
\end{theorem}
\begin{proof} 
(1), (2), and (3) are Theorem \ref{forcingpreservesZFCU} and Theorem \ref{forcingpreservesreplacement}.

(4) is clear because $M \subseteq M[G]$ and $M[G]$ and $M$ has the same ordinals. 

(5) follows easily from Lemma \ref{McoversM[G]}.

For (6), suppose that $M \models $ Plenitude$^+$ and $\dotx_G \in M[G]$. Then in $M$, there is a bijection $f$ from $dom(\dotx)$ to a set of urelements. Using $f$ we can code an injective function in $M[G]$ from $\dotx_G$ to $\A$. So Plenitude$^+$ holds in $M[G]$.

For (7), suppose that $M\models$ Tail and $A \subseteq \A$ is in $M[G]$. Then let $A' \in M$ be a set of urelements containing $A$ and $B'$ be a tail of $A'$. It is not hard to check that $A'\setminus A \cup B'$ is a tail of $A$ in $M[G]$.

(8) Suppose that $M \models $ DC$_{<Ord}$. It is a standard result that $\forall \kappa \text{DC}_\kappa$ implies AC, so $M \models$ AC and hence $M[G] \models$ AC. Since  DC$_{<Ord}$ implies that either $\A$ is a set, or Plenitude holds, it follows that $M[G] \models (\A \text{ is a set} \lor \text{Plenitude})$ by (4). By Theorem \ref{maintheorem1}, we have $M[G] \models$ DC$_{<Ord}$. 

(9) Suppose that $M \models$ RP$^-$ and $M[G] \models \varphi(\dotx{_1}_G, ..., \dotx{_n}_G)$. Let $p\in G$ be such that $p \forces \varphi(\dotx_1, ..., \dotx_n)$. By RP$^-$ in $M$, there is a transitive set $m$ containing $\P$ and $\dotx_1, ..., \dotx_n$ such that $(p \forces \varphi(\dotx_1, ..., \dotx_n))^m$ and $m$ satisfies some finite fragment of $\ZFUR$ that suffices for the construction of $\P$-names inside $m$. It then follows that $m[G] \models \varphi(\dotx{_1}_G, ..., \dotx{_n}_G)$. $m[G]$ is a transitive set in $M[G]$ because $\dot{m}= \{\<\doty, 1_\P> : \doty \in m \cap M^\P\}$ is a $\P$-name for $m[G]$. Therefore, $M[G] \models$ RP$^-$.

(10) Suppose that $M \models$RP. Given a formula $\varphi (v_1, ..., v_n)$ and some $\dotu_G \in M[G]$, let $\psi (p, \P, v_1, ..., v_n)$ be the formula asserting that $p$ is a forcing condition in $\P$ and $p \forces \varphi(v_1, ..., v_n)$ for $\P$-names $v_1, ..., v_n$. By RP in $M$, there will a transitive set $m$ containing $\{\P, \dotu \}$ that reflects $\psi$ and satisfies some finite fragment of $\ZFUR$ sufficient for forcing. Then as in the last paragraph, $m[G]$ is a transitive set containing $\dotu_G$ in $M[G]$. If $M[G] \models \varphi (\dotx{_1}_G, ..., \dotx{_n}_G)$ for some $\dotx{_1}_G, ..., \dotx{_n}_G$ in $m[G]$, then there will be $p \in G$ such that $(p \forces \varphi(\dotx{_1}, ..., \dotx{_n}))^m$, and so $m[G] \models \varphi (\dotx{_1}_G, ..., \dotx{_n}_G)$. And if $M[G] \models \varphi (\dotx{_1}_G, ..., \dotx{_n}_G)^{m[G]}$, then there is some $p \in G$ such that $(p\forces \varphi(\dotx{_1}, ..., \dotx{_n}))^m$, so $(p\forces \varphi(\dotx{_1}, ..., \dotx{_n}))^M $ and hence $M[G] \models \varphi (\dotx{_1}_G, ..., \dotx{_n}_G)$. This shows that $M[G] \models $ RP.

(11)  Suppose that in $M$, Closure holds and every set of urelements is well-orderable. Let $x \in M[G]$ be a set of realized cardinals whose supermum is some limit cardinal $\lambda$. It suffices to show that in $M$, every cardinal $\kappa < \lambda$ is realized. Since $\lambda$ is a limit cardinal in $M$, for every $\kappa < \lambda$, there is some cardinal $\kappa'$ in $M[G]$ with $\kappa < \kappa' < \lambda$ that is realized by some $B\subseteq \A$ in $M[G]$; so by Lemma \ref{McoversM[G]} $B \subseteq B'$ for some set of urelements $B' \in M$. Then $(\kappa' \leq |B'|)^M$ since otherwise it would contradict the fact that $\kappa'$ is a cardinal in $M[G]$. Thus, $\kappa$ is realized in $M$.
\end{proof}
\noindent It is unclear if forcing preserves Closure if the ground model does not have $\ACA$.

\subsection{Destroying the $\DCK$ and recovering Collection}\label{dck&collection}
I now move on to the preservation of the $\DCK$. A forcing poset $\P$ is $\kappa$-closed if in $\P$ every infinite descending chain of length less than $\kappa$ has a lowerbound. It is a text-book result that $\kappa$-closed forcing posets preserve cardinalities $\leq \kappa$.
\begin{theorem}\label{kclosedforcingpreservedck}
Let $M$ be a countable transitive model of $\ZFCUR$ + DC$_\kappa$-scheme, $\P \in M$ be such that $(\P \text{ is } \kappa^+\text{-closed})^M$ and $G$ be an $M$-generic fitler over $\P$. Then $M[G] \models$ $\ZFCUR$ + DC$_\kappa$-scheme.
\end{theorem}
\begin{proof}

We first make some definitions. For every $\alpha$-sequence $s$ of $\P$-names, let $\dot{s}^{(\alpha)}$ denote the canonical $\P$-name such that $\dot{s}^{(\alpha)}_G$ is an $\alpha$-sequence in $M[G]$ with $\dot{s}^{(\alpha)}_G(\eta) = s(\eta)_G$ for all $\eta < \alpha$. Given a $p \in \P$ and a suitable formula $\varphi$, a $\kappa$-sequence of the form $\<\<p_\alpha, \dot{x}_\alpha> : \alpha < \kappa>$, where $\<p_\alpha, \dot{x}_\alpha> \in \P \times M^\P$, is said to be a \textit{ $\varphi$-chain below} $p$ if $\<p_\alpha : \alpha < \kappa>$ is a descending chain below $p$ and for every $\alpha < \kappa$, $p_\alpha \forces \varphi(\dot{s}^{(\alpha)}, \dot{x}_{\alpha+1})$ where $s = \<\dot{x}_\eta : \eta < \alpha>$.

Suppose that $M[G] \models \forall x \exists y \varphi(x. y, u)$. There is some $p \in G$ such that $p \forces \forall x \exists y \varphi (x, y, \dot{u})$. Let $D$ be the set of forcing conditions that are a lower bound of some $\varphi$-chain below $p$. We claim that $D$ is dense below $p$. If $r \leq p$, let $\psi(x, y, r, \P)$ be the formula defined as follows.
\begin{itemize}
    \item [] $\psi(x, y, \P, \dot{u}) =_{df}$ if $x = \<\<p_\eta, \dot{x}_\eta> : \eta < \alpha>$, where  $\<p_\eta : \eta < \alpha>$ is a descending chain of length $\alpha$ for some $\alpha < \kappa$, then $y = \<q, \dot{x}> \in \P \times M^\P$ such that $q$ bounds $\<p_\eta : \eta < \alpha>$ and $q \forces \varphi(\dot{s}^{(\alpha)}, \dot{x}, \dot{u})$.
\end{itemize}
Let $\P\downarrow r$ denote the set of conditions in $\P$ below $r$. In $M$, for every $x \in (\P \downarrow r \times M^\P)^{<\kappa}$, since $\P$ is $\kappa$-closed, there is some $y \in \P \downarrow r \times M^\P$ such that $\psi(x, y, \P. \dot{u})$. By DC$_\kappa$-scheme in $M$, there exists a $\varphi$-chain $\<\<p_\alpha, \dot{x}_\alpha> : \alpha < \kappa>$, where $\<p_\alpha : \alpha < \kappa>$ is below $r$ and hence below $p$. $\P$ is $\kappa^+$-closed, so there is some $q$ that bounds this $\varphi$-chain below $p$. Thus, $D$ is dense below $p$. It then follows that there is $q \in G$ that bounds a $\varphi$-chain, $\<\<p_\alpha, \dot{x}_\alpha> : \alpha < \kappa>$, below $p$. Let $s = \<\dot{x}_\alpha : \alpha < \kappa>$ and $f = \dot{s}^{(\kappa)}_G$. $f$ is then a $\kappa$-sequence in $M[G]$ and as $\P$ is $\kappa$-closed, $\kappa$ is the same cardinal in $M[G]$ as in $M$. Moreover, $M[G] \models \varphi (f\restriction \alpha, f(\alpha), u)$ for all $\alpha < \kappa$ because $q \forces \varphi(\dot{s}^{(\alpha)}, \dot{x}_\alpha, \dot{u})$. 
\end{proof}
For any infinite cardinals $\kappa$ and $\lambda$ with $\kappa < \lambda$, $\textup{Col}(\kappa, \lambda)$ is the forcing poset consisting of all partial functions from $\kappa$ to $\lambda$ of size less than $\kappa$ (ordered by reverse inclusion). Forcing with $\textup{Col}(\kappa, \lambda)$ collapses $\lambda$ to $\kappa$.
\begin{theorem}
Forcing with urelements does not preserve the DC$_{\omega_1}$-scheme in general even if the ground model satisfies ZFCU. 
\end{theorem}
\begin{proof}
Consider a countable transitive model $M$ of $\ZFCUR$ where every set of urelements has tail cardinal $\omega_1$. By Theorem \ref{tail->collection} and Lemma \ref{Tailkappa->DCkappa}, both Collection and the DC$_{\omega_1}$-scheme hold in $M$.  In $M$, let $\P = \textup{Col}(\omega, \omega_1)$ and $G$ be $M$-generic over $\P$. Then in $M[G]$, every set of urelements is countable, because every $A \in M[G]$ is a subset of some $A' \in M$ such that $|A'| \leq \omega_1{^M}$ but $\omega_1{^M}$ is collapsed to $\omega$ in $M[G]$. As a result, every set of urelements will have tail cardinal $\omega$. By an usual argument as in Theorem \ref{zfcurindependece}, this implies that the DC$_{\omega_1}$-scheme fails in $M[G]$.
\end{proof}
Note that since ZFCU proves the DC$_\omega$-scheme (Theorem \ref{maintheorem1}), forcing over ZFCU preserves DC$_\omega$-scheme as it preserves ZFCU.
\begin{question}
Does forcing over $\ZFCUR$ preserve the DC$_\omega$-scheme?
\end{question}

\begin{lemma}\label{recovercollection}
Let $M$ be a countable transitive model of $\ZFCUR$ where for every set of urelements, there is another infinite disjoint set of urelements. Then there is a forcing extension of $M$ which satisfies ZFCU.
\end{lemma}
\begin{proof}
By Theorem \ref{maintheorem1} and \ref{forcingpreservesZFCU}, we may assume that in $M$, there is a least cardinal $\kappa$ not realized since otherwise Collection holds in every forcing extension of $M$. Let $G$ be an $M$-generic filter over $\textup{Col}(\omega, \kappa)$. As $\kappa$ is collapsed to $\omega$ in $M[G]$, every set of urelements in $M[G]$ is countable. If $A$ is  a set of urelements in $M[G]$, let $A' \in M$ be such that $A \subseteq A'$. By assumption, there is another infinite $B \in M$ disjoint from $A'$. Since $B$ has size $\omega$ in $M[G]$, Tail holds in $M[G]$ and hence Collection holds in $M[G]$ by Lemma \ref{tail->collection}. \end{proof}
The next theorem says that ZFCU (in particular, Collection) is necessarily ``forceble'' when the ground model satisfies $\ZFCUR$ + DC$_\omega$-scheme.
\begin{theorem}
If $M$ is countable transitive model of $\ZFCUR$ + DC$_\omega$-scheme and $M[G]$ is a forcing extension of $M$, then $M[G]$ has a forcing extension that satisfies ZFCU.
\end{theorem}
\begin{proof}
First, we may assume that in $M$ $\A$ is a proper class. By the DC$_\omega$-scheme in $M$, for every set of urelements in $M$, there is an infinite set of urelements disjoint from it. But notice that this fact is preserved by forcing by Lemma \ref{McoversM[G]}. So we can apply Lemma \ref{recovercollection} to $M[G]$.
\end{proof}
\noindent Not every model of $\ZFCUR$ has a forcing extension which satisfies ZFCU. For example, if in $M$ every set of urelements is finite but there is a proper class of them, then this will remain the case in every forcing extension of $M$.

\subsection{Ground model definability}
Laver \cite{Laver2007-LAVCVL} and Woodin \cite{woodin2011continuum} proved independently the ground model definability for ZFC: every model of ZFC is definable in its forcing extensions with parameters. Here I first show that the ground model definability fails badly when the gound model contains a proper class of urelements. I then generalize Laver's argument, which is also attributed to Hamkins \cite{Hamkins2003ExtensionsWT}, to consider when ground model definability will hold.

For any infinite set of $x \in M$, $\textup{Fn}(x, 2)$ is the forcing poset consisting of all finite partial functions from $x$ to $2$ ordered by reversed inclusion. If $G$ is an $M$-generic filter over $\textup{Fn}(x, 2)$, then for every set $y \in M$ that is equinumerous with $x$, $M[G]$ contains a new subset of $y$.
\begin{theorem}
Let $M$ be a countable transitive model of $\ZFUR$.
\begin{enumerate}
    \item If $M \models$ DC$_\omega$-scheme + ``$\A$ is a proper class'', then $M$ has a forcing extension in which $M$ is not definable with parameters;
    \item if $M \models $ Plentitude + AC, then $M$ is not definable in any of its non-trivial forcing extensions. 
\end{enumerate}
\end{theorem}
\begin{proof}
(1) Suppose that $M \models$ DC$_\omega$-scheme + ``$\A$ is a proper class''. Let $\P \in M$ be $\textup{Fn}(\omega, 2)$ and $G$ be an $M$-generic filter over $\P$. Suppose \textit{for reductio} that $M$ is definable in $M[G]$ with a parameter $\dot{u}_G \in M[G]$ such that $M = \{ x \in M[G] : M[G] \models \varphi (x,\dot{u}_G)\}.$ Let $B' \in M$ be an infinite set of urelements disjoint from $ker(\dot{u})$, which exists by the DC$_\omega$-scheme. Since DC$_\omega$ implies that $B'$ must have an $\omega$-subset, $M[G]$ contains a new countable subset $B$ of $B'$ which is not in $M$. Fix another $\omega$-set of urelements $C \in M$ disjoint from $ker(\dot{u}) \cup B'$. In $M[G]$, there will be an automorphism that swaps $C$ and $B$ while point-wise fixing $ker(\dot{u})$. Since $M[G] \models \neg \varphi (B,\dot{u}_G)$ and $ker(\dot{u}_G) \subseteq ker(\dot{u})$, it follows that $M[G] \models \neg \varphi (C,\dot{u}_G)$ and hence $C \notin M$, which is a contradiction.

(2) Suppose that $M \models$ Plentitude$^+$ and consider any $M[G]$ such that $M \subsetneq M[G]$. First observe that there must be some set of urelements $B$ such that $B \in M[G] \setminus M$. Fix some $\dot{x}_G \in M[G] \setminus M$ of the least rank so that $\dot{x}_G \subseteq M$. Let $A = ker(\dot{x})$. It follows that $\dot{x}_G \subseteq V_\alpha (A)^M$ for some $\alpha$. By Plenitude and AC in $M$, there is a bijection $f$ from $V_\alpha (A)^M$ to a set of urelements. $f[\dot{x}_G]$ will then be a new set of urelements in $M[G]$.

For \textit{reductio}, suppose that $M = \{ x \in M[G] : M[G] \models \varphi(x,  \dot{u}_G)\}$ for some formula $\varphi$ with parameter $ \dot{u}_G$. Fix some $B \in M[G]\setminus M$ and $B' \in M$ such that $B \subseteq B'$. In $M$, $B'$ has a duplicate $E$ that is disjoint from $ker(\dot{u})$. Then $E$ has a new subset $D$ in $M[G]$ that is disjoint from $ker(\dot{u})$. By AC and Plenitude in $M$, we can again find a duplicate $C \in M$ of $D$ that is disjoint from $ker(\dot{u})$. So there will be an automorphism in $M[G]$ that swaps $C$ and $D$ while point-wise fixing $ker(\dot{u})$. As $M[G] \models \neg \varphi (D, \dot{u}_G)$, it follows that $M[G] \models \neg \varphi(C, \dot{u}_G)$ and hence $C \notin M$, which is a contradiction. 
\end{proof}

The ground model definability holds for models of ZFCU when there is only a set of urelements. This is proved by making some adjustments to the ZFC arguments as in \cite{Laver2007-LAVCVL} and \cite{Reitz2007-JONTGA}. I shall include a proof of this for completeness.
\begin{definition}
 Let $\mathcal{L}_\delta$ be the language extending the language of urelement set theory with a constant symbol $\delta$. ZFCU$_\delta$ is the theory in $\mathcal{L}_\delta$ consisting of ZU, AC, the axiom that $\delta$ is a regular cardinal, $\leq \delta$-Replacement, which states that Replacement holds for every definable function with a domain of size less than $\delta$, and the following axiom.
\begin{itemize}
    \item [](*) Every well-ordering is isomorphic to $\langle \alpha, \in \rangle$ for some ordinal $\alpha$.
\end{itemize}
\end{definition}
\noindent Given a set of urelements $A$, a regular cardinal $\delta$ and a Beth-fixed point $\lambda$ with cf$(\lambda) > \delta$, $V_\lambda(A) \models$ ZFCU$_\delta$.
\begin{lemma}\label{M=M'}
Let $M$ and $M'$ be two transitive models of ZFCU$_\delta$. If $M$ and $M'$ have the same sets of ordinals and the same sets of urelements, then $M = M'$.
\end{lemma}
\begin{proof}
Fix any set $x \in M$. Since $M \models$ZFCU$_\delta$, every set in $M$ has a transitive closure and hence a kernel. By AC in $M$, there is some $\alpha$ and $R \subseteq \alpha \times \alpha$ such that $\langle \alpha, R \rangle$ is isomorphic to $\langle trc(\{x\}), \in \rangle$. Let $W$ well-order $(\alpha +1  )\times (\alpha + 1)$ by G\"odel's pairing function, which will be isomorphic to some ordinal $\beta$ by (*). Let $G$ be the injective function on $\alpha \times \alpha$ such that $G(\eta, \zeta)$ is the order-type of the initial segment of $W$ up to $\langle \eta, \zeta \rangle$. Now we can code $\langle \alpha, R \rangle$ into a set of ordinals $\overline{x} = \{G(\alpha, G(\eta, \zeta)) : \langle \eta, \zeta \rangle \in R\}$. Then $x$ is coded by some $\langle \overline{x}, ker(x) \rangle$, which is also in $M'$ by assumption. Since the definition of $G$ is absolute, we can then decode $\langle \overline{x}, ker(x) \rangle$ in $M'$ and hence $x \in M'$.\end{proof}

\begin{definition}[Hamkins]
Let $N \subseteq M$ be transitive models of $\ZFCUR$ and $\delta$ is a regular cardinal in $M$. 
\begin{enumerate}
    \item ( $\<N, M>$ has the $\delta$-cover property if and only if for each $ x \in M$ with $x \subseteq N$ and $(|x| < \delta)^M$, there is a $y \in N$ such that $x \subseteq y$ and $(|y| < \delta)^N$.
    \item $\<N, M>$ has the $\delta$-approximation property if and only if for each $y \in M$ with $y \subseteq N$, if $y \cap z \in N$ for every $z \in N$ with $(|z|<\delta)^N$, then $y \in N$.
\end{enumerate}
\end{definition}
\begin{theorem}\label{Hamkinsuniquenessthm}
Let $N, N' \subseteq M$ be transitive models of ZFCU$_\delta$ for some regular cardinal $\delta$ in $M$. If both $\langle N, M\rangle$ and $\langle N', M\rangle$ have the $\delta$-cover and $\delta$-approximation properties and $(\delta^+)^N = (\delta^+)^{N'} = (\delta^+)^M$ and $P(\delta)^N = P(\delta)^{N'}$, then $N$ and $N'$ have the same sets of ordinals.
\end{theorem}
\begin{proof}
The proof is the same as in Laver \cite[Theorem 1]{Laver2007-LAVCVL}, attributed to Hamkins.
\end{proof}
\noindent When $N$ and $N'$ are transitive models of ZFC$_\delta$, we will in fact have $N' = N$ as every pure set can be coded by a set of ordinals. But it is possible for $N$ and $N'$ to have the same sets of ordinals but different urelements (or sets of urelements).
\begin{lemma}\label{M[G]coversM}
Let $M$ be a transitive model of $\ZFCUR$, $\P \in M$ be such that $(|\P| < \delta)^M$ for some regular cardinal $\delta$ in $M$ and $G$ be an $M$-generic filter over $\P$. Then $\langle M, M[G] \rangle$ has the $\delta$-cover and $\delta$-approximation properties. Moreover, for any set of urelemetns $A \in M$, there are unboundedly many Beth-fixed points $\lambda$ such that $\langle V_\lambda(A)^M, V_\lambda(A)^{M[G]} \rangle$ has $\delta$-cover and $\delta$-approximation properties.
\end{lemma}
\begin{proof}
For the $\delta$-cover property, suppose that $x \in M[G]$, $x \subseteq M$ and $(|x| < \delta)^{M[G]}$. Then $ker(x) \subseteq A$ for some $A \in M$. So $x \subseteq V_\alpha (A)^M$ for some $\alpha$. In $M[G]$, fix a bijection $f$ from $\kappa$ to $x$ for some $\kappa < \delta$. Since $\P \text{ has } \delta\text{-}c.c.$,\footnote{A forcing poset $\P$ has $\kappa\textup{-}c.c.$ if every antichain in $\P$ has size less than $\kappa$.} by the standard argument there is some function $g \in M$ from $\kappa$ to $P(V_\alpha (A)^M)$ such that $f(\beta) \in g(\beta)$ and $|g(\beta)| < \delta$ for every $\beta < \kappa$. Then in $M$, $\bigcup_{\beta < \kappa} F(\beta)$ has size $< \delta$ and covers $x$. 

For the $\delta$-approximation property, by following Laver's argument in \cite{Laver2007-LAVCVL} we can first show that the $\delta$-approximation property holds sets of ordinals in $M[G]$. Fix an $x \in M[G] \setminus M$ with $x \subseteq M$. Since $x$ is a subset of some $V_\lambda(A)^M$ for some $A \in M$, in $M$ we can fix a bijection $f$ between some ordinal $\alpha$ and $V_\lambda(A)^M$. $f[x]$ is then a set of ordinals not in $M$, so there is a set of ordinals $z \in M$ with $(|z| < \delta)^M$ such that $f[x] \cap  z \notin M$. Using $f$ and $z$ we can find a $y \in M$ such that $(|y| < \delta)^M$ and $y\cap x \notin M$. Hence, $\langle M, M[G] \rangle$ has the $\delta$-approximation property.

Moreover, fix some $\alpha$ and let $\lambda > \alpha$ be a Beth-fixed point with $ \textup{cf} (\lambda) > \delta$. We show that $ V_\lambda(A)^M$ and $V_\lambda(A)^{M[G]}$ have the $\delta$-cover and $\delta$-approximation properties. For the $\delta$-cover property, if $x \in V_\lambda(A)^{M[G]}$ is such that $(|x| < \delta)^{M[G]}$ and $x \subseteq  V_\lambda(A)^M$, then $x \in V_\beta (A)^{M[G]}$ for some $\beta < \lambda$ and so  $x \subseteq V_\beta (A)^M$. By the $\delta$-cover property of $\langle M, M[G] \rangle$, $x$ is covered by some $y \in M$ with $(|y| < \delta)^M$. As a result, $y \cap V_\beta (A)^M$ covers $x$ in $V_\lambda(A)^M$.

For the $\delta$-approximation property, let $x \in V_\lambda(A)^{M[G]}$ be such that $x \subseteq  V_\lambda(A)^M$ and suppose that for every $y \in  V_\lambda(A)^M$  with $(|y| < \delta)^M$, $x \cap y \in V_\lambda(A)^M$. For any $\overline{y} \in M$ with cardinality less than $\delta$, as cf$(\lambda) > \delta$, $\overline{y} \cap V_\lambda(A)^M \in V_\lambda(A)^M$; so $\overline{y} \cap x = \overline{y} \cap V_\lambda(A)^M \cap x$ is in $ V_\lambda(A)^M$ and hence in $M$. By the $\delta$-approximation property of $\langle M, M[G] \rangle$, this implies that $x \in M$. Since $x \subseteq  V_\beta (A)^M$ for some $\beta < \lambda$, $x \in V_\lambda(A)^M$. \end{proof}

\begin{theorem}
Let $M$ be a transitive model of ZFCU + ``$\A$ is a set'', $\P \in M$ and $G$ be an $M$-generic filter over $\P$. $M$ is definable in $M[G]$ with parameters.
\end{theorem}
\begin{proof}
In $M[G]$, let $\A$ be the set of all urelements, $\gamma = |\P|$ and $\delta$ be $\gamma^+$. For any set $N$ and cardinal $\lambda$ in $M[G]$, $N$ is said to be is \textit{good for} $\lambda$ if 
\begin{itemize}
    \item [] (i) $N$ is transitive and $\langle N, \in \rangle \models$ZFCU$_\delta$;
    \item [] (ii) $\delta^+{^N} = \delta^+{^{M[G]}}$ and $P(\delta)^N = P(\delta)^M$; 
    \item [] (iii) $\lambda = Ord \cap N$;
    \item [] (iv) $\langle N, V_\lambda (\A)^{M[G]} \rangle \text{ has the } \delta \text{-cover and } \delta \text{-approximation properties}$.
\end{itemize}
For any cardinal $\lambda$, we say $\lambda$ is \textit{tall} if $\lambda$ is a Beth-fixed point with cf$(\lambda) > \delta$.
Let $\varphi (x, P(\A)^M, P(\delta)^{M})$ be the following formula with parameters $P(\A)^M$ and $P(\delta)^{M}$.
\begin{itemize}
    \item [] $\varphi (x, P(\A)^M, P(\delta)^{M})$ if and only if $\exists N, \lambda ( x \in N \land N \text{ is good for } \lambda \land \lambda \text{ is tall} \land P(\A)^N = P(\A)^M )$ 
\end{itemize}
We claim that $M = \{x \in M[G] : M[G] \models \varphi (x, P(\A)^M, P(\delta)^{M})\}$. 

If $x \in M$, then $x$ is in some tall $\lambda$ and so $V_\lambda(\A)^M \models$ZFCU$_\delta$. $\langle V_\lambda(\A)^M , V_\lambda (\A)^{M[G]} \rangle$ has the $\delta$-cover and $\delta$-approximation properties by Theorem \ref{M[G]coversM}, and $P(\A) \cap V_\lambda(\A)^M = P(\A)^M$. Therefore,  $M[G] \models \varphi (x, P(\A)^M, P(\delta)^{M})$.

Suppose that $M[G] \models \varphi (x, P(\A)^M, P(\delta)^{M})$. Then $x$ is in some $N$ such that $N$ is good for $\lambda$, so $N$ and $V_\lambda (\A)^{M[G]}$ have the $\delta$-cover and $\delta$-approximation properties. Since $\lambda$ is tall, $V_\lambda(\A)^M$ and $V_\lambda (\A)^{M[G]}$ satisfy ZFCU$_\delta$ and have the $\delta$-cover and $\delta$-approximation properties. $V_\lambda(\A)^M$ computes $\delta^+$ correctly with respect to $V_\lambda (\A)^{M[G]}$ because $\P$ has $\delta$-c.c. in $M$ and $\lambda$ is tall; it is also clear that $P(\delta)^{V_\lambda(\A)^M} = P(\delta)^M$. Theorem \ref{Hamkinsuniquenessthm} applies, so $V_\lambda(\A)^M$ and $N$ have the same sets of ordinals. Since they also have the same sets of urelements, it follows from Lemma \ref{M=M'} that $V_\lambda(\A)^M = N $. Therefore, $x$ is in $M$. \end{proof}

\begin{corollary}
Let $M$ be a transitive model of $\ZFCUR$ where some cardinal $\kappa$ is not realized. If $\P \in M$ is such that $(|\P|$ is $\kappa$-closed $)^M$ and $G$ is $M$-generic over $\P$, $M$ is definable in $M[G]$.
\end{corollary}
\begin{proof}
First, $M$ and $M[G]$ have the same sets of urelements because a $\kappa$-closed forcing does not add small subsets. Hence, $P(A)^M = P(A)^{M[G]}$ for every $A \in M$.
Let $\delta$ be a regular cardinal in $M$ such that $(|\P| < \delta)^M$. Then we can verify that $M = \{x \in M[G] : M[G] \models \exists A \varphi (x, P(A), P(\delta)^{M})\}$, where $\varphi (x, P(A), P(\delta)^{M})$ is defined as in the last theorem with parameter $P(\delta)^M$.
\end{proof}

\subsection{$M[G] = M[G]_\#$}\label{subsection:m[g]=m[g]sharp}
Finally, I show that the two ways of defining $\P$-names, as in Definition \ref{oldpnames} and \ref{newpnames}, give rise to the same forcing extension when we force over a countable transitive model of $\ZFCUR$.

\begin{definition}
Let $M$ be a countable transitive model of $\ZFCUR$ and $\P \in M$ be a forcing poset. We define a map $\sim : M^\P_\# \rightarrow M^\P$ by recursion as follows. For every $\dotx \in M^\P_\#$,
\begin{equation*}
    \widetilde{\dotx} =
    \begin{cases*}
     \{\<a, 1_\P>\} & if $\A(\dotx)$  \\
     \{\<\widetilde{\doty}, p> : \<\doty, p> \in \dot{x}\}        & otherwise 
    \end{cases*}
  \end{equation*}
\end{definition}

\begin{lemma}\label{tilde1-1}
For any $M$-generic filter $G$ over $\P$ and $\dot{y}, \dot{x} \in M^\P_\# $, $\dot{y}_G = \dot{x}_G$ if and only if $\widetilde{\dot{y}}_{G} = \widetilde{\dot{x}}_{G}$.
\end{lemma}
\begin{proof}
We prove it by an induction on the rank of $\dot{y}$ and $\dot{x}$. The lemma holds easily when $\dot{y}$ and $\dot{x}$ are urelements. Suppose $\dot{x}$ and $\dot{y}$ are sets. Then $\widetilde{\dot{x}}$ and $\widetilde{\dot{y}}$ don't contain any urelements in their domains, so $\widetilde{\dot{y}}_{G}$ and $\widetilde{\dot{x}}_{G}$ must be sets. If $\dot{y}_G = \dot{x}_G$, then for any $\widetilde{\dotz}_{G}  \in \widetilde{\dot{y}}_{G}$, $\dotz_G \in \dot{x}_G$ so $\dotz_G = \dot{v}_G$ for some $\dot{v} \in dom(\dot{x})$; by the induction hypothesis, $\widetilde{\dotz}_{G} = \widetilde{\dot{v}}_{G} \in \widetilde{\dot{x}}_{G}$ so $\widetilde{\dot{y}}_{G} \subseteq \widetilde{\dot{x}}_{G}$ and hence $\widetilde{\dot{y}}_{G} = \widetilde{\dot{x}}_{G}$ by the same argument. If $\widetilde{\dot{y}}_{G} = \widetilde{\dot{x}}_{G}$, then for any $\dotz_G \in \dot{y}_G$, $\widetilde{\dotz}_{G} \in \widetilde{\dot{x}}_{G}$ so $\widetilde{\dotz}_{G} = \widetilde{\dot{v}}_{G}$ for some $\dot{v} \in dom(\dot{x})$; by the induction hypothesis, $\dotz_G = \dot{v}_G \in \dot{x}_G$; so $\dot{y}_G \subseteq \dot{x}_G$ and hence the same argument shows that $\dot{y}_G = \dot{x}_G$.
\end{proof}
\noindent The next lemma shows that every $\P$-name in $M^\P$ is a mixture of the $\sim$-image of some $\P$-names in $M^\P_\#$.

\begin{lemma}\label{MPbarmixMP}
Let $M$ be a countable transitive model of $\ZFCUR$ and $\P \in M$ be a forcing poset. For every $\P$-name $\dot{x}$ in $M^\P$,  there is a function $f : dom(f) \rightarrow M^\P_\# $ in $M$ such that (i) $ker(f) \subseteq ker(\dot{x}) \cup ker(\P)$; (ii) $dom(f)$ is a maximal antichain of $\P$; and (iii) for every $p \in dom(f)$, $p \forces \dot{x} = \widetilde{f(p)}$.
\end{lemma}
\begin{proof}
By induction on the rank of $\dot{x}$. Suppose the lemma holds for all the $\P$-names in the domain of $\dot{x}$. Condition (i) allows us to find (without using Collection) some $\alpha$ that is big enough such that for every $\doty \in dom(\dot{x})$, there is some $f$ as in the lemma that lives in $V_\alpha(ker(\dot{x}) \cup ker(\P))$. Then by AC in $M$, we can choose an $f_{\doty}$ for each $\doty \in dom(\dot{x})$. In $M$, define
\begin{align*}
    \dot{w} = \{\<f_{\doty} (p), r> : \doty \in dom(\dot{x}) \cap M^\P \land \exists q (\<\doty, q> \in \dot{x} \land p \in dom(f_{\doty}) \land r \leq p, q)\}.
\end{align*}
It is clear that $\dot{w} \in M^\P_\# $ and $ker(\dot{w}) \subseteq ker(\dot{x}) \cup ker(\P)$. Define $Z = \{ p \in \P : \exists a, q \in \mathcal{A} (\<a, q> \in \dot{x} \land p \leq q) \}$. Let $Y$ be a maximal antichain in $Z$ and let $X$ be a maximal antichain in $\P$ extending $Y$. Note that for every $p \in Y$, there is a unique urelement $a_p \in dom(\dot{x})$ such that $p \leq q$ and $\<a_p, q> \in \dot{x}$ for some  $q$. Now we define $f: X \rightarrow (\mathcal{A} \cap dom(\dot{x})) \cup \{\dotw \}$ as follows.
 \begin{equation*}
    f(p) =
    \begin{cases*}
       a_p & if $p \in Y$  \\
       \dot{w}    & otherwise 
    \end{cases*}
\end{equation*}
It is clear that $ker(f) \subseteq ker(\dot{x}) \cup ker(\P)$. 

It remains to show that for every $p \in X$, $p \forces \dot{x} = \widetilde{f(p)}$. Fix a $p \in X$ and an $M$-generic filter $G$ over $\P$ that contains $p$.

Case 1: $p \in Y$. Then $\widetilde{f(p)} = \{\<a_p, 1_\P>\}$. And since there is a $q$ such that $\<a_p, q> \in \dot{x}$ and $p \leq q$, it follows that $\dot{x}_{G} = a_p = \widetilde{f(p)}_{G}$.

Case 2: $p \notin Y$.
\begin{claim}
$\dot{x}_{G}$ is a set.
\end{claim}
\begin{claimproof}
Suppose $\dot{x}_{G}$ is an urelement. Then for some urelemen $a$ and $q \in G$, $\<a, q> \in \dot{x}$. So there is a $r$ which extends both $p$ and $q$; as $r \in Z$, there is some $s \in Y$ such that $s$ and $r$ are compatible because $Y$ is maximal in $Z$. But this means that $p$ is compatible with some $s \in Y$, which is a contradiction because $X$ is an antichain.
\end{claimproof}

\noindent Then $\widetilde{f(p)}_{G} = \widetilde{\dot{w}}_{G}$. Note that $\widetilde{\dot{w}}_{G}$ is a set by the construction of $\widetilde{\dot{w}}$. So it remains to show that $\dot{x}_{G} \subseteq \widetilde{\dot{w}}_{G}$ and $\widetilde{\dot{w}}_{G} \subseteq \dot{x}_{G}$. Consider any $\doty_{G} \in \dot{x}_{G}$ with $\<\doty, q> \in \dot{x}$ and $q \in G$. Since $dom(f_{\doty})$ is a maximal antichain, there is some $p' \in dom(f_{\doty})$ and $r \in G$ such that $p' \in G$ and $r \leq q, p'$. So $\<f_{\doty} (p'), r> \in \dot{w}$ and $p' \forces \doty = \widetilde{f_{\doty}(p')}$. It follows that $ \doty_{G} = \widetilde{f_{\doty}(p')}_{G} \in \widetilde{\dot{w}}_{G}$ and hence $\dot{x}_{G} \subseteq \widetilde{\dot{w}}_{G}$.

To show that $\widetilde{\dot{w}}_{G} \subseteq \dot{x}_{G}$, fix some $\widetilde{f_{\doty}(p')}_{G} \in \widetilde{\dot{w}}_{G}$ such that $\doty \in dom(\dot{x})$, $p' \in dom(f_{\doty})$ and  $\<f_{\doty}(p'), r> \in \dot{w}$ for some $r \in G$. Then there is some $q$ such that $\<\doty, q> \in \dot{x}$ and $r \leq p', q$, which implies $\doty_{G} \in \dot{x}_{G}$. As $p' \forces \doty = \widetilde{f_{\doty}(p')}$, we have $\widetilde{f_{\doty}(p')}_{G} = \doty_{G} \in \dot{x}_{G}$, as desired.
\end{proof}

\begin{theorem}
Let $M$ be a countable transitive model of $\ZFCUR$, $\P \in M$ be a forcing poset and $G$ be an $M$-generic filter over $\P$. There is an elementary embedding from $M[G]_\#$ to $M[G]$. Hence, $M[G] = M[G]_\#$.
\end{theorem}
\begin{proof}
 We prove that the map $\dot{x}_G  \mapsto \widetilde{\dot{x}}_{G}$ is elementary by an induction on formulas. Lemma \ref{tilde1-1} shows that this map is well-defined and 1-1. It is easy to check that the map preserves membership. Also, it is clear that $\dot{x}_G$ is an urelement just in case $\widetilde{\dot{x}}_{G}$ is. The Boolean cases are trivial. If $M[G] \models  \exists x \varphi (x)$, then $M[G] \models \varphi (\dotx_{G})$ for some $\dotx \in M^\P$. Fix a function $f$ for $\dotx$ as in Lemma \ref{MPbarmixMP}. Then for some $p \in dom(f)$, $p \in G$ and $p \forces \dot{x} = \widetilde{f(p)}$, and so $\dot{x}_{G} = \widetilde{\dot{y}}_{G}$ where $\dot{y} = f(p) \in M^\P_\# $. By the induction hypothesis, $M[G]_\# \models \varphi (\dot{y}_G)$ and hence $M[G]_\# \models \exists x \varphi(x)$. Therefore, $M[G] \models $ ZFCU$_R$, and by the minimality of $M[G]_\#$ and $M[G]$ it follows that $M[G]_\# = M[G]$.\end{proof}
The assumption $M \models$ AC is not necessary for the conclusion that $M[G]_\# = M[G]$. This is because one can show that $M[G]_\# \models \ZFUR$ whenever $M$ does (the argument is the same as the proof of Theorem \ref{forcingpreservesreplacement}), and so $M[G]_\# = M[G]$ by the minimality of both forcing extensions. However, the proof presented here clarifies the relationship between these two kinds of $\P$-names.

\section{Boolean-valued models with urelements}\label{section:forcingBVM}
In this section, I discuss Boolean-valued models of set theory with urelements. I first summarize some key results regarding Boolean-valued models of $\ZFCUR$ proved in my joint work with Wu \cite{wu2022}. Based on these results, I further investigate Boolean-valued ultrapowers of models of ZFCU. Basic knowledge of Boolean-valude models of ZFC, which is covered in the first three chapters of \cite{LBell2005-LBESTB}, will be assumed.

\subsection{An overview of $\UB$}
If $V$ is a model of ZFC, then given a complete Boolean algebra $\B \in V$, by transfinite recursion in $V$ a $\B$-name is defined to be a function from a set of $\B$-names to $\B$. As in Definition \ref{oldpnames}, there is also a straightforward generalization of $V^\B$, adopted in  \cite{blass1989freyd}, in urelement set theory. Namely, we treat each urelement as its own $\B$-name. And given what have seen earlier, it is unsurprising that this approach faces the same problem as $\P$-names$_\#$.
\begin{definition}
    A Boolean-valued model for a language $\mathscr{L}$ is said to be \textit{full} just in case for any formula $\varphi$ in $\mathscr{L}$ and $\tau_1, ..., \tau_n \in M^\B$, there is some $x \in M^\B$ such that $\llbracket \exists v \varphi(v, \tau_1, ..., \tau_n) \rrbracket = \llbracket  \varphi(x, \tau_1, ..., \tau_n) \rrbracket$.  
\end{definition}
\noindent It is straightforward to show that if urelements are treated as their own names, almost all Boolean-valued models of ZFCU are not \textit{full} (see \cite{wu2022}). However, fullness, as I have mentioned in the beginning of this chapter, is crucial for applications of Boolean-valued models in set theory. In the joint work \cite{wu2022} with Wu, we provided the following new definition of Boolean-valued models with urelements.
\begin{definition}[$\ZFUR$ \cite{wu2022}]\label{newbnames}
Let $\B$ be a complete Boolean algebra.
\begin{enumerate}
    \item A function $\tau: dom(\tau) \rightarrow \mathbb{B}$ is a $\BB$-name if and only if for any $x \in dom(\tau)$, $x$ is either an urelement or a $\BB$-name, and for any urelement $a \in dom(\tau)$ and $x \in dom(\tau)$, $\tau(a) \land \tau(x) = 0_\B$ whenever $x \neq a$.
    
    \item Let $\tau$ be a $\BB$-name. 
    \subitem $dom^{\mathscr{A}}(\tau) = \{ a \in dom(\tau) : a \in \mathscr{A} \}$;
    \subitem $dom^{\mathbb{B}}(\tau) = \{ \eta \in dom(\tau) : \eta \text{ is a $\BB$-name}\}$;
    \subitem $\tau(a) = 0_\B$ whenever $a$ is an urelement not in $dom^{\mathscr{A}}(\tau)$.
    
    \item $\UBB = \{ \tau \in U : \tau \text{ is a } \BB \text{-name}\}$.
    
     \item $\pazocal{L}_{\BB}$ contains $\{\subseteq, =, \in, \A, \overset{\mathscr{A}}{=}\}$ as the non-logical symbols and each $\B$-name as a constant symbol. $\pazocal{AL}_{\BB}$ is the class of all atomic formulas in $\mathcal{L}_{\BB}$. The Boolean evaluation function $\llbracket\ \  \rrbracket: \pazocal{AL}_{\BB} \rightarrow \B$ is defined as follows.
     \begin{align*}
     & \llbracket \mathscr{A}(\tau) \rrbracket = \bigvee\limits_{a \in \mathscr{A}} \tau(a) \\
     & \llbracket \tau \overset{\mathscr{A}}{=} \sigma \rrbracket = \bigwedge\limits_{a \in \mathscr{A}} (\tau(a) \Leftrightarrow \sigma (a) )\\
     & \llbracket \tau \in \sigma \rrbracket = \bigvee\limits_{\mu \in dom^{\mathbb{B}}(\sigma)}\llbracket \tau = \mu \rrbracket \land \sigma(\mu) \\
    & \llbracket \tau \subseteq \sigma \rrbracket= \bigwedge\limits_{\eta \in dom^{\mathbb{B}}(\tau)} \tau(\eta) \Rightarrow \llbracket \eta \in \sigma \rrbracket\\
     & \llbracket \tau = \sigma \rrbracket = \llbracket \tau \subseteq \sigma \rrbracket \land \llbracket \sigma \subseteq \tau \rrbracket \land \llbracket \tau \overset{\mathscr{A}}{=} \sigma \rrbracket
     \end{align*}
     \item $\llbracket\ \  \rrbracket$ is extended to $\mathcal{L}_\B$ in the standard way. Namely,
\begin{enumerate}
    \item [] $\llbracket \varphi \land \psi \rrbracket = \llbracket \varphi \rrbracket\land \llbracket \psi \rrbracket$;
    \item [] $\llbracket \neg \varphi \rrbracket = \neg \llbracket \varphi \rrbracket$;
    \item []  $\llbracket \exists x \varphi \rrbracket = \bigvee\limits_{\tau \in \UB}\llbracket \varphi(\tau) \rrbracket$.
\end{enumerate}
\end{enumerate}
\end{definition}
\noindent Note that no urelement is a $\BB$-name since every $\BB$-name is a set, and each urelement $a$ will be represented canonically by $\{\<a, 1_\B>\}$ in $\UBB$ instead of itself. For any $\tau \in \UBB$ and urelement $a$, $\tau(a)$ will be the $\B$-degree to which $\{\<a, 1_\B>\}$ \textit{is identical }to $\tau$, as opposed to the value of $\{\<a, 1_\B>\}$'s membership to $\tau$. This motivates the incompatibility condition: if $a, b \in dom(\tau)$ are two urelements, then $\tau(a) \land \tau(b)$ must be $0_\B$ because this is the degree to which $\tau$ is both of them; if $a, \sigma \in dom(\tau)$, where $\sigma$ is a $\BB$-name, then $\tau(a) \land \tau(\sigma)$ must be $0_\B$ as well because this is the degree to which $\tau$ is an urelement with a member. In fact, it is this restriction that ensures $\llbracket \text{no urelement has any members} \rrbracket =1_\B$. To see this, consider any $\tau \in \UB$. Since for any urelement $a$ and $\B$-name $\mu \in dom(\tau)$, $\tau(a) \leq \neg \tau(\mu)$, we have
\begin{align*}
   \llbracket \mathscr{A}(\tau) \rrbracket & = \bigvee\limits_{a \in \mathscr{A}} \tau(a) \\
    & \leqslant \bigwedge\limits_{\sigma \in \UB}\bigwedge\limits_{\mu \in dom^{\mathbb{B}}(\tau)} \llbracket \sigma \neq \eta \rrbracket \lor \neg\tau(\mu) \\
    & = \llbracket \forall y (y \notin \tau) \rrbracket.
\end{align*}
Finally, $ \llbracket \tau \overset{\mathscr{A}}{=} \sigma \rrbracket^{\UBB}$ is the degree to which $\tau$ and $\sigma$ are identical when they are taken as urelements.
\begin{theorem}[\cite{wu2022}]\label{collection<->fullness}
Over $\textup{ZFCU}_\text{R}$, the following are equivalent.
\begin{enumerate}
    \item Collection.
    \item For every complete Boolean algebra $\mathbb{B}$, $\UBB$ is full.\qed
\end{enumerate}
\end{theorem}
\begin{theorem}[The Fundamental Theorem of $\UB$ [\cite{wu2022}]]\label{fundamentalthm}
Let $U$ be a model of $\ZFCUR$ and $\B$ be a complete Boolean-algebra in $U$. Then ($\UB \models \varphi$ abbreviates $\llbracket \varphi \rrbracket = 1_\B$)
\begin{enumerate}
    \item $\UB \models $ $\ZFCUR$;
    \item $\UB \models $ Collection if $U \models$ Collection.\qed
\end{enumerate}
\end{theorem}

\subsection{Boolean ultrapowers with urelements}
In this section, I consider how the construction of Boolean ultrapowers studied in \cite{hamkins2012well} can be carried out for models of ZFCU. Let $U$ be any model of ZFCU and $\B$ be a complete Boolean algebra in $U$. Theorem \ref{collection<->fullness} and \ref{fundamentalthm} allow us to transform $U^\B$ into a classical two-valued model with respect to \textit{any} ultrafilter $F$ on $\B$. For every $\sigma, \tau \in \UBB$, we define:
\begin{itemize}
    \item [] $\sigma =_F \tau$ if and only if  $\llbracket \sigma =\tau \rrbracket \in F$;
    \item []  $\sigma \in_F \tau$ if and only if  $\llbracket \sigma \in \tau \rrbracket \in F$;
    \item [] $\A_F(\sigma)$ if and only if $\llbracket \A(\sigma) \rrbracket \in F$.
\end{itemize}
For each $\tau \in U^\B$, let $[\tau]_F$ be the equivalence class $\{\sigma \in \UBB : \sigma =_F \tau \}$ and $\UBB/F = \{[\tau]_F : \tau \in \UBB \}$. $\in_F$ and $\A_F$ are well-defined on the corresponding equivalent classes because $=_F$ is a congruence with respect to them. This generates a two-valued model $\<\UBB/F, \A_F, \in_F>$ for the language of urelement set theory (which will also be denoted by $\UBB/F$.)

\begin{theorem}[Łoś Theorem]\label{losthm}
Let $U$ be any model of ZFCU, $\B$ be a complete Boolean algebra in $U$, and $F$ be an ultrafilter on $\B$. For every $\tau_1, ..., \tau_n \in \UB$, $\UBB/F \models \varphi([\tau_1]_F, ..., [\tau_n]_F)$ if and only if $\llbracket \varphi(\tau_1, ..., \tau_n) \rrbracket \in F$. Hence, $\UB/F \models$ ZFCU.
\end{theorem}
\begin{proof}
Atomic formulas and Boolean connectives are easy to check. For the quantifier case, suppose $\varphi$ is some $\exists x \psi$. The left-to-right dIrection is immediate. If $\llbracket \exists x \psi (x, \tau_1, ..., \tau_n) \rrbracket$ is in $F$, since  $\UBB$ is full by Theorem \ref{collection<->fullness}, it follows that $\llbracket \psi (\tau, \tau_1, ..., \tau_n) \rrbracket \in F$ for some $\tau \in \UB$ and hence $\UBB/F \models \exists x \psi (x, [\tau_1], ..., [\tau_n])$ by the induction hypothesis. $\UBB/F \models$ ZFCU by Theorem \ref{fundamentalthm}.
\end{proof}
\noindent Therefore, to show that a certain statement $\varphi$ is not provable from ZFCU it suffices to find a $\B$ such that $\llbracket \varphi \rrbracket \neq 0_\B$ and then consider an ultrafilter $F \subseteq \B$ that contains $\llbracket \varphi \rrbracket$. And one can further establish the mutual interpretability of various extensions of ZFCU. However, we should note that since forcing preserves Plenitude, this way of taking quotient models directly from a model of ZFCU is much less flexible than the method described in Theorem \ref{zfcumutualinter}.

Following \cite{hamkins2012well}, let us further consider how $U$ sits inside $U^\B$. Define $\Ucheck$ as a new unary predicate such that $\llbracket \tau \in \Ucheck \rrbracket = \bigvee\limits_{x \in U}\llbracket \tau =  \xcheck \rrbracket$, which clearly obeys the law of identity. Define $\tau \in_F \Ucheck_F$ as $\llbracket \tau \in \Ucheck \rrbracket \in F$. $=_F$ remains a congruence with respect to $\Ucheck_F$. In $\UB$, $\Ucheck$ represents the class of all objects the ground model $U$. For any formula $\varphi$ in the language of urelement set theory, let $\varphi^{\Ucheck}$ be the result of restricting all the quantifiers to $\Ucheck$. Then by the same argument as in \cite[Lemma 4]{hamkins2012well}, one can show that $U \models \varphi(x_1, ..., x_n)$ if and only if $\llbracket \varphi^{\Ucheck} (\xcheck_1, ..., \xcheck_n) \rrbracket = 1_\B$. Let $\Ucheck_F = \{[\tau]_F : \llbracket \tau \in \Ucheck \rrbracket \in F\}$. As a submodel of $U^\B/F$, by the same argument as in \cite[Lemma 12]{hamkins2012well}, one can show that $\Ucheck_F \models \varphi ([\tau_1]_F, ..., [\tau_n]_F)$ if and only if $\llbracket \varphi^{\Ucheck}(\tau_1, ..., \tau_n)\rrbracket \in F$. This yields the Boolean ultrapower embedding.

\begin{theorem} \label{ultrapowermap}
Let $U$ be a model of ZFCU, $\B$ be a complete Boolean algebra in $U$ and $F$ be an ultrafilter on $\B$. The Boolean ultrapower map $j_F : U \rightarrow \Ucheck_F $ defined by $j_F(x) = [\xcheck]_F$ is an elementary embedding. \qed
\end{theorem}
\noindent Note that the map $j_F$ is not necessarily onto: the same argument in \cite[Theorem 16]{hamkins2012well} will show that $j_F$ is onto just in case $F$ is $U$-generic.

Furthermore, following \cite{hamkins2012well}, we can show that the quotient structure $U^\B/F$ is in fact a forcing extension of $\Ucheck_F$ in the sense of Definition \ref{m[g]def}. Let $\dot{G} = \{\<\check{p}, p> : p \in \B\}$. It is not hard to check that (see \cite[Lemma 8]{hamkins2012well})
$$U^\B \models \dot{G} \text{ is a } \Ucheck\text{-generic ultrafilter on } \check{\B}.$$
If $\tau \in U^\B$ and $F$ is a filter on $\B$, we can define the $F$-valuation of $\tau$, $\tau_F$, as follows (similar to Definition \ref{m[g]def}). 
\begin{itemize}
    \item [] $\tau_F = a$ if $\mathcal{A}(a)$ and  $\langle a, p \rangle \in \tau$ for some $p \in F$;
    \item [] $\tau_F = \{ \sigma_F: \langle \sigma , p \rangle \in \tau \text{ for some } \sigma \in U^\B \text{ and } p \in F \}$ otherwise.
\end{itemize}
 \noindent Note that this is well-defined by the incomptability condition in Definition \ref{newbnames}. Then by induction one can check that for any $\tau \in U^\B$,
$$\UB \models \tau \text{ is the } \dot{G}\text{-valuation of } \check{\tau}.$$
In other words, $\UB$ sees itself as the forcing extension $\Ucheck[\dot{G}]$. Let $G$ be $[\dot{G}]_F$. By Theorem \ref{losthm}, $\UB/F \models \forall x \exists y \in \check{U}_F(x = y_G)$ and hence $\UB/F=\Ucheck_F[G]$. 

I conclude this chapter with an application of this machinery.
\begin{theorem}
If $U$ is a model of ZFCU + $\neg$Plenitude, then there is a model $W$ of ZFCU in which every set of urelements is countable, and $W$ is a forcing extension of an elementary extension of $U$.
\end{theorem}
\begin{proof}
Let $\kappa$ be the least cardinal not realized in $U$ and $\B = RO(\kappa^{\omega})$ be the complete Boolean-algebra which consists of all the regular open sets of the product topology $\kappa^\omega$, where $\kappa$ is assigned the discrete topology. It is a classical result that $\UB \models \check{\kappa} \sim \omega$. For every $\tau \in \UB$, observe that $\llbracket \tau \subseteq \A \rrbracket = \llbracket  \tau \subseteq \check{A_\tau} \rrbracket$ and $\UB \models |\check{A_\tau}| \leq |\check{\kappa}|$, where $A_\tau$ is the set of urelements in $dom(\tau)$. It follows that $\llbracket \tau \subseteq \A \rrbracket \leq \llbracket |\tau| \leq \omega \rrbracket$ for every $\tau \in \UB$, i.e., $\UB \models $ ``every set of urelements is countable''. Let $F$ be an ultrafilter on $\B$. By Theorem \ref{losthm}, it follows that $\UB/F \models $ ZFCU + ``every set of urelements is countable''. And by the previous disucssion, $\UB/F$ is a forcing extension of $\Ucheck_F$, which is an elementary extension of $U$ by Theorem \ref{ultrapowermap}.
\end{proof}
\chapter{Class Theory with Urelements}
This chapter studies class theory with urelements. In Section \ref{section:classwithurs}, I first introduce the axioms for class theory and show how urelement class theory can be interpreted in pure class theory. Then by generalizing the construction of class permutation models due to Felgner \cite{felgner1976choice}, I isolate a hierarchy of axioms in class theory with urelements. Section \ref{section:RP2withurelements} is concerned with the second-order reflection principle (RP$_2$) in Kelley-Morse class theory with urelements. I first show that when there are no more urelements than the ordinals, RP$_2$ with urelements is bi-interpretable with RP$_2$ in pure class theory. I then introduce a new form of accumulative hierarchy, $U_{\kappa, A}$, and prove a generalized version of Zermelo's Quasi-Categoricity Theorem with urelements. Assuming the consistency of a $\kappa^+$-supercompact cardinal, I construct a $U_{\kappa, A}$-model of RP$_2$ where the urelements are more numerous than the pure sets. At the end, I discuss how this result, together with my recent joint work with Hamkins \cite{HamkinsForthcoming-HAMRIS}, might challenge the doctrine of limitation of size.

\section{Class theory with urelements}\label{section:classwithurs}

\subsection{Axioms} \label{subsection:axiomsinclasstheory}
The language of \textit{ urelement class theory} is a two-sorted language extending the language of urelement set theory, with the first-order variables quantifying over sets and urelements, and the second-order variables quantifying over classes. Proper classes are classes that are not co-extensional with any set. A model for this language is of the form $\<M, \A^M, \in^M, \mathscr{M}>$, where the first-order part $\<M, \A^M, \in^M >$ is a model for the lanauge of urelement set theory and the second-order part $\mathscr{M} \subseteq P(M)$ serves as the domain for classes. In class theory, the axiom schemes in first-order set theory can be now written as single axioms.
\begin{itemize}
    \item [] (Separation) For every class $Y$ and set $x$, $Y \cap x$ is a set.
    \item [] (Replacement) If $F$ is a class function on a set $x$, $F[x]$ is a set.
    \item [] (Collection) If $R$ is a class relation on a set $w$ such that $\forall x \in w\ \exists y R(x, y)$, then there is a set $v$ such that for every $x \in w$ there is a $y \in v$ such that $R(x, y)$.
\end{itemize}
Similarly, for every $\kappa$, the DC$_\kappa$-scheme can be formulated as a single axiom, which I shall call ``$\kappa$-DC''.
\begin{itemize}
    \item [] ($\kappa$-DC) For every relation $R \subseteq U \times U$ that has no terminal nodes, there exists a function $f$ on $\kappa$ such that $R(f\restriction \alpha, f(\alpha))$ for all $\alpha < \kappa$.
 \end{itemize}
As in Proposition \ref{prop:DCkappaVariants}, $\kappa$-DC is equivalent to the assertion that for every class $X$, if every $s \in X^{<\kappa}$ has some $y \in X$ such that $R(s, y)$, then there is an $f \in X^\kappa$ such that $R(f\restriction \alpha, \alpha)$ for every $\alpha < \kappa$. $Ord$-DC is the statement that $\forall \kappa (\kappa$-DC). The same argument as in Lemma \ref{easyimplication} will show that $Ord$-DC holds when there is only a set of urelements over the theory $\GBUR$, which is defined below. The first-order reflection principle will now allow class parameters.
\begin{itemize}
    \item [] (RP) For every $X_1$, ..., $X_n$, there is a transitive set $t$ extending any given set such that for every $x_1, ..., x_m \in t$, $$\varphi(X_1, ..., X_n, x_1, ..., x_m) \leftrightarrow \varphi^t(X_1 \cap t, ..., X_n \cap t, x_1, ..., x_m),$$ where $\varphi$ contains only first-order quantifiers.
 \end{itemize}

Two standard theories of classes are G\"odel-Bernays class theory (GB) and Kelley-Morse class theory (KM). In addition to ZF for sets, GB adopts the following first-order comprehension axiom.
\begin{itemize}
\item [] (First-order Comprehension) For every formula $\varphi$ which contains only first-order quantifiers (but possibly with class parameters), $\{x : \varphi(x)\}$ is a class.
\end{itemize}
KM extends GB by adopting the full comprehension axiom.
\begin{itemize}
\item [] (Full Comprehension) For every formula $\varphi$, possibly with class parameters, $\{x : \varphi(x)\}$ is a class.
\end{itemize}
\noindent GB and KM  are often formulated with a second-order version of AC. However, in the presence of urelements different formulations of this principle can come apart. We shall discuss the following three versions.
\begin{itemize}
	\item[] (Limitation of Size) All proper classes are equinumerous.\footnote{It is immediate that Limitation of Size is equivalent to the following principle. \begin{itemize}
	    \item [] $X$ is a proper class if and only if it is equinumerous with the universe $U$.
	\end{itemize}
And this principle is precisely what the limitation-of-size conception of size (see \ref{subsection:philosophy}) is asserting, so the terminology here is justified.}
	\item[] (Global Well-Ordering) There is a well-ordering of the universe $U$.
	\item[](Global Choice) There is a class function $F$ such that for every non-empty set $x$, $F(x) \in x$.
\end{itemize}
Standard arguments show that over a suitable theory (such as the theory GBU$_\text{R}$ defined below), Limitation of Size $\rightarrow$ Global Well-Ordering $\rightarrow$ Global Choice. However, neither of the implications can be reversed, as proved in \cite{howard1978independence} (see also Theorem \ref{thm:KMURnvdashCollection} and Lemma \ref{UKA}). For this reason, it is useful to isolate urelement class theories that do not include any second-order choice principles.
\begin{definition}
\ \newline
GBU$_\text{R}$ $=$ ZU + Class Extensionality + Replacement + First-Order Comprehension.\\
KMU$_\text{R} =$ GBU$_\text{R}$ + Full Comprehension.\\
$\GBUR =$ GBU$_\text{R}$ + AC.\\
$\KMUR =$ KMU$_\text{R}$ + AC.\\
GBCU $=$ GBU$_\text{R}$ + Global Well-Ordering.\\
KMCU $=$ KMU$_\text{R}$ + Global Well-Ordering.\\
GBc $=$ $\GBUR + \forall x \neg \A (x) $\\
KMc $=$ $\KMUR + \forall x \neg \A (x)$\\
GBC $=$ GBCU  + $\forall x \neg \A (x) $.\\
KMC $=$ KMCU  + $\forall x \neg \A (x) $.
\end{definition}

Theories with a lowercase c do not include any second-order choice principle but only AC for sets. The subscript R indicates again that the theory is formulated with Replacement. As we will see in \ref{section:independenceinClassTheory}, neither $\omega$-DC nor Collection is provable from $\KMUR$, and Collection does not imply $\omega$-DC over $\KMUR$ (although Collection does imply the DC$_\omega$-scheme in $\ZFCUR$). Note that Limitation of Size was not included as an axiom for either KMCU or GBCU, and its philosophical status will be discussed in \ref{section:amodelofKMU+RP2+notLS}. Finally, Gitman and Hamkins \cite{GitmanHamkins:Kelley-MorseSetTheoryAndChoicePrinciplesForClasses} observed that a robust class theory should include the following principle of Class Choice. 
\begin{itemize}
    \item [] (CC) $\forall Z (\forall x \exists X \varphi(x, X, Z) \rightarrow \exists Y \subseteq U \times U \ \forall x \varphi(x ,Y_x, Z))$, where $Y_x = \{ y : \<x, y> \in Y \}$.
\end{itemize}
That is, if every $x$ has a class witness for some relation $\varphi$, there will be a two-dimensional class $Y$ whose $x$-slice is a class witness for $x$. For any cardinal $\kappa$, $\kappa$-CC is the following restricted version of CC.
\begin{itemize}
    \item [] ($\kappa$-CC) $\forall Z (\forall \alpha < \kappa\ \exists X \varphi(\alpha, X, Z) \rightarrow \exists Y \subseteq \kappa \times U \ \forall \alpha <\kappa\ \varphi(\alpha ,Y_\alpha, Z))$.
\end{itemize}
\noindent It is shown in \cite{GitmanHamkins:Kelley-MorseSetTheoryAndChoicePrinciplesForClasses} that even $\omega$-CC is not provable in KMC and that KMC and KMC + CC are mutually interpretable. 

I shall end this subsectoin with some useful facts, all of which should be well-known.
\begin{prop}\label{prop:GC<->GWO<->LS}
 Over GBU$_\text{R}$ + $\A$ is a set. The following are equivalent.
 \begin{enumerate}
     \item Limitation of Size.
     \item Global Well-Ordering.
     \item Global Choice.
 \end{enumerate}
\end{prop}
\begin{proof}
It is clear that over GBU$_\text{R}$, (1) $\rightarrow$ (2) $\rightarrow$ (3). To show (3) $\rightarrow$ (1), let $F$ be a global choice function. Define $G : Ord \rightarrow U$ as 
\begin{itemize}
    \item [] $G(0) = 0$;
    \item [] $G(\alpha) = F(V_\beta(\A) \setminus G[\alpha])$, where $\beta$ is the least ordinal such that $V_\beta(\A) \setminus G[\alpha]$ is non-empty.
\end{itemize}
It is not hard to verify that $G$ is a bijection and hence a global well-ordering. Now let $X$ be a proper class. As $\A$ is a set, the map from $X$ to $Ord$ that maps each $x \in X$ to its rank (i.e., the least $\beta$ such that $x \in V_\beta(\A)$) must be unbounded, which produces an onto map from $X$ to $Ord$. Since we can well-order $X$, there is an injective map from $Ord$ to X; and by a class version of the Cantor-Bernstein Theorem, it follows that $X$ and $Ord$ are equinumerous. Therefore, all proper classes are equinumerous. \end{proof}

\begin{prop}\label{prop:rp+gbur<->collection+dcs}
The following are equivalent over $\GBUR$.
\begin{enumerate}
    \item RP.
    \item Collection $\land \ \omega$-DC.
\end{enumerate}
\end{prop}
\begin{proof}
 To show (1) $\rightarrow$ (2), assume RP. Suppose that $\forall x \in w\ \exists y\  R(x, y)$ for some class relation $R$ and set $w$. We can then reflect it down to a transitive set extending $w$, which will give us a collection set. And suppose that $R$ is a class relation without terminal nodes; we can then reflect this fact down to some transitive set $t$. The existence of a desired $\omega$-sequence then follows from AC. (2) $\rightarrow$ (1) is proved in the same way as in Theorem \ref{collection+dc->rp} since the class parameters do not raise any difficulties.
\end{proof}

\begin{prop}\label{prop:OrdDC->RP}
  GBU$_\text{R} + Ord$-DC $\vdash$ RP.
\end{prop}
\begin{proof}
Assume $Ord$-DC. By a standard argument as in \cite[Theorem 8.1]{jech2008axiom}, $Ord$-DC implies AC. Also, by Proposition \ref{prop:rp+gbur<->collection+dcs}, it suffices to show that Collection holds. Let $w$ be an infinite set and $R$ be a class relation such that $\forall x \in w\ \exists y\  R(x, y)$. We may enumerate $w$ with $\{x_\alpha : \alpha < \kappa\}$ by some infinite cardinal $\kappa$. Define $X = \{ \<x_\alpha, y> : \alpha < \kappa \land R(x_\alpha, y) \}$. For every $s$ and $z$, define
\begin{itemize}
    \item []  $R^*(s, z)$ if and only if, whenever $s$ is function on $\alpha$ for some $\alpha < \kappa$, $z = \<x_\alpha, y>$ and $z \in X$.
\end{itemize}
 By $Ord$-DC, there is an $f \in X^\kappa$ such that $R^*(f\restriction \alpha, f(\alpha))$ for every $\alpha < \kappa$. Then the transitive closure of $f$ will be a desired collection set.  
\end{proof}

\begin{prop}\label{prop:GBU->RPandDCord}
GBCU $\vdash$ $Ord$-DC $\land$ RP.
\end{prop}
\begin{proof}
Given a class relation without terminal nodes, for every $\kappa$ we can use the global well-ordering to construct a desired $\kappa$-sequence. Then RP holds by Proposition \ref{prop:OrdDC->RP}.
\end{proof}
To conclude, the following diagram holds in $\GBUR$, whose completeness will be discussed in \ref{section:independenceinClassTheory}.
\begin{figure}[hbt!]
\begin{center}
\begin{tikzpicture}
\begin{scope}[every node/.style={}]
    \node (A) at (8, -2.5) {Collection};
    \node (B) at (9, 2){$\A$ is a set};
    \node (C) at (3, 1.5){Global Choice};
    \node (D) at (6,1.5) {$Ord$-DC };
    \node (E) at (6, -0.5) {$\kappa$-DC};
    \node (I) at (8, -0.5) {RP};

    \node (L) at (6, 0.7) {.};
    \node (M) at (6, 0.5) {.};
    \node (N) at (6, 0.3) {.};
    \node (P) at (6, -1.7) {.};
    \node (Q) at (6, -1.5) {.};
    \node (R) at (6, -1.3) {.};
    \node (S) at (6, -2.5) {$\omega$-DC};
    \node (T) at (6, 3) {Global Well-Ordering};
    
\end{scope}

\begin{scope}[>={stealth},
              every node/.style={fill=white,circle},
              every edge/.style={draw=black}]

    \path [->] (T) edge (D);
\path [->] (T) edge (C);
     \path [->] (B) edge (D);

    \path [->] (I) edge (A);
      \path [->] (I) edge (S);
    \path [->] (D) edge (I);
    \path [->] (D) edge (L);
    \path [->] (N) edge (E);
    \path [->] (E) edge (R);
     \path [->] (P) edge (S);
       
\end{scope}
\end{tikzpicture}
 \caption{Implication diagram in $\GBUR$}
 \label{GBUdiagram}
\end{center}
\end{figure}
\FloatBarrier
\subsection{Interpreting $\mathcal{U}$ in $\mathcal{V}$}\label{section:interpretingUrelementsinClassTheory}
The construction of $V\llbracket X \rrbracket$ introduced in Section \ref{section:interpretingUinV} can be easily generalized to interpreting urelement class theory in pure class theory.
\begin{definition}\label{barwiseinterpretation2}
Let $\<V, \in , \mathscr{V}>$ be a model of GBc and $X \in \mathscr{V}$. $\<V \llbracket X \rrbracket, \bar{A}, \bar{\in}>$ is then defined as in Definition \ref{barwiseinterpretation1}. Let $\mathscr{V}\llbracket X \rrbracket = \{Y \in \mathscr{V} : Y \subseteq  V \llbracket X \rrbracket\}$. $\mathscr{V}\llbracket X \rrbracket$ denotes the model $\< V \llbracket X \rrbracket, \bar{A}, \bar{\in}, \mathscr{V}\llbracket X \rrbracket>$.\footnote{Note that when $Y \in \mathscr{V}\llbracket X \rrbracket$, $\mathscr{V}\llbracket X \rrbracket \models x \in Y$ if and only if $\mathscr{V} \models x \in Y$.}
\end{definition}

\begin{theorem}\label{Con(KM)->Con(KMU)}
Suppose $\mathscr{V} \models $ GBc and $X$ is a class in $\mathscr{V}$. Then
\begin{enumerate}
    \item $\mathscr{V}\llbracket X \rrbracket \models $ $\GBUR$ + Collection;
    \item $\mathscr{V}\llbracket X \rrbracket \models $ $\KMUR$ + Collection if $\mathcal{V} \models$ KMc;
    \item $\mathscr{V}\llbracket X \rrbracket \models $ CC if $\mathscr{V} \models$ CC;
    \item $\mathscr{V}\llbracket X \rrbracket \models $ Limitation of Size (and hence GBCU) if $\mathcal{V} \models$ GBC.
\end{enumerate}
\end{theorem}

\begin{proof}

For (1) and (2), $\mathscr{V}\llbracket X \rrbracket \models$ ZU by Theorem \ref{thm:V[X]modelsZFU}, so it remains to show that $\mathscr{V}\llbracket X \rrbracket$ satisfies Collection, Class Extensionality, and First-Order Comprehension (or Full Comprehension), all of which follow easily from the fact that $\mathscr{V} \models$ GBc (or KMc). For example, to show $\mathscr{V}\llbracket X \rrbracket \models$ Collection, suppose that in $\mathscr{V}\llbracket X \rrbracket$ for every $x \in \barw$ there is some $y$ such that $x \bar{R} y$ for some $\barw, \bar{R} \in \mathscr{V}\llbracket X \rrbracket$, where $\barw = \<1, w>$. Then in $V$ there is some $v \subseteq V \llbracket X \rrbracket$ such that for every $\barx \in w$, there is some $\bary \in v$ such that $\mathscr{V}\llbracket X \rrbracket \models \<\barx, \bary> \in \bar{R}$. Then $\bar{v} = \<1, v>$ is a desired collection set in $\mathscr{V}\llbracket X \rrbracket$.

(3) Suppose that $\mathscr{V} \models$ CC AND for every $\barx \in V\llbracket X \rrbracket$, there is some class $\bar{X} \in \mathscr{V}\llbracket X \rrbracket$ with $\varphi(\barx, \bar{X}, \bar{Z})^{ \mathscr{V}\llbracket X \rrbracket}$. By CC in $\mathscr{V}$, there is a class $Y \subseteq V\llbracket X \rrbracket \times U$ such that for every $\barx \in V\llbracket X \rrbracket$, $\varphi(\barx, Y_{\barx}, \bar{Z})^{ \mathscr{V}\llbracket X \rrbracket}$ and $Y_{\barx} $ is a class of $\mathscr{V}\llbracket X \rrbracket$. Define $\overline{Y} = \{\overline{\<\barx, \bary>} : \bary \in Y_{\barx} \land \barx \in \mathscr{V}\llbracket X \rrbracket \}$ (where $\overline{\<\barx, \bary>}$ codes the ordered-pair as in Theorem \ref{thm:ACholdsViffACholdsinV[X]}), which is a class of $\mathscr{V}\llbracket X \rrbracket$. Since for every $\barx \in \mathscr{V}\llbracket X \rrbracket$, $\mathscr{V}\llbracket X \rrbracket \models Y_{\barx} = \overline{Y}_{\barx}$, it follows that CC holds in $\mathscr{V}\llbracket X \rrbracket$.

(4) Suppose that $\mathscr{V}$ Global Well-Ordering. Note that if $\bar{Y}$ and $\bar{Z}$ are two proper classes in $\mathscr{V}\llbracket X \rrbracket$ , they must be equinumerous proper classes in $\mathscr{V}$ by Propostion \ref{prop:GC<->GWO<->LS}. So in $\mathscr{V}$ there must be a bijection $F$ between $\bar{Y}$ and $\bar{Z}$. Then $\bar{F} = \{ \overline{\<\bary, \barz>} : F(\bary) = \barz \}$ will be a bijection between $\bar{Y}$ and $\bar{Z}$ in $\mathscr{V}\llbracket X \rrbracket$.\end{proof}

\begin{theorem}\label{KMU+LSbiintepsKM+LS}
The following pairs of theories are bi-interpretable with parameters.
\begin{enumerate}
    \item GBC and GBCU + $\A \sim \omega$;
    \item KMC and KMCU + $\A \sim \omega$;
    \item GBC and GBCU + Limitation of Size;
    \item KMC and KMCU + Limitation of Size;
    \item GBC and GBCU + Limitation of Size + Plenitude;
    \item KMC and KMCU + Limitation of Size + Plenitude.
\end{enumerate}
\end{theorem}
\begin{proof}
Working in GBC (KMC), we can form either $\mathscr{V}\llbracket \omega \rrbracket$ (or $V\llbracket Ord\rrbracket$), and the map $y \mapsto \hat{y}$ and $Y \mapsto \hat{Y} = \{\hat{y} : y \in Y\}$ will be a definable isomorphism between $\mathcal{V}$ and $V^{\mathcal{V}\llbracket \omega \rrbracket}$ (or $V^{\mathcal{V}\llbracket Ord \rrbracket}$). In KMCU + Limitation of Size, there will be an injective map $F$ from $\A$ to $V$. For every urelement $a$, let $\tilde{a} = \<0, F(a)>$; and for every set $x$, we let $\tilde{x} = \<1, \{\tilde{y} : y \in x\}>$. It follows by an easy induction that the map $x \mapsto \tilde{x}$ and $X \mapsto \tilde{X} = \{\tilde{x} : x \in X\}$ is an isomorphism between $U$ and $V\llbracket F[\A] \rrbracket$.
\end{proof}

\begin{corollary}\label{corollary:LS<->AlessthanV}
Over KMCU, the following are equivalent.
\begin{enumerate}
    \item Limitation of Size.
    \item There is an injective map from $\A$ to $V$.
\end{enumerate}
\end{corollary}
\begin{proof}
(1) $\rightarrow$ (2) is clear. For (2) $\rightarrow$ (1), first observe that over KMCU, Limitation of Size holds for the pure classes $\mathcal{V}$. For, given a global choice function, its restriction to $V$ is a pure class so $\mathcal{V} \models $ Global Choice, which means $\mathcal{V} \models $ Limitation of Size by Proposition \ref{prop:GC<->GWO<->LS}. Now suppose that $F$ is an injective map from $\A$ to $V$. Since $U$ and $V\llbracket F[\A] \rrbracket$ are equinumerous, every proper class is equinumerous with a pure proper class. Therefore, all proper classes are equinumerous.
\end{proof}

\subsection{Independence results}\label{section:independenceinClassTheory}
In this subsection, I discuss several independence results concerning the completeness of Diagram \ref{GBUdiagram} over $\KMUR$. To begin with, arguments in \ref{subsection:implicationdiagram} which appeal to homogeneity can no longer go through in the context of class theory. For example, one might attempt to show that $\KMUR$ + Plenitude implies $\omega$-DC by the same argument as in Theorem \ref{Plenitude->DCS}. But the problem is that the ``kernel'' of a class relation $R$ might be a proper class, in which case we cannot find some a set of urelements that is big enough to ``fix'' $R$. In fact, we will see that over $\KMUR$,
\begin{enumerate}
    \item $\kappa$-DC $\nrightarrow$ Collection;
    \item Global Choice  $\nrightarrow$ ($\omega$-DC $\lor$ Collection);
    \item Collection $\nrightarrow$ $\omega$-DC;
    \item RP $\nrightarrow$ $\omega_1$-DC.
\end{enumerate}
Since it is well-known that KMc cannot prove Global Well-Ordering, it follows that Diagram \ref{GBUdiagram} is indeed complete over $\KMUR$.

Let me first discuss a general method of constructing \textit{class permutation models} of $\KMUR$ used in \cite{felgner1976choice}. Given a model $\<U, \A, \in, \U>$ of KMCU, we might fix an $\A$-ideal $\I$ as in Definition \ref{normalideal} and consider $U^\I$, the class of all first-order objects whose kernel is small in the sense of $\I$. This will give us a model of $\ZFCUR$ as before. To have a model of $\KMUR$, however, we cannot take all subclasses of $U^I$ as the second-order part of $U^\I$. To see this, suppose that $\A$ is an infinite set and $\I$ is its ideal of finite subsets; then an injection $F$ from $\omega$ to $\A$ will be a subclass of $U^\I$, which means having $F$ as a class of the model we intend to build would violate Replacement. In other words, we need to throw out some subclasses of $U^\I$. This is done by finding some suitable group $\G$ of permutations on $\A$ and only keep those subclasses of $U^\I$ that are \textit{symmetric} with respect to $\I$ and $\G$.

\begin{definition}[KMCU]
Let $\G$ be a group of permutations of $\A$.\footnote{It is understood that every $\pi \in \G$ is a permutation of a \textit{set} of urelements. Whenever $\pi, \sigma \in \G$, $\pi \circ \sigma$ is taken as the composition of their canonical extensions, which point-wise fix the urelements not in their original domains.} For any $x \in U$, $sym(x)$ and $fix(x)$ are defined as in Definition \ref{permutationmodeldef}. Let $\I\in \U$ be an $\A$-ideal as in Definition \ref{normalideal}. $\I$ is said to be $\G$-\textit{flexible} if 
\begin{enumerate}
    \item for every $\pi \in \G$ and $A \in \I$, $\pi A \in \I$;
    \item $\I$ has a \textit{basis} such that for every $B$ in the basis and $ A \in \I$ disjoint from $B$, there is some $\pi \in fix(B)$ such that $\pi A \neq A$.
\end{enumerate}
A class $X$ is \textit{symmetric} (w.r.t. $\G$ and $\I$) if there is some $A \in \I$, called a \textit{support} of $X$, such that $fix(A) \subseteq sym(X)$, where $sym(X) = \{\pi \in \G : \pi X = X\}$.  $U^\I =\{ x\in U : ker(x) \in \I\}$. $\W = \{ X \subseteq U^\I : X \text{ is symmetric}\}$. The model $\<U^\I, \in, \A, \W>$ is also denoted by $\W$.
\end{definition}
\noindent As we shall see later, $\G$-flexibility ensures that enough classes are thrown out so that Replacement can hold in $\W$. It is easy to check that every $\pi \in \G$ is an automorphism of $\W$ because $\I$ is $\G$-flexible.

\begin{theorem}[KMCU]\label{classpermutationmodel}
Let $\G$ be a group of permutations on $\A$ and $\I$ be an $\A$-ideal that is $\G$-flexible. Then $\W \models \KMUR$. 
\end{theorem}
\begin{proof}
$\W$ is a transitive class containing all the urelements and pure sets. And by the same argument as in Theorem \ref{smallkernelmodel}, it follows that $\W \models $ ZU + AC + Separation + Class Extensionality.

To show $\W \models$ Replacement, suppose that $x \in \W$ be a set and $F \in \W$ be a class function on $x$ with a support $A \in \I$. Let $\mathcal{J}$ be a basis of $\I$ which witnesses its flexibility and $B$ be a set of urelements in $\J$ with $ker(x) \cup A \subseteq B$. It suffices to show that $ker(F[x]) \subseteq B$. If not, then there is some $y \in F[x]$ such that $ker(y) \setminus B$ is not empty. Since $ker(y) \setminus B \in \I$, it follows that there is some $\pi \in fix(B)$ such that $\pi (ker(y) \setminus B) \neq ker(y) \setminus B$, and as a result, $\pi y \neq y$. But $F(z) = y$ for some $z \in x$ so $F(z) = \pi y$, which is a contradiction.

To show that Full Comprehension holds in $\W$, fix a formula $\varphi$ in the language of urelement class theory and classes $X_0, .., X_n \in \W$, in $\U$ let $X = \{ x \in  \W: \W \models \varphi (x, X_0, ... ,X_n \}$. It suffices to show that $sym(X_0) \cap ... \cap sym(X_n) \subseteq Sym(X)$. Fix a $\pi \in sym(X_0) \cap ... \cap sym(X_n)$. For every $x \in X$, since $\pi$ is an automorphism of $\W$, it follows that $\W \models \varphi (\pi x, X_0, ..., X_n)$ and hence $\pi x \in X$. Therefore, $\pi X = X$.
\end{proof}
\begin{theorem}\label{thm:KMUR+kDCvdashCollection}
Assume the consistency of KM. Let $\kappa$ be any infinite cardinal. There is a model of $\KMUR$ in which
\begin{enumerate}
    \item $\kappa$-DC holds;
    \item Collection fails.
\end{enumerate}
\end{theorem}
\begin{proof}
Let $\U$ be a model of KMCU + $\A \sim \aleph_{\kappa^+}$. Let $\I$ be the ideal of all sets of urelements of size less than $\aleph_{\kappa^+}$ and $\G$ be the group of all permutations of $\A$. It is clear that $\I$ is $\G$-flexible, so the resultant model $\W$ satisfies $\KMUR$. Collection fails in $\W$ because every cardinal below $\aleph_{\kappa^+}$ is realized while $\aleph_{\kappa^+}$ is not. To show $\kappa$-DC holds in $\W$, suppose that $R$ is a class relation in $\W$ without terminal nodes. Since the first-order domain of $\W$, $U^\I$, is closed under $\kappa$-sequences, $\U$ thinks that $\forall s \in (U^\I)^{<\kappa} \exists y \in U^\I R(x, y)$; by $Ord$-DC in $\U$, it follows that there is an $f \in (U^\I)^{\kappa}$ such that $R(f\restriction \alpha, f(\alpha))$ for every $\alpha < \kappa$. $f$ lives in $U^\I$, so $\kappa$-DC holds in $\W$.
\end{proof}

\begin{theorem}[Felgner \cite{felgner1976choice}]\label{thm:KMURnvdashCollection}
Assume the consistency of KM. There is a model of $\KMUR$ in which
\begin{enumerate}
    \item Global Choice holds;
    \item $\omega$-DC fails;
    \item Collection fails.
\end{enumerate}
\end{theorem}
\begin{proof}
Let $\U$ be a model of KMCU + $\A \sim \omega$, in which we identify $\A$ with the rationals $\<\Q, <_\Q>$. Let $\G$ be the group of permutations of $\A$ that preserves $<_\Q$ and $\I$ be the ideal of finite subsets of $\A$. $\I$ is $\G$-flexible: if $A, B \in \I$ are disjoint, for any $b \in B \setminus A$, there is an open interval containing $b$ that is disjoint from $A \cup B \setminus \{b\}$; so we can permute this interval in an order-preserving way and leave $A$ point-wise fixed. By Theorem \ref{classpermutationmodel}, it follows that the resultant class permutation model $\W$ satisfies $\KMUR$. Moreover, both Collection and $\omega$-DC fail in $\W$ since there is a proper class of urelements but every set of them is only finite. For a proof of $\W \models$ Global Choice, see \cite[pp. 249-250]{felgner1976choice} or \cite[Lemma 2.3]{yao2022reflection}. \end{proof}

\begin{theorem}\label{kmudoesnotproveDC}
Assume the consistency of KM. There is a model of $\KMUR$ in which 
\begin{enumerate}
    \item Collection holds;
    \item Plenitude holds;
    \item $\omega$-DC fails.
\end{enumerate}
\end{theorem}

\begin{proof}
The model used here is also due to Felgner\cite{felgner1976choice}. The point here is that the model also satisfies Collection, which is not discussed in Felgner's paper. Let $\U$ be a model of KMCU with an enumeration of $\A$ with the tree $Ord^{<\omega}\setminus \{\emptyset\}$ consisting of all non-empty finite sequences of ordinals. Each urelement is identified with a node on the tree and define $a \lhd b$ as $a \subsetneq b$. $b$ is said to be an \textit{immediate descendant} of $a$ if $b$ extends $a$ by one digit. $b$ and $b'$ are \textit{siblings} if either they are both top nodes, or they are an immediate descendant of the same node. A set $t \subseteq \A$ is a \textit{tree} if it is closed under initial segment (i.e, if $b \in t$ and $a \lhd b$, then $a \in t$). A \textit{path} of a tree $t$ with length $\alpha$ is a function $f: \alpha \rightarrow Ord$ such that $f\restriction \beta \in t$ for all $\beta < \alpha$. A \textit{branch} is a maximal path, i.e., it is not properly extended by any path of the tree. A tree $t$ is \textit{small} if it has no infinite branch. Let $\T$ be the class of all small trees, which forms a basis for an ideal $\I$. Let $\G$ be the group of  permutations of $\A$ that preserves $\lhd$.
\begin{lemma}\label{treeflexible}
$\I$ is $\G$-flexible with respect to the basis $\T$.
\end{lemma}
\begin{proof}
If $a$ and $b$ are siblings with domain $n + 1$, then there is a natural permutation $\pi^a_b \in \G$ such that for every node $c$ with $dom(c) = j +1$ and $i < j+1$,
\begin{equation*}
   \pi^a_b c (i)  =
    \begin{cases*}
      b(n) & if $i = n$, $c(n) = a(n)$ and $c \restriction n = a \restriction n$  \\
      a(n) & if $i = n$, $c(n) = b(n)$ and $c \restriction n = a \restriction n$  \\
      c(i) & otherwise
    \end{cases*}
 \end{equation*}
 $\pi^a_b$ thus swaps only $a$ and $b$ and their descendants. Given a small tree $t$ and some $A \in \I$ disjoint from $t$, since every node has  $Ord$-many siblings we can find a node $a$ in $A$ and a sibling $b$ of $a$ such that $b \notin t \cup A$. $\pi^a_b$ will then leave $t$ point-wise fixed because $t$ is a tree.
\end{proof}
\noindent Therefore, the class permutation model $\W$ given by $\G$ and $\I$ satisfies $\KMUR$. $\W$ clearly satisfies Plenitude because for every $\kappa$, there are $\kappa$-many top nodes on the tree $Ord^{<\omega}\setminus \{\emptyset\}$. Suppose \textit{for reductio} that $\omega$-DC holds in $\W$. Since in $\W$ every $a \in \A $ has some $b \in A$ such that $a \lhd b$, then in $\W$ there is an infinite branch $s = \<a_n : n < \omega>$ such that $a_n \lhd a_{n+1}$ for every $n$. Let $t$ be a small tree such that $fix(t) \subseteq sym(s)$. Fix some $a_n$ not in $t$ and some sibling $b$ of $a_n$ such that $b \notin t$. Then $\pi_{a_n, b} \in fix(t)$ but $\pi_{a_n, b} (s) \neq s$, which is a contradiction.

It remains to show that $\W \models$ Collection. For any two small tress $t$ and $t'$, we say that $t$ \textit{mildly extend} $t'$ if $t' \subseteq t$ and no branch of $t$ properly extends a branch of $t'$.

\begin{lemma}\label{treelemma}
Let $t_0, t_1 \in \T$ be such that $t_0 \subseteq t_1$ and every terminal node in $t_0$ has a descendent in $t_1$. Then for every $t \in \T$, there is a $\pi \in fix(t_0)$ such that $ t_1 \cup \pi t $ mildly extends $t_1$.  
\end{lemma}
\begin{proof}
Define $M = \{ a \in t_1 \setminus t_0 : a \text{ is an initial node in } t_1 \setminus t_0 \text{ and } a \text{ has a descendant in } t\}$, where ``$a$ initial in $t_1 \setminus t_0$'' means there is no node $b \in t_1 \setminus t_0$ such that $b \lhd a$. Since every node has $Ord$-many siblings and Global Well-Ordering holds in $\U$, for every $a \in M$ we can pick a sibling $a'$ of $a$ such that $a' \notin t_1 \cup t$, and we can ensure that $a_1' \neq a_2'$ for any distinct  $a_1, a_2 \in M$. Let $\pi = \bigcup_{a \in M}\pi^a_{a'}$, which is in $\G$. $\pi \in fix(t_0)$, because no node in $t_0$ is a descendant of any node in $M$ and $\pi$ only moves nodes in $M$ and their descendants.

To show that $t_1 \cup \pi t $ mildly extends $t_1$, consider any branch $f$ of $t_1$. Suppose \textit{for reductio} that $f$ is properly extended by a branch $g$ of $\pi t$. Note that $f$ must contain a node not in $t_0$ since otherwise $f$ would be a branch of $t_0$, which is impossible because every branch of $t_0$ is properly extended by a branch of $t_1$. So let $a$ be the least such node. There is a node $b$ on $g$ such that $a \lhd b$, where $b =\pi c$ for some $c \in t$. It follows that $a$ must be in $M$. If not, then $\pi a = a$ so $a \lhd c$ and hence $a$ is in $M$ after all. Thus, $\pi a = a'$ for some $a' \notin t$, but then $a' \lhd c$ so $a' \in t$---contradiction.
\end{proof}

\begin{lemma}
For any infinite cardinal $\kappa$, if $\<t_\alpha : \alpha < \kappa>$ is a sequence of small trees such that $t_\alpha$ mildly extends $\bigcup_{\beta < \alpha} t_\beta$ for every $\alpha < \kappa$, then $\bigcup_{\alpha < \kappa} t_\alpha$ is a small tree.
\end{lemma}
\begin{proof}
Let $t= \bigcup_{\alpha < \kappa} t_\alpha$. Suppose \textit{for reductio} that $f$ is an infinite branch of $t$. There will be some $t_\alpha$ with some $0 < m < \omega$ such that $f\restriction m$ is a branch of $t_\alpha$. Then for some $\beta > \alpha$, $t_\beta$ contains a branch that extends $f \restriction m$, which contradicts the assumption.
\end{proof}

Now suppose that $\W \models \forall x \in w \exists y \<x, y> \in R$ for some $w, R \in \W$. Let $t_0$ be a small tree that includes $ker(w)$ and some support of $R$, and enumerate $w$ by $\{x_\alpha : \alpha <\kappa\}$ for some $\kappa$. In $U$, we define a $\kappa$-sequence of small tress $\<t_\alpha : \alpha < \kappa>$ such that 
\begin{itemize}
    \item [] (i) $t_\alpha$ mildly extends $\bigcup_{\beta < \alpha} t_\beta$ for every $\alpha < \kappa$;
    \item [] (ii) for each $x_\alpha$, there is some $y \in \W$ such that $\<x_\alpha, y> \in R$ and $ker(y) \subseteq  t_\alpha$.
\end{itemize}
This is possible, because for every $x_\alpha$, fix some $y'$ with $\<x, y'> \in R$. Since $ker(y')$ is a subset of some small tree $t$, by Lemma \ref{treelemma}, there is a $\pi \in fix (t_0)$ such that $\bigcup_{\beta < \alpha} t_\beta \cup \pi t$ mildly extends $\bigcup_{\beta < \alpha} t_\beta$. Thus $\<x_\alpha, \pi y'> \in R$ and $ker(\pi y') \subseteq \pi t$. Let $t_\kappa = \bigcup_{\alpha < \kappa} t_\alpha$, which is a small tree. It follows that $\forall x \in w \exists y \in V(t_\kappa) \ (\<x, y> \in R)$, which suffices for Collection in $\W$.
\end{proof}
\noindent By using the same argument at the previous theorem, it is not hard to show that $\W \models \kappa$-CC for every infinite cardinal $\kappa$. However, it is unclear if the model $\W$ in Theorem \ref{kmudoesnotproveDC} satisfies CC if we assume $\U \models$ CC.

\begin{theorem}\label{thm:RPnvdashOmega1Dc}
Assume the consistency of KM. There is a model of $\KMUR$ in which 
\begin{enumerate}
    \item RP holds;
    \item Plenitude holds;
    \item $\omega_1$-DC fails.
\end{enumerate}
\end{theorem}
\begin{proof}
Let $\U$ be a model of KMCU with an enumeration of $\A$ with the tree $Ord^{<{\omega_1}}\setminus \{\emptyset\}$ consisting of all non-empty \textit{countable} sequences of ordinals. A set $t$ of urelements is an $\omega_1$-small tree if it is a tree without any $\omega_1$-branch. Let $\I$ be the ideal generated by all the $\omega_1$-small trees and $\G$ be the group of permutations of $\A$ that preserve the tree structure as in the proof of Theorem \ref{kmudoesnotproveDC}. Since the same arguments in Lemma \ref{treeflexible} and \ref{treelemma} still go through, it follows that $\I$ is $\G$-flexible and that Collection holds in the resultant model $\W$.
\begin{claim}
$\W \models $ RP. 
\end{claim}
\begin{claimproof}
By Proposition \ref{prop:rp+gbur<->collection+dcs}, it is enough to show that $\omega$-DC holds in $\W$. So it suffices to show that the ideal $\I$ is countably closed. This is simply because the union countably many $\omega_1$-small trees, $\{t_n : n< \omega \}$, is an $\omega_1$-small tree. 
\end{claimproof}
\begin{claim}
$\W \models $ $\neg$($\omega_1$-DC).
\end{claim}
\begin{claimproof}
Suppose \textit{for reductio} that $\omega_1$-DC holds in $\W$. Say that a sequence $s \in \A^\alpha$ is a \textit{chain} if $s(\beta) \lhd s(\beta')$ for every $\beta < \beta' < \alpha$; and $s$ is said to be \textit{bound} by $a$ if $s(\beta) \lhd a$ for every $\beta < a$. In $\W$, every chain $s\in \A^{<\omega_1}$ has a bound $a \in \A$. By  $\omega_1$-DC, in $\W$ there is an $f \in \A^{\omega_1}$ that is a chain of length $\omega_1$. Then $ker(f)$ must be contained in some $\omega_1$-small tree, which is impossible.
\end{claimproof}

\end{proof}

\subsection{Open questions}
It is a classic result that GBC is a conservative extension of ZFC (e.g., see \cite{felgner1976choice}). But we know that GBCU is \textit{not} a conservative extension of ZFCU: GBCU proves that either $\A$ is a set, or Plenititude holds, which is not provable in ZFCU.
\begin{question}
\
\begin{enumerate}
    \item Is $\GBUR$ + Global Choice conservative over $\ZFCUR$?
    \item Is $\GBUR$ + Collection + Global Choice conservative over ZFCU?
\end{enumerate}
\end{question} 
\noindent Given Theorem \ref{kmudoesnotproveDC} and the fact that CC is a stronger version of Collection, it is natural to ask the following.
\begin{question}
Does $\KMUR$ prove any of the following?
\begin{enumerate}
    \item Collection $\land$ Global Choice $\rightarrow$ $\omega$-DC.
    \item CC $\rightarrow$ $\omega$-DC.
    \item CC $\land$ Global Choice $\rightarrow$ $\omega$-DC.
\end{enumerate}
\end{question}
A natural question arises at this point as in Section \ref{section:WhatisZFCU}: what is KMc (or, GBc) class theory with urelements if we only wish to have AC for sets? Notably, in both class and set theory with urelements, Collection tends to lose its strength without enough choice. As previously conjectured, $\ZFUR$ + Collection does not prove RP, and in fact, $\ZFUR$ + DC should not be able to prove the DC$_\omega$-scheme. Theorem \ref{kmudoesnotproveDC} confirms that this is indeed the case in urelement class theory. Thus, although Collection is still strictly stronger than Replacement in class theory with urelements, adding it into $\KMUR$ cannot produce a theory with sufficient strength. So perhaps a robust version of KMc with urelements \textit{should} include RP as an axiom since it implies $\omega$-DC (Proposition \ref{prop:rp+gbur<->collection+dcs}). However, RP cannot exclude \textit{all} pathological models: in the proof of Theorem \ref{thm:RPnvdashOmega1Dc}, the model satisfies RP but contains a $Ord$-splitting tree without any $\omega_1$-branch. That said, it seems that bringing urelements back to the picture inevitably invites \textit{axiomatic freedom}.

\section{Second-order reflection with urelements}\label{section:RP2withurelements}

\subsection{Bi-interpretabtion with few urelements}
The second-order reflection principle (first introduced by Bernays \cite{bernays1976problem}) is the  scheme

\begin{itemize}
    \item [] (RP$_2$) $\forall X [\varphi(X) \rightarrow \exists t( t \text{ is transitive} \land \varphi^t(X \cap t))]$,
\end{itemize}
where $\varphi$ can be any formula in the language of class theory, and $\varphi^t$ is the result of restricting all the first-order quantifiers to the members of $t$ and all the second-order quantifiers to the subsets of $t$. Thus, $\varphi^t(X \cap t)$ is simply the assertion $\<t, \in, P(t)> \models \varphi(X \cap t)$. As observed in \cite{bernays1976problem} and \cite{tait2005constructing}, RP$_2$ is able to ``bootstrap''. For example, with the Axiom of Separation, Foundation and Extensionality, RP$_2$ can recover the remaining axioms of KMC.

\begin{prop}
RP$_2$ + Separation + Extensionality + AC + Foundation $\vdash$ KMCU + CC. 
\end{prop}
\begin{proof}
Note that RP alone implies that every $x_1, ..., x_n$ will be contained in some transitive set since we can reflect the formula $\exists y (x_1 = y) \land ... \land \exists y (x_n =y)$. This implies Pairing and Union given Separation. Then we can reflect the assertion ''for every $x$,  $x \cup \{x\}$ exists" down to a transitive set to get Infinity. Collection (and hence Replacement) follows by Proposition \ref{prop:rp+gbur<->collection+dcs}. To get Powerset, note that for every set $u$, by Separation we have ``for every class $X$ that is a subclass of $u$, there is a set $x$ that is co-extensional with $X$''. So by RP$_2$ we can reflect this assertion down to some transitive set $t$ containing $u$. Accordingly, $t$ contains every subset of $u$ as a member, which suffices for Powerset given Separation.

For Class Choice, suppose \textit{for reductio} that $\forall x \exists X \varphi(x, X, Z)$ for some class $Z$ but there is no $Y \subseteq U \times U$ such that $\varphi(x, Y_x, Z)$. By RP$_2$, there is some transitive set $t$ such that $\forall x \in t \exists X \subseteq t \varphi^t(x, X, Z \cap t)$ and there is no $Y \subseteq t\times t$ such that $\forall x \in t \varphi^t (x, Y_x, Z \cap t)$. Since there is a well-ordering of $P(t)$, for every $x \in t$ we can choose a $y_x \in P(t)$ such that $\varphi^t(x, y_x, Z\cap t)$. Let $Y = \bigcup_{x \in t}\{\<x, z> : z \in y_x\}$. It follows that $\forall x \in t \varphi^t (x, Y_x, Z \cap t)$, which is a contradiction.

Similarly, to show that there is a global well-ordering, we suppose \textit{for reductio} that there is no global well-ordering and reflect this statement to some transitive set $t$ that is closed under pairs. Since there is a well-ordering of $t$, the reflected statement will yield a contradiction.

Finally, we also get Full Comprehension. This is because if there is a failure of Full Comprehension of the form $\neg \exists X \forall z (z \in X \leftrightarrow \varphi(x, P))$, then we can reflect it down to some transitive $t$ to get $\neg \exists x \subseteq t \forall z (z \in x \leftrightarrow \varphi^t(x, P \cap t))$, which will then contradict Separation.\footnote{As a consequence, note that RP$_2$ is equivalent to the following scheme, which asserts that every statement, possibly with class parameters, is absolute to some transitive set.
\begin{itemize}
    \item [] (RP$_2^+$) $\forall X \exists \text{ transitive } t (\varphi(X) \leftrightarrow \varphi^t(X \cap t))$.
\end{itemize}
This is because given $\varphi$ and some class $X$, by Full Comprehension we can form the class $Y$ such that $\forall y (y \in Y \leftrightarrow \varphi(X))$; then by RP$_2$, there will be a non-empty transitive set $t$ such that $\forall y \in t (y\in Y\cap t \leftrightarrow \varphi^t(X \cap t))$. It follows that $\varphi(X) \leftrightarrow \varphi^t (X \cap t)$.}\end{proof}

In pure class theory, the bootstrapping of RP$_2$ goes beyond KMC as it yields large cardinals. In particular, RP$_2$ implies the exsitence of a proper class of inaccessible cardinals, Mahlo cardinals and weakly compact cardinals (see \cite{tait2005constructing} for more on this). However, the consistency strength of RP$_2$ is bounded by ZFC + an $\omega$-Erd{\"o}s cardinal (see \cite[Exercise 9.18]{kanamori2008higher}), which is consistent with $V = L$. Some natural question arise in the context of urelements: What is the consistency strength of RP$_2$ in urelement class theory? Could it be somehow affected by urelements? 

The next lemma shows that the $\mathcal{V} \llbracket X \rrbracket$ construction introduced in Definition \ref{barwiseinterpretation2} preserves second-order reflection.
\begin{lemma}\label{lemma:KM+RP2interpretsKMU+RP2}
Let $\mathcal{V} \models$ KM + RP$_2$ and $W \in \mathcal{V}$ be a class. Then $\mathcal{V} \llbracket W \rrbracket \models $ KMU + RP$_2$ + Limitation of Size.
\end{lemma}
\begin{proof}
Since RP$_2$ + KM proves that there is a global well-ordering, which implies Limitation of Size over KM, by Theorem \ref{Con(KM)->Con(KMU)} it follows that $\mathcal{V} \llbracket W \rrbracket \models $ KMU + Limitation of Size. So it remains to show that every instance of RP$_2$ holds in $\mathcal{V} \llbracket W \rrbracket$. For every transitive set $t \in V$, let $t\llbracket W \rrbracket = \<1, V\llbracket W \rrbracket \cap t>$, which is a transitive set in $V\llbracket W \rrbracket$.

\begin{claim}\label{permu}
Let $t$ be a transitive set in $\mathcal{V}$. Then for any $x_1, ... x_n \in t\llbracket W \rrbracket$ and $X_1, ... ,X_m \in \mathcal{V} \llbracket W \rrbracket$, $ \mathcal{V} \models (\varphi^{\mathcal{V} \llbracket W \rrbracket})^t \leftrightarrow (\varphi^{t\llbracket W \rrbracket})^{\mathcal{V} \llbracket W \rrbracket}$ for any suitable formula $\varphi$ in the language of urelement class theory.
\end{claim}
\begin{claimproof}
If $\varphi$ is an atomic formula, then the claim holds because the definition of $\bar{\in}$ and $\bar{\A}$ (see Definition \ref{barwiseinterpretation1}) is absolute for transitive sets. Boolean cases commute. And if $\varphi$ is $\exists x \psi$, we have
\begin{align*}
    (\varphi^{\mathcal{V} \llbracket W \rrbracket})^t  &= (\exists x \in \mathcal{V} \llbracket W \rrbracket  \psi^{\mathcal{V} \llbracket W \rrbracket})^t \\
                &\Leftrightarrow \exists x \bar{\in} t\llbracket W \rrbracket (\psi^{\mathcal{V} \llbracket W \rrbracket})^t \\
                & \Leftrightarrow \exists x \bar{\in} t\llbracket W \rrbracket (\psi^{t\llbracket W \rrbracket})^{\mathcal{V} \llbracket W \rrbracket} \tag*{(by induction hypothesis)} \\
                & = (\varphi^{t\llbracket W \rrbracket})^{\mathcal{V} \llbracket W \rrbracket}.
\end{align*}
Similarly, if $\varphi$ is $\exists X \psi$, then we have
\begin{align*}
    (\varphi^{\mathcal{V} \llbracket W \rrbracket})^t  &= (\exists X \subseteq V \llbracket W \rrbracket  \psi^{\mathcal{V} \llbracket W \rrbracket})^t \\
                &\Leftrightarrow \exists X \subseteq t\llbracket W \rrbracket (\psi^{\mathcal{V} \llbracket W \rrbracket})^t \\
                & \Leftrightarrow \exists  X \subseteq t\llbracket W \rrbracket (\psi^{t\llbracket W \rrbracket})^{\mathcal{V} \llbracket W \rrbracket} \tag*{(by induction hypothesis)}\\
                & = (\varphi^{t\llbracket W \rrbracket})^{\mathcal{V} \llbracket W \rrbracket}.
\end{align*}
This proves the claim.                                                  \end{claimproof}

Now if $\mathcal{V} \llbracket W \rrbracket \models \varphi$, then we can reflect $\varphi^{\mathcal{V} \llbracket W \rrbracket}$ in $\mathcal{V}$ down to some transitive set $t$. By the claim, it follows that $\mathcal{V} \llbracket W \rrbracket \models \varphi ^{t\llbracket W \rrbracket}$. Hence, $\mathcal{V} \llbracket W \rrbracket \models$ RP$_2$. \end{proof}
\begin{theorem}\label{KM + RP <-> KMU + RP + LS}
KM + RP$_2$ and KMU + RP$_2$ + Limitation of Size are bi-interpretable with parameters. \qed
\end{theorem}
\begin{proof}
First note that KMU + RP$_2$ also interprets KM + RP$_2$, because if $\mathcal{U} \models$ KMU + RP$_2$, then its pure part $\mathcal{V} \models$ KM + RP$_2$. For, as in Lemma \ref{lemma:KM+RP2interpretsKMU+RP2}, given a transitive $t \in U$, we can show that $(\varphi^{\mathcal{V}})^t \leftrightarrow (\varphi^{V\cap t})^\mathcal{V}$, where $V\cap t$ is a transitive pure set. So if $\mathcal{V} \models \varphi$, then we can reflect $\varphi^\mathcal{V}$ to some transitive set $t$, which implies $\mathcal{V} \models \varphi^{V\cap t}$. And given Lemma \ref{lemma:KM+RP2interpretsKMU+RP2}, it follows that KM + RP$_2$ and KMU + RP$_2$ + Limitation of Size are mutually interpretable. And their bi-interpretability (with parameters) follows from Theorem \ref{KMU+LSbiintepsKM+LS}.
\end{proof}
 As a consequence, KMU + RP$_2$ also implies the existence of a proper class of inaccessible cardinals, Mahlo cardinals and weakly compact cardinals. By Corollary \ref{corollary:LS<->AlessthanV}, it follows that when the urelements are \textit{few}, i.e., no more numerous than the pure sets, RP$_2$ has the same strength as in pure class theory.

\subsection{A model of RP$_2$ with many urelements}\label{section:amodelofKMU+RP2+notLS}
What if there are more urelements than the pure sets? Or, does KMU + RP$_2$ prove Limitation of Size? In this final section, I construct a model of KMU + RP$_2$ where the urelements are more numerous than the pure sets by assuming the consistency of a $\kappa^+$-supercompact cardinal. To begin with, there is an alternative accumulative hierarchy that can produce natural models of KMCU where Limitation of Size fails.
\begin{definition}
Let  $\kappa$ be an infinite cardinal. For any set $x$, $P_\kappa(x)$ is the set of all subsets of $x$ of size less than $\kappa$. For any set of urelements $A$, $U_{\kappa, A} = \bigcup_{B \in P_\kappa (A)} V_\kappa (B)$. $\mathcal{U}_{\kappa, A}$ denotes the model $\langle U_{\kappa, A}, A, \in, P(U_{\kappa, A})\rangle$ .
\end{definition}
\noindent The $U_{\kappa, A}$-hierarchy is a generalization of the $V_\kappa(A)$-hierarchy: $U_{\kappa, A} = V_\kappa(A)$ when the size of $A$ is no greater than $\kappa$. While every $A$ appeas as a set in $V_\kappa(A)$, $A$ would be a proper class in $U_{\kappa, A}$ when its size is greater than $\kappa$.  Another useful stratification is the $H_\kappa(A)$-hierarchy (used in \cite{HamkinsForthcoming-HAMRIS}), where $H_\kappa(A) = \{x \in U : ker(x) \subseteq A \land |trc(\{x\})| < \kappa \}$. Note that when $\kappa$ is inaccessible and $|A| > \kappa$, $H_\kappa(A) = U_{\kappa, A}$.
\begin{lemma}[ZFCU]\label{UKA}
For any transitive set $t$, the following are equivalent.
\begin{enumerate}
    \item $t = U_{\kappa, A}$, where $\kappa$ is inaccessible and $A \subseteq \A$ .
    \item $\<t, ker(t), \in, P(t)> \models$ $\KMUR$.
\end{enumerate}
In fact, $\mathcal{U}_{\kappa, A} \models $ KMCU + CC whenever $\kappa$ is inaccessible; and Limitation of Size fails in $\mathcal{U}_{\kappa, A}$ when $A$ has size greater than $\kappa$.
\end{lemma}
\begin{proof}
(1) $\rightarrow$ (2). $\mathcal{U}_{\kappa, A}\models$ ZU + Class Extensionality since it is transitive and sufficiently tall. For example, to show $\mathcal{U}_{\kappa, A}\models$ Powerset, fix some $x \in V_\alpha(B)$ for some $B \in P_\kappa(A)$ and $\alpha < \kappa$. Then $|V_\alpha(B)| < \kappa$ since $\kappa$ is a strong limit; so $P(x)$ is a subset of $V_\kappa(B)$ of size less than $\kappa$, and it will be contained in $V_\beta(B)$ for some $\beta <\kappa$ as $\kappa$ is regular. $\mathcal{U}_{\kappa, A}\models$ Global Well-Ordering since the well-ordering of  $U_{\kappa, A}$ in $U$ is a class of $U_{\kappa, A}$. To show $\mathcal{U}_{\kappa, A}\models$ Class Choice, suppose that $\mathcal{U}_{\kappa, A}\models \forall i \in I \exists X \varphi(i, X)$ for some $I \subseteq U_{\kappa, A}$. In $U$, we can well order $P(U_{\kappa, A})$ and then for each $i \in I$, choose some $X_i \in P(U_{\kappa, A})$ such that $\mathcal{U}_{\kappa, A}\models  \varphi(i, X_i)$. $Y = \bigcup_{i \in I}\{\<i, x> : x \in X_i\}$ will then be a desired class of $U_{\kappa, A}$. Note that Class Choice implies Collection and hence Replacement, so $\mathcal{U}_{\kappa, A}\models$ KMCU. And clearly, when $A$ has size greater than $\kappa$, $A$ and $\kappa$ are two proper classes in $\mathcal{U}_{\kappa, A}$ that are not equinumerous.

(2) $\rightarrow$ (1). Suppose that $ \mathcal{T} = \<t, ker(t), \in P(t)> \models$ $\KMUR$. Let $\kappa = Ord \cap t$. $\kappa$ must be a regular cardinal. Suppose not. Then given a cofinal sequence $f$ on $\kappa$ with length $\alpha$, where $\alpha < \kappa$, since $t$ is closed under ordered-pairs $f$ is a class function in $\mathcal{T}$ on $\alpha$. By Replacement in $\mathcal{T}$, it follows that $\kappa \in t$, which is a contradiction. To show $\kappa$ is a strong limit. First note that for every set $x \in t$, $P(x) = P^\mathcal{T}(x)$. This is because for every set $y \subseteq x$, $y$ is a class in $\mathcal{T}$ and so by Separation in $\mathcal{T}$, $y \cap x = y \in t$. So if $\alpha < \kappa$, then $P(\alpha)$ is in $t$ and by AC in $\mathcal{T}$, it is equinumerous with some $\beta < \kappa$. Clearly, $\omega < \kappa$, so $\kappa$ is inaccessible.

Now let $A = ker(t)$. It remains to show that $t = U_{\kappa, A}$. First note that $P_\kappa(t) \subseteq t$. For, any enumeration of $x$ with some ordinal $\alpha < \kappa$ is a class in $\mathcal{T}$, so $x \in t$ by Replacement in $\mathcal{T}$. If $B \subseteq A$ is of size less than $\kappa$, then by an easy induction $V_\alpha(B)$ has size less than $\kappa$ for all $\alpha < \kappa$. This shows that $U_{\kappa, A} \subseteq t$. For every set $x \in t$, let $B = ker(x)$. Since $\mathcal{T} \models \KMUR$, $B \in t$. Then $B$ must have size less than $\kappa$ because $\mathcal{T} \models B \sim \alpha$ for some $\alpha < \kappa$. Let $\beta$ be the least ordinal such that $x \in V_\beta(B)$. As $x \in t$, it is clear that $\beta < \kappa$. Therefore, $x \in U_{\kappa, A}$. This shows that $t = U_{\kappa, A}$.                                      
\end{proof}

 Recall Zermelo's Quasi-Categoricity Theorem: any full second-order model of second-order ZF is isomorphic to some $V_\kappa$, where $\kappa$ is inaccessible. Now let ZFCU$_2$ be the corresponding version of ZFCU formulated in the second-order language. We then have the following generalized quasi-categoricity theorem in ZFCU + Plenitude.

\begin{theorem}[ZFCU + Plenitude]
For every full second-order structure $\M$, $\M \models$ ZFCU$_2$ if and only if $\M$ is isomorphic to some $\mathcal{U}_{\kappa, A}$, where $A \subseteq \A$ and $\kappa$ is an inaccesible cardinal.
\end{theorem}
\begin{proof}
Since $\M \models$ ZFCU$_2$ and it is a full second-order model, a standard argument shows that $\in^\M$ is well-founded. By AC and Plenitude, we can then fix a bijection $i$ from $\A^\M$, the class of urelements in $\M$, to a set of urelements $A$. $i$ can then be extended to $\M$ by letting $i (x) = \{i(y) : y \ \in^\M x\}$ as in Mostowski collapse. $i[\M]$ is then a transitive set $t$ such that $\<t, ker(t), \in, P(t)> \models \KMUR$. So $t = U_{\kappa, A}$ for some $A \subseteq \A$ and inaccessible cardinal $\kappa$ by Lemma \ref{UKA}.
\end{proof}
Now I proceed to prove the following.
\begin{theorem}\label{thm:RP2nvdashLS}
Assume the consistency of  $\text{ZFC} + \exists \kappa (\kappa \text{ is } \kappa^+ \text{-supercompact})$. There is a model of KMCU in which
\begin{enumerate}
    \item RP$_2$ holds;
    \item Limitation of Size fails.
\end{enumerate}
\end{theorem}
\begin{proof}
Let $V \models \text{ZFC} + \exists \kappa (\kappa \text{ is } \kappa^+ \text{-supercompact})$, where $\kappa^+$-supercompactness is defined as having a normal fine measure on $P_{\kappa}(\kappa^+)$. Note that by class forcing we can add a global well-ordering to $V$ without adding any new sets, which yields a model $\mathcal{V} \models$ GBC + Limitation of Size +  $\exists \kappa (\kappa \text{ is } \kappa^+ \text{-supercompact})$. By Theorem \ref{Con(KM)->Con(KMU)}, this gives us a model $\U \models$ GBCU + Limitation of Size + Plenitude + $\exists \kappa (\kappa \text{ is } \kappa^+ \text{-supercompact})$ (e.g., consider $\<V\llbracket Ord\rrbracket, \{0\}\times Ord, \in, \mathcal{V}\llbracket Ord \rrbracket >$). Working in $\mathcal{U}$, let $F$ be a normal fine measure on $P_{\kappa}(\kappa^+)$. For any functions $f$ and $g$ on $P_{\kappa}(\kappa^+)$,  define the equivalence relation
$$f =_F g \text{ if and only if } \{ x \in P_{\kappa}(\kappa^+) : f(x) = g(x) \} \in F.$$ 
Global Well-Ordering then allows us to pick a unique $f$ from each equivalence class $[g]_{=_F}$ and then form an internal ultrapower $U/F$ as in the beginning of Section \ref{section:WhatisZFCU}, which is a class in $\U$. Note that since the first-order part $U$ satisfies ZFCU (in particular, Collection), \L{}o\'s's Theorem holds for $U/F$ (see Theorem \ref{thm:collection<->losthm}). That is, for every $f_1, ... f_n \in U/F$, $ U/F \models \varphi (f_1..., f_n)$ if and only if $\{x \in P_{\kappa}(\kappa^+): \varphi (f_1(x)...,f_n(x))\} \in F$. 

\begin{lemma}
$\in_F$ is a well-founded and set-like relation on $U/F$.
\end{lemma}
\begin{proof}
$\in_F$ is well-founded because $F$ is $\kappa$-complete. To show it is set-like, fix any $f \in U/F$ and let $X = \{g \in U/F : g \in_F f \}$. We may assume that the set $ \bar{y}= \{x \in P_{\kappa}(\kappa^+) : f(x) \neq \emptyset \}$ is in $F$. Let $\bar{z}$ be the set of all functions from $\bar{y}$ to $(\bigcup f[\bary]) \cup \{\emptyset\}$. For each $g \in_F f$, we define a function $g'$ on $P_{\kappa}(\kappa^+)$ as follows.
\begin{equation*}
    g' (x) =
    \begin{cases*}
      g(x) & if $g(x) \in f(x)$ \\
     \emptyset        & otherwise 
    \end{cases*}
  \end{equation*}
$g' \restriction \bary$ is in $\barz$ for every $g$ such that $g \in_F f$. It suffices to show that the map $g \mapsto g' \restriction \bary$ is 1-1 from $X$ into $\barz$. Consider two $g_1, g_2 \in X$ .  Since $\bary \cap \{x \in P_{\kappa}(\kappa^+) : g_1(x) \neq g_2(x) \land g_1(x) \in f(x) \land g_2(x) \in f(x)\}$ is in $F$, there must be some $x \in \bary$ such that $g_1'(x) \neq g_2'(x)$ and hence $g_1' \restriction \bary \neq g_2' \restriction \bary$.
\end{proof}

Now we wish to collapse $U/F$ into a transitive class $M$, which yields an elementary embedding from $U$ to $M$. For reasons that will be clear, it is useful to have the elementary embedding fix $\kappa^+$-many urelements.\footnote{Note that we cannot expect the resulting elementary embedding $j$ to fix all the urelements. Otherwise, let $A$ be a set of urelements of size $\kappa$; then $j(A) = A$ but $|j(A)| = j(\kappa) > \kappa $.} So let $A$ be a set of urelements in $U$ enumerated by $\langle a_\alpha : \alpha < \kappa^+\rangle$. For every $y \in U$, let $C_y \in U/F$ be the function that is $=_F$-equivalent to the constant function that maps everything to $y$. Since all proper classes are equinumerous in $\mathcal{U}$, there is a one-one mapping $G$ from $\A_F \setminus \{C_{a_\alpha}: \alpha < \kappa^+ \}$ into $\A \setminus A$, where $\A_F$ is the class of urelements in $U/F$. We then define the collapsing function $\pi$ as follows. For every $f \in \A_F$, 
\begin{equation*}
    \pi (f) =
    \begin{cases*}
      a_\alpha & if $f = C_{a_\alpha}$, for some $a_\alpha \in A$ \\
      G(f)        & otherwise 
    \end{cases*}
  \end{equation*}
And for $f \in U/F \setminus \A_F$, we let $\pi(f) = \{\pi(g): g \in_F f \}$, which is well-defined by the previous lemma. 

\begin{definition}
Let $M = \pi[U/F]$, $i : U \rightarrow U/F$ be such that $i(y) = C_y$, and $j = \pi \circ i$.
\end{definition}
By \L{}o\'s's Theorem, $j$ is an elementary embedding from $U$ to $M$. Note that $j$ fixes every urelement in $A$ because for every $a_\alpha \in A$, $j(a_\alpha) = \pi (C_{a_\alpha}) = a_\alpha$. Therefore, $A \subseteq j(A)$.

\begin{lemma} Let $\kappa, M, j$ be defined as above.
\begin{itemize}
    \item [] (i) $j(\gamma) = \gamma $ for all $\gamma < \kappa$;
    \item [] (ii) $j(\kappa) > \kappa^+$;
    \item [] (iii) $M^{\kappa^+} \subseteq M$.
\end{itemize}
\end{lemma}
\begin{proof}
All by standard text-book arguments.
\end{proof}
\noindent In particluar, $A \in M$. By Lemma \ref{UKA}, $\mathcal{U}_{\kappa, A}$ is a model of KMCU where Limitation of Size fails. It remains to show $\mathcal{U}_{\kappa, A} \models$ RP$_2$.

\begin{lemma}\label{j}
\singlespacing
For every $x \in U_{\kappa, A}$ and  $y \subseteq  U_{\kappa, A}$, $j(x) =x$ and $y = j(y) \cap  U_{\kappa, A}$.
\end{lemma}
\begin{proof}
First observe that for every set $x$ with $|x| < \kappa$, $j(x) = j[x] = \{j(y) : y \in x \}$.
Let $f:\alpha \rightarrow  x$ be a surjection, where $\alpha < \kappa$. $j(f)$ is then a surjection from $\alpha$ onto $j(x)$. It suffices to show that that $jf[\alpha] = j[x]$. If $y \in x$, then $y= f(\beta)$ for some $\beta < \alpha$ so $j(y) = jf(j(\beta))=jf(\beta) \in jf[\alpha]$. On the other hand, for $\beta < \alpha$, $jf(\beta) = jf(j(\beta)) = jf(\beta) = j(f(\beta)) \in j[x]$. Now given any $x \in U_{\kappa, A}$, since $|x| < \kappa$ and $j$ fixes all the urelements in $A$ it follows that $j(x)=x$.
\end{proof}

\begin{lemma}\label{M}
$U_{\kappa, A} = U_{\kappa, A}^M$, and $P (U_{\kappa, A}) = P (U_{\kappa, A})^M$.
\end{lemma}
\begin{proof}
Since $M$ is transitive and closed under $\kappa^+$-sequences, for every $x \in M$, $M \models |x| < \kappa$ if and only if $|x| < \kappa$. This shows that $P_{\kappa}(A)^M = M \cap P_{\kappa}(A)$ so $U_{\kappa, A}^M = U_{\kappa, A} \cap M$. But $U_{\kappa, A} \subseteq M$ by Lemma \ref{j}; thus, $U_{\kappa, A} = U_{\kappa, A}^M $. If $y \subseteq U_{\kappa, A}$, by Lemma \ref{j} $y$ is in M. Therefore, $P (U_{\kappa, A}) = P (U_{\kappa, A})^M$.
\end{proof}

\begin{lemma}
$\mathcal{U}_{\kappa, A}\models$ RP$_2$.
\end{lemma}

\begin{proof}
Suppose that $ \mathcal{U}_{\kappa, A}\models \varphi(x, Y)$, where $x \in U_{\kappa, A}$ and  $Y \subseteq U_{\kappa, A}$. By Lemma \ref{j}, we have
\begin{align}
\varphi(j(x), j(Y) \cap U_{\kappa, A})^{\mathcal{U}_{\kappa, A}}
\end{align}
It then follows from Lemma \ref{M} that 
\begin{align}
M \models \varphi(j(x), j(Y) \cap U_{\kappa, A} )^{\mathcal{U}_{\kappa, A}}
\end{align}
Since $A \subseteq j(A)$ and $M \models |A| < j(\kappa)$, it follows that
\begin{align}
M \models \exists \lambda < j(\kappa) \exists B \subseteq j(A) [|B| = \lambda \land \varphi(j(x), j(Y) \cap U_{\lambda, B})^{\mathcal{U}_{\lambda, B}}]
\end{align}
By the elementarity of $j$, we have
\begin{align}
\exists \lambda < \kappa \exists B \subseteq A [|B| = \lambda \land \varphi(x,  Y\cap U_{\lambda, B})^{\mathcal{U}_{\lambda, B}}]
\end{align}
Fix such $\lambda$ and $B$. $U_{\lambda, B} = \bigcup \{V_\lambda (C) : C \in P_{\lambda}(B)\}$ is a subset of $V_{\kappa} (B)$ with size less than $\kappa$, so $U_{\lambda, B} \in V_{\kappa} (B)$ and hence $U_{\lambda, B} \in U_{\kappa, A}$. Therefore,
\begin{align}
\mathcal{U}_{\kappa, A} \models \exists t [ t \text{ is transitve} \land \varphi(x, Y\cap t)^t].
\end{align}
\end{proof}
\noindent This completes the proof of Theorem \ref{thm:RP2nvdashLS}. 
\end{proof}
The result here is extended and improved in \cite{HamkinsForthcoming-HAMRIS}, where KMU + RP$_2$ + more than $Ord$-many urelements is shown to be consistent relative to ZFC with a \textit{nearly} $\kappa^+$-supercompact cardinal $\kappa$. The notion of nearly $\lambda$-supercompact cardinals is first studied by Schanker in \cite{Schanker2011:Dissertation}, \cite{Schanker2011:WeaklyMeasurableCardinals} and \cite{Schanker2013:PartialNearSupercompactness}. Although a nearly $\kappa^+$-supercompact cardinal is strictly weaker than a $\kappa^+$-supercompact cardinal, it remains a strong large cardinal axiom as shown in \cite{Schanker2011:WeaklyMeasurableCardinals}. Moroever, in \cite{HamkinsForthcoming-HAMRIS} we prove that if there are \textit{abundant} urelements in some second-order sense, then KMU + RP$_2$ is bi-interpretable with KMC plus a supercompact cardinal. These results together reveal an interesting interaction between \textit{limitation of size} and \textit{reflection}, the two philosophical conceptions of set mentioned in Section \ref{section:UrelementsinSetTheory}. Limitation of size is often viewed as a \textit{maxiamality principle} (see G\"odel \cite{godel1986collected}) because it asserts that any collection of objects that is not ``too big'' can form a set. However, with urelements, this view is challenged since limitation of size is precisely the reason why reflection has little strength. On the one hand, Theorem \ref{KM + RP <-> KMU + RP + LS} shows that under limitation of size, second-order reflection is still a weak large cardinal axiom. On the other hand, according to Theorem \ref{thm:RP2nvdashLS} and the results in \cite{HamkinsForthcoming-HAMRIS}, a strong \textit{violation} of Limitation of Size can dramatically increase the strength of second-order reflection. Thus, under the reflection conception it is the \textit{violation} of limitation of size that maximizes.


\appendix
\include{Chapters/appendix}


\backmatter              

\bibliographystyle{nddiss2e}
\bibliography{masterbib}           

\end{document}